\RequirePackage{fix-cm}
\documentclass[smallcondensed]{svjour3}     % onecolumn (ditto)
\smartqed  % flush right qed marks, e.g. at end of proof
\usepackage{graphicx}
\usepackage{bbm}
\usepackage{bm}
\usepackage{color}
\usepackage{xcolor}
\usepackage{subfig}
\usepackage{epstopdf}
\usepackage{psfrag}
\usepackage{amssymb}
\usepackage{amsmath}
\usepackage{dsfont}
\usepackage{rotating}
\usepackage{pgf,tikz}

\newcommand{\Ss}{ \mathbf{\mathcal{S}}}
\newcommand{\Nn}{ \mathbf{\mathcal{N}}}
\newcommand{\Ff}{ \mathbf{\mathcal{F}}}
\newcommand{\D}{\mathbf{\mathcal{D}}}
\newcommand{\RM}{\mathbf{\mathcal{R}}}
\newcommand{\LM}{\mathbf{\mathcal{L}}}

\newcommand{\Mpsi}{\vec{\mathcal{M}}}
\newcommand{\Q}{\vec{\mathcal{Q}}}
\newcommand{\F}{\mathbf{F}}

\newcommand{\bphi}{ \boldsymbol{\phi}}
\newcommand{\bpsi}{ \boldsymbol{\psi}}
\newcommand{\etah}{ \hat{\mathbf{p}}}
\newcommand{\vvh}{\hat{\mathbf{v}}}
\newcommand{\diff}[2]{\frac{\partial {#1} }{\partial {#2} } }

\spnewtheorem{rem}{Remark}[section]{\bfseries}{}
\journalname{Communications on Applied Mathematics and Computation}
\begin{document}
\title{A novel staggered semi--implicit space--time discontinuous Galerkin method for the incompressible Navier-Stokes equations}
\titlerunning{Novel staggered semi-implicit space-time discontinuous Galerkin schemes}        % if too long for running head
\author{F. L.  Romeo  \and M. Dumbser  \and
	M. Tavelli  %etc.
}
\institute{Francesco Lohengrin Romeo \at
            Department of Civil, Environmental and Mechanical Engineering, University of Trento, Via Mesiano 77, I-38123 Trento, Italy \\
            \email{francesco.romeo-4@unitn.it}           %  \\
 			\and     
           	Michael Dumbser \at
			Department of Civil, Environmental and Mechanical Engineering, University of Trento, Via Mesiano 77, I-38123 Trento, Italy \\
			\email{michael.dumbser@unitn.it}
			\and
			Maurizio Tavelli \at
			Department of Civil, Environmental and Mechanical Engineering, University of Trento, Via Mesiano 77, I-38123 Trento, Italy \\
			\email{m.tavelli@unitn.it}       
}
\date{Received: 31 December 2019}%/ Accepted: date}
\maketitle
\begin{abstract}
A new high order accurate staggered semi-implicit space--time discontinuous Galerkin (DG) method is presented for the simulation of viscous incompressible flows on unstructured triangular grids in two space dimensions. The staggered DG scheme defines the discrete pressure on the primal triangular mesh, while the discrete velocity is defined on a staggered edge-based dual quadrilateral mesh. In this paper, a new pair of  equal-order-interpolation velocity-pressure finite elements is proposed. On the primary triangular mesh (the pressure elements) the basis functions are piecewise polynomials of degree $N$ and are allowed to jump on the boundaries of each triangle. 
On the dual mesh instead (the velocity elements), the basis functions consist in the union of piecewise polynomials of degree $N$ on the two subtriangles that compose each quadrilateral and are allowed to jump 
only on the dual element boundaries, while they are continuous inside. In other words, the basis functions on the dual mesh are built by continuous finite elements on the subtriangles. 
This choice allows the construction of an efficient, quadrature-free and memory saving algorithm. In our coupled space-time pressure correction formulation for the incompressible Navier-Stokes equations, arbitrary high order of accuracy in time is achieved through the use of time-dependent test and basis functions, in combination with simple and efficient Picard iterations. Several numerical tests on classical benchmarks confirm that the proposed method outperforms existing staggered semi-implicit space-time DG schemes, not only from a computer memory point of view, but also concerning the computational time. 
\end{abstract}
\section{Introduction}
\label{intro}
%Our papers
In this article we propose a new pair of velocity-pressure elements in the framework of the family of arbitrary high-order accurate staggered semi-implicit space--time Discontinuous Galerkin (DG) methods for the two-dimensional incompressible Navier-Stokes equations on unstructured meshes, extending the works in \cite{STINS2D,STINS2DTri} by Tavelli and Dumbser.\\
%Navier-Stokes
The discretization of the incompressible Navier-Stokes equations was mainly carried out in the past using finite difference methods \cite{markerandcell,patankar,patankarspalding,vanKan} or continuous finite element schemes \cite{SUPG,Fortin,Rannacher1,Rannacher3,SUPG2,TaylorHood,Verfuerth}.
%DG 
While until a few years ago finite volume methods were more popular, nowadays DG schemes, first introduced by Reed and Hill in \cite{reed}, are now widely applied in several  different fields. In a collection of seminal works \cite{cbs3,cbs2,cbs1,cbs0,cbs4}, Cockburn and Shu provided a rigorous formulation of these methods, contributing to their widespread use nowadays. For their applications to the Navier-Stokes equations, let us mention the pioneering works of Bassi and Rebay \cite{BassiRebay}, and Baumann and Oden \cite{Baumann199979}. Later, several high order DG methods for the incompressible and compressible Navier-Stokes equations have been proposed, see for example \cite{Bassi2006,Bassi2007,Crivellini2013,DumbserNSE,Ferrer2011,MunzDiffusionFlux,HartmannHouston1,HartmannHouston2,KleinKummerOberlack2013,Nguyen2011,Rhebergen2012,Rhebergen2013,Shahbazi2007}.\\
%Semi-implicit
The typical saddle point problem of the incompressible Navier-Stokes equations is conveniently managed through the well-known \textit{semi-implicit} technique on staggered meshes introduced by Harlow and Welch in \cite{markerandcell}. Further related semi-implicit methods were developed for the incompressible Navier-Stokes equations in \cite{Bell1989,chorin1,chorin2,HirtNichols,patankar,vanKan}. For asymptotic preserving semi-implicit schemes, also in the context of all Mach number flow solvers, see e.g.  \cite{KleinMach,Munz2003,ParkMPV,CordierDegond,DumbserCasulli2016,RussoSI,RussoSI2,RussoAllMach}. 

A whole family of staggered finite difference methods for free-surface flows has been developed in the last decades by Casulli et al., see  \cite{Casulli2009,CasulliVOF,casulli:1984,CasulliStelling2011,CasulliWalters2000}. A theoretical analysis of this approach can be found in \cite{BrugnanoCasulli,BrugnanoCasulli2,BrugnanoSestini,CasulliCattani,CasulliCheng1992,CasulliZanolli2012}.  
All of the mentioned works are characterized by the property that  the nonlinear convective terms are discretized explicitly, while the pressure terms are discretized implicitly. Thanks to their high computational efficiency, these methods have been extended to many different applications in science and engineering: let us mention for example blood flow in compliant vessels \cite{CasulliDumbserToro2012,FambriDumbserCasulli} or compressible gas dynamics in elastic tubes \cite{DumbserIbenIoriatti}. 
The first direct extension of staggered semi-implicit schemes to the DG framework has been derived in \cite{DumbserCasulli,2DSIUSW} for the shallow water equations on Cartesian and unstructured triangular grids, respectively. The resulting staggered semi-implicit DG method has furthermore been extended to the incompressible and compressible Navier-Stokes equations in two and three space dimensions, see  \cite{FambriDumbser,STINS2D,STINS2DTri,STINS3D,STCNS}, thus constituting a solid and well-established 
starting point for our research.\\ 
%Staggered grids
In our \textit{staggered} approach, the spatial discretization is performed through the use of two unstructured grids: a primary triangular grid for the discrete pressure, and a dual edge-based quadrilateral grid for the discrete velocity, see also \cite{Bermudez1998,Bermudez2014,USFORCE,USFORCE2}. While the use of staggered grids is a very common practice in the finite difference and finite volume framework (see e.g. \cite{casulli:1984,markerandcell,patankarspalding}), its use is not so common in the framework of DG methods. The first \textit{vertex-based} staggered DG schemes were proposed in \cite{CentralDG1,CentralDG2}. Other recent \textit{edge-based} staggered DG schemes have been proposed in \cite{StaggeredDGCE1,StaggeredDG}. The main advantage of using edge-based staggered grids is the sparsity pattern of the final linear system that has to be solved for the pressure. Alternatively, semi-implicit DG schemes on \textit{collocated} grids have been presented, for example, in \cite{Dolejsi1,Dolejsi2,Dolejsi3,GiraldoRestelli,TumoloBonaventuraRestelli}.\\
%Riemann solver
Another consequence of the use of a staggered formulation is that, since all quantities are readily defined where they are needed, we do not need any Riemann solvers (numerical flux functions), apart from the nonlinear convective terms. This special feature is not standard for DG schemes, which typically require numerical fluxes or penalty terms, see \cite{CBS-convection-diffusion,CBS-convection-dominated,MunzDiffusionFlux,YanShu}. For the nonlinear convective part of the incompressible Navier-Stokes equations, we use a standard DG scheme  based on the Rusanov flux \cite{Rusanov:1961a}. Also the viscous terms are discretized explicitly. The DG discretization of the viscous fluxes is based on the formulation of Gassner et al. \cite{MunzDiffusionFlux}, which was developed for the solution of the Generalized Riemann Problem (GRP) of the diffusion equation.\\
%Lagrangian-ALE
The introduction of our new pair of quadrature-free equal-order-interpolation velocity-pressure elements has the objective to reduce memory requirements and CPU cost of the previous staggered-semi-implicit DG schemes introduced in \cite{STINS2D,STINS3D}. While the scheme presented in this paper is still limited to fixed meshes (Eulerian approach), the proposed quadrature-free formulation has also the potential to be extended to moving unstructured meshes in the \textit{Arbitrary-Lagrangian-Eulerian} (ALE) context in the future. A lot of research about Lagrangian methods has been carried out in the last decades \cite{Benson1992a,Caramana1998,Despres2009,Maire2007,munz94,vonneumann50,Smith1999}, due to their excellent properties in the resolution of moving material interfaces and contact waves. In ALE schemes, see e.g. \cite{Feistauer4,Feistauer1,Feistauer3,Feistauer2,Hirt74,Peery2000,Smith1999}, the computational grid moves with an arbitrary mesh velocity that does not necessarily coincide with the real fluid velocity. For some examples of the so-called indirect cell-centered ALE algorithms, the reader is referred to \cite{MaireMM2,ShashkovRemap1,ShashkovRemap5,ShashkovRemap3,ShashkovRemap4,ShashkovCellCentered}. Many multi-phase and multi-material flow problems are solved relying on this approach \cite{ShashkovMultiMat3,ShashkovMultiMat1,Hirt74,Peery2000,ShashkovMultiMat4,Smith1999,ShashkovMultiMat2}. In a recent series of works \cite{ALEMQF,LagrangeMDRS,Lagrange2D,Lagrange3D,LagrangeDG,LagrangeQF,LagrangeMHD,ALELTS2D,ALEMOOD2,ALEMOOD1,ALELTS1D,Lagrange1D,KlingenALE2017}, a new family of high order accurate finite volume and DG schemes has been proposed in the ALE framework.\\
%Overview
The rest of the paper is organized as follows. In Section \ref{gov_eqn}, we briefly introduce the notation for the Navier-Sokes equations; Section \ref{new_stag} is dedicated to the presentation of the new staggered space-time discontinuous Galerkin finite elements; in Section \ref{SIDG_method} we provide a detailed description of the final solving algorithm; and in Section \ref{N_tests} the validation and the performance analysis of the method are illustrated through several numerical tests.

\section{The governing equations}
\label{gov_eqn}
The conservative form of the Navier-Stokes equations for an incompressible Newtonian fluid with constant 
density $\rho$ and constant dynamic viscosity coefficient $\mu$ in the domain $\Omega \subset \mathds{R}^d$ 
read as follows: 
\begin{eqnarray}
& &\nabla \cdot \mathbf{v}  =  0,  \label{eq:CS_2}, \\
& &\frac{\partial \mathbf{v}}{\partial t}+\nabla \cdot \mathbf{F}_c + \nabla p = \nabla \cdot \left( \nu \nabla  \mathbf{v} \right) + \mathbf{S},  \label{eq:CS_2_2_0}  
\end{eqnarray}
where $\mathbf{v}=\left( u , v \right)^T$ is the velocity vector, $u$ and $v$ are the velocity components in the $x$ and $y$ direction, respectively, $\mathbf{F}_c=\mathbf{v} \otimes \mathbf{v}$ is the flux tensor containing the nonlinear convective terms, and which in two space dimensions reads  
$$ \mathbf{F}_c=\left(\begin{array}{cc} u^2 & uv \\ uv & v^2 \end{array} \right); $$ 
then $p=P/\rho$ is the normalized fluid pressure, $P$ is the physical pressure, $\nu = \mu / \rho$ is the kinematic viscosity coefficient and, finally,
$\mathbf{S}=\mathbf{S}(\mathbf{v},x,y,t)$ is a generic source term representing the volumetric forces\footnote{In the present work, the source term is always ignored, therefore $\mathbf{S}=\mathbf{0}$. Nevertheless, we decide to keep the source term in the following in order to have a non-zero right hand side, see e.g. Equation \eqref{eq:CS_2_2}.}.\\
Following the same idea of \cite{DumbserNSE,MunzDiffusionFlux},
the viscosity term is grouped together with the nonlinear convective term.
In this way the momentum equation \eqref{eq:CS_2_2_0} can be rewritten as:
\begin{equation}
\frac{\partial \mathbf{v}}{\partial t}+\nabla \cdot \mathbf{F}_{cv} + \nabla p= \mathbf{S} 
\label{eq:CS_2_2},
\end{equation}
where $\mathbf{F}_{cv}=\mathbf{F}_{cv}(\mathbf{v},\nabla \mathbf{v})=\mathbf{F}_c(\mathbf{v})-\nu \nabla \mathbf{v}$ is a nonlinear tensor that depends not only on the velocity, but also on its gradient:
$$
\F_{cv} = \left( \begin{matrix}u^2 - \nu \diff{u}{x} & uv - \nu \diff{u}{y}  \\  uv - \nu \diff{v}{x}  & v^2 - \nu \diff{v}{y}  \end{matrix} \right).
$$
For a clear presentation of our numerical method, which involves a weak space--time formulation of the Navier-Stokes equations, we restyle equations \eqref{eq:CS_2} and \eqref{eq:CS_2_2} in the following
 equivalent, but more compact, form:
\begin{eqnarray}
& &\tilde{\nabla} \cdot \tilde{\mathbf{v}} =  0 \label{eq:CS_2_div}, \\
& &\tilde{\nabla} \cdot \tilde{\F}_{cv} + \tilde{\nabla} \cdot \tilde{\F}_{p}  =  \mathbf{S} \label{eq:CS_2_2_0_div}, 
\end{eqnarray}
by introducing the space--time divergence operator $\tilde{\nabla} = \left( \diff{}{t}, \diff{}{x}, \diff{}{y} \right)^T $, the vector $\tilde{\mathbf{v}} = \left( 0, u, v \right)^T $ and the tensors:
\begin{equation}
\label{fluxesT}
\tilde{\F}_{cv} = \left( \begin{matrix} u & u^2 - \nu \diff{u}{x} & uv - \nu \diff{u}{y} \\ v & uv - \nu \diff{v}{x} & v^2 - \nu \diff{v}{y}  \end{matrix} \right),
\quad 
\tilde{\F}_{p} = \left( \begin{matrix} 0 & p & 0 \\ 0 & 0 & p \end{matrix} \right) .
\end{equation}

\section{The new staggered space--time basis functions}
\label{new_stag}
The main novelty of our new numerical method with respect to previous similar work presented in \cite{STINS2D,STINS2DTri} consists in the use of a different kind of finite elements in the framework of a 
staggered space-time DG method.
An \textit{edge-based staggered} approach for two-dimensional domains is considered, see \cite{Bermudez1998,StaggeredDG,2DSIUSW,STINS2D,STINS2DTri,STCNS,USFORCE}, but here we propose a new pair of finite elements that is in principle also suitable for staggered ALE methods for the incompressible Navier-Stokes equations on moving meshes. \\
In paragraph \ref{stagg_grid} we briefly summarize the used notation; then paragraphs \ref{spatial_bf} and \ref{spcetime_bf} illustrate the spatial and the temporal discretization of the problem, respectively, and, finally, the space--time finite elements are defined in paragraph \ref{subparam}.
\subsection{Unstructured staggered grids}
\label{stagg_grid}
\begin{figure}[h]
	\begin{center}
		\includegraphics[width=0.95\textwidth]{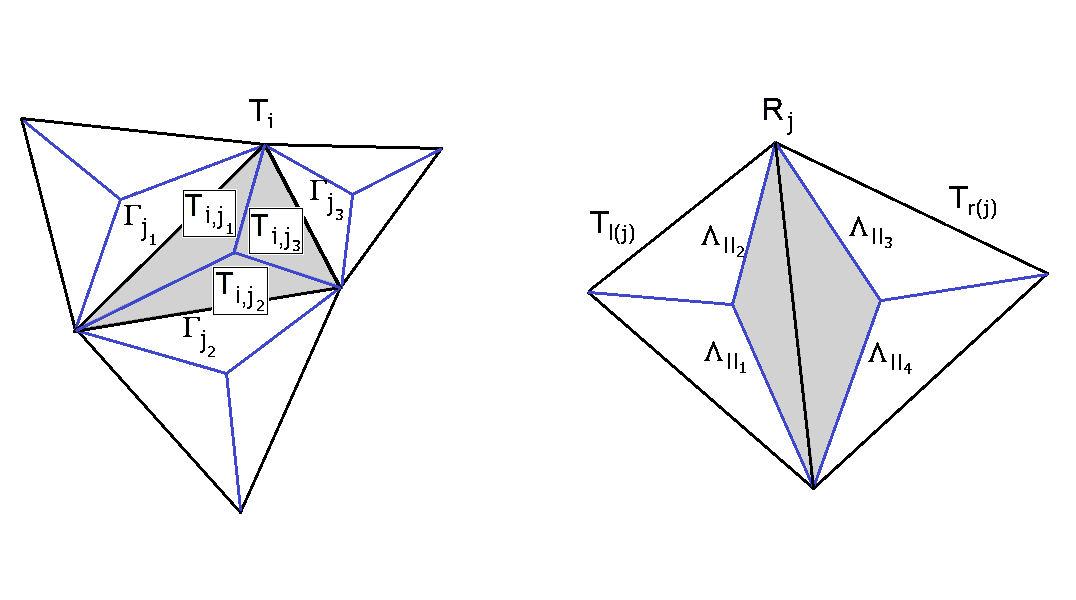}
		\caption{The notation used for the primary triangular mesh (on the left) and for the associated dual quadrilateral mesh (on the right).}
		\label{not_st_grid}
	\end{center}
\end{figure}

Let us consider the domain of the problem $\Omega\subset\mathbb{R}^2$ and its boundary $\partial\Omega$. For a graphical illustration of the following discussion, the reader is addressed to Figure \ref{not_st_grid}.\\
The computational domain $\Omega_{h}$ is covered with a triangulation $\mathcal{T}_{h}$ of $N_e$ non-overlapping triangular elements, denoted by $T_i$:
$$
\Omega_{h}=\bigcup_{i=1,\ldots,N_e}{T}_i \: .
$$
Let $N_d$ be the total number of edges of the triangles, and let us call $\Gamma_j$ the $j-$th edge of the primary grid, with $j=1, \ldots, N_d$. The subset of edges $\mathcal{B}(\Omega_{h})$ is defined as the set of the indices $j$ corresponding to the boundary edges:
$$
\mathcal{B}(\Omega_{h}) = \{j \in [1,N_d] : \Gamma_{j} \subset \partial\Omega  \} .
$$
For every $i=1,\ldots,N_e$, we denote by $S_i$ the set of the three edges belonging to triangle ${T}_i$: 
$$S_i=\{j \in [1,N_d]  : \Gamma_j \subset \partial{T}_i \} \quad \Rightarrow \quad |S_{i}|=3.$$ 
For every internal edge with index $j\in [1, N_d]\setminus\mathcal{B}(\Omega_h)$, the corresponding dual quadrilateral element is called ${R}_j$ and is built by taking as vertices the two barycenters of the triangles sharing the edge (the one on the left is called ${T}_{l(j)}$ and the one on the right is called ${T}_{r(j)}$) and the two vertices of $\Gamma_j$. \\
Let us denote by ${T}_{i,j}={R}_j \cap {T}_i$ the sub-triangle which is in common between element ${T}_i$ and element ${R}_j$, for every $i=1,\ldots,N_e$ and for every $j \in S_i$.
Notice that, if $j \in \mathcal{B}(\Omega_h)$, then a special triangular element ${R}_j$ will be part of the dual staggered grid: it is defined by the two vertices of $\Gamma_{j}$ and the barycenter of triangle ${T}_{i}$, with $j\in S_{i}$. In particular, in this case they are equivalent: ${R}_j = {T}_{i,j}  \: .$ \\
According to \cite{STINS2D}, the mesh of triangular elements $\{ {T}_i \}_{i \in [1, N_e]}$ is called the \textit{main grid} or \textit{primary grid} and will be used for the approximation of the pressure $p_h$, while the set of the elements $\{{R}_j \}_{j \in [1, N_d]}$ is called the \textit{dual grid} and will be used for the approximation of the velocity $\mathbf{v}_{h}$.\\
The notation for the analogous quantities of the dual grid is the following: $N_l$ indicates the total number of edges; $\Lambda_{ll}$ indicates the $ll$-th physical edge; and for every $j=1,\ldots,N_d$, we denote by $E_j$ the set of the edges belonging to element ${R}_j$: 
\begin{equation*}
E_j=\{ll \in [1,N_l]  : \Lambda_{ll} \subset \partial{R}_j \} \quad \Rightarrow \quad \begin{cases}
|E_{j}|=4 \quad \text{if } j \notin \mathcal{B}(\Omega_{h}) \\
|E_{j}|=3 \quad \text{if } j \in \mathcal{B}(\Omega_{h})
\end{cases} 
.
\end{equation*}
\subsection{Spatial basis functions}
\label{spatial_bf}
As usual in the framework of finite element methods, we first define the basis functions on a reference element and then we can extend them to the physical space of the whole computational domain. We employ the standard  \textit{nodal} or \textit{"Lagrangian"} approach for the definition of the basis functions.

\begin{figure}[htbp] 
	\begin{center}
		\subfloat[$\phi_{1}(\xi,\eta)$]{\includegraphics[width=0.35\textwidth]{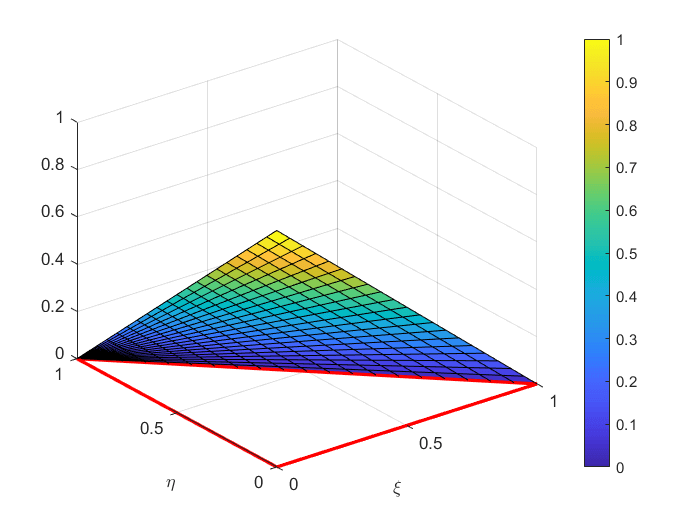}}
		\subfloat[$\phi_{2}(\xi,\eta)$]{\includegraphics[width=0.35\textwidth]{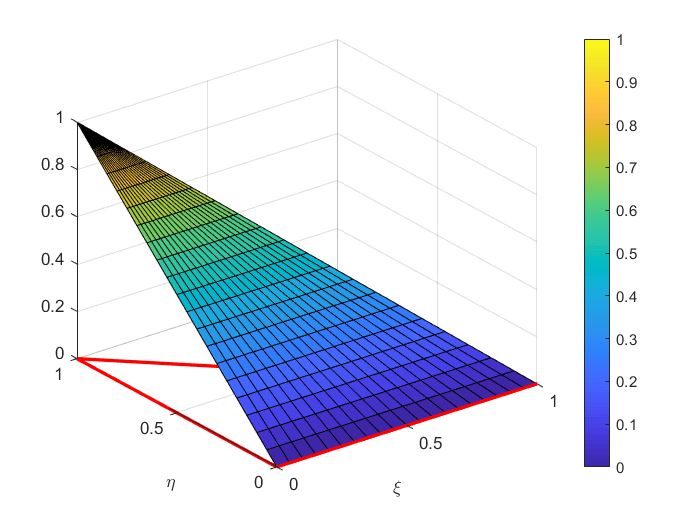}}
		\\
		\subfloat[$\phi_{3}(\xi,\eta)$]{\includegraphics[width=0.35\textwidth]{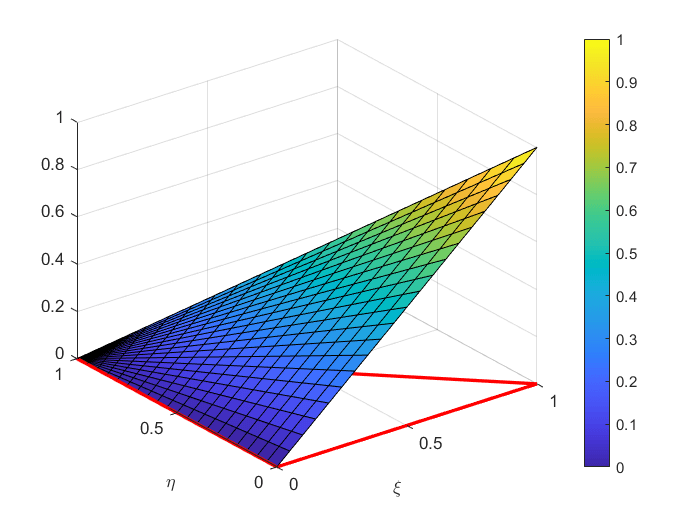}}
		\caption{The $ N=1 $ basis functions of the primary grid.}
		\label{p1-phi}
	\end{center}
\end{figure}

\begin{figure}[htbp] 
	\begin{center}
		\subfloat[$\phi_{1}(\xi,\eta)$]{\includegraphics[width=0.35\textwidth]{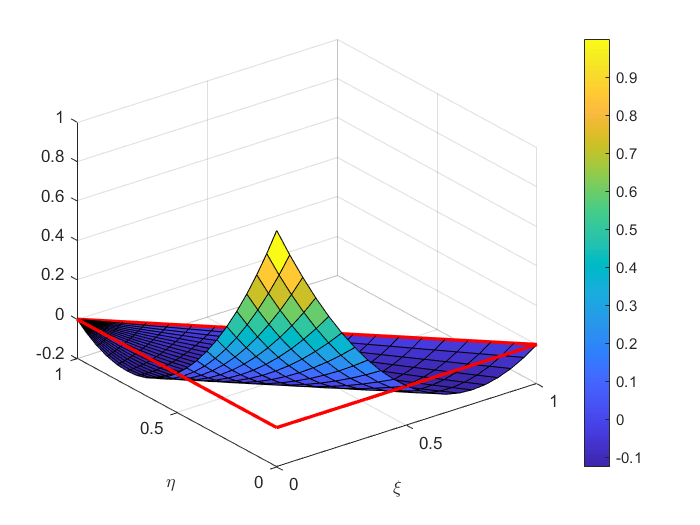}}
		\subfloat[$\phi_{2}(\xi,\eta)$]{\includegraphics[width=0.35\textwidth]{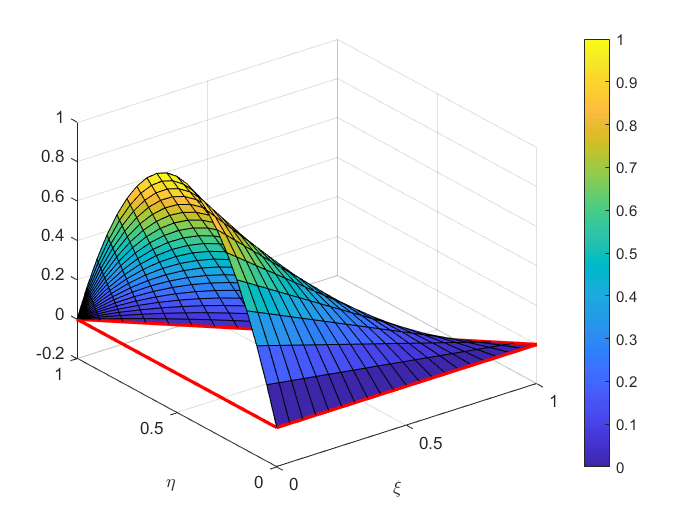}}
		\\
		\subfloat[$\phi_{3}(\xi,\eta)$]{\includegraphics[width=0.35\textwidth]{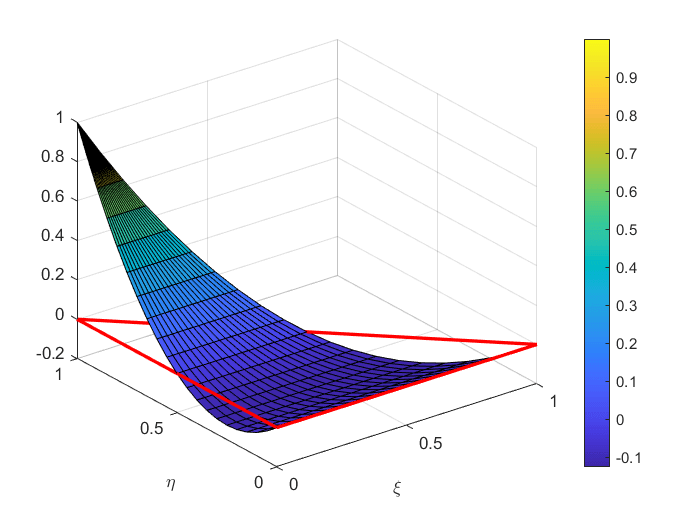}}
		\subfloat[$\phi_{4}(\xi,\eta)$]{\includegraphics[width=0.35\textwidth]{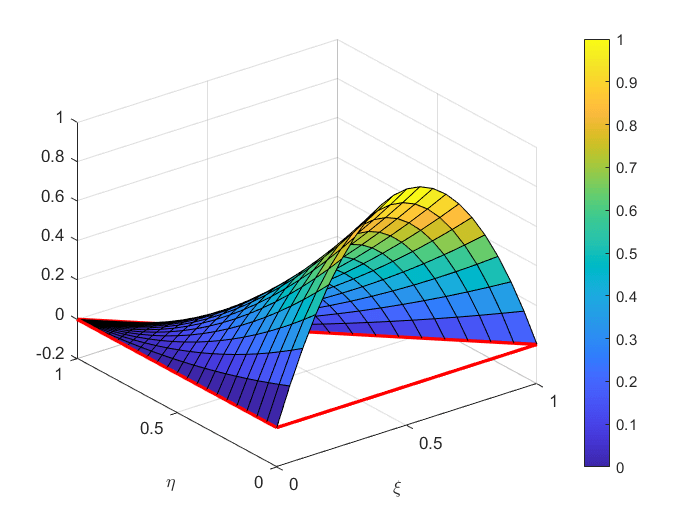}}
		\\
		\subfloat[$\phi_{5}(\xi,\eta)$]{\includegraphics[width=0.35\textwidth]{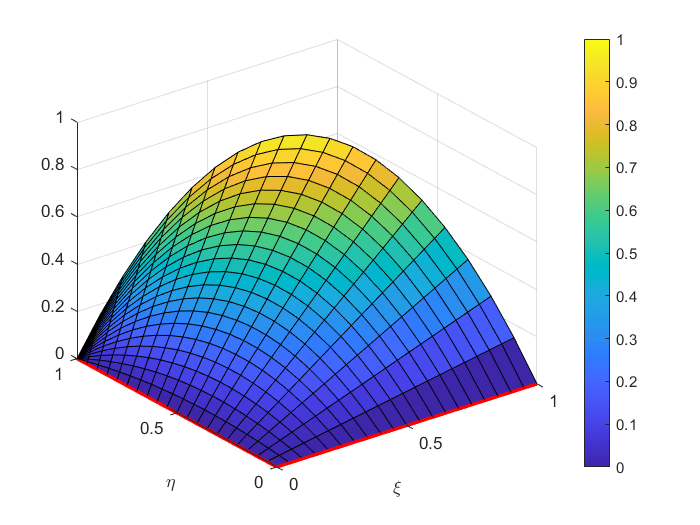}}
		\subfloat[$\phi_{6}(\xi,\eta)$]{\includegraphics[width=0.35\textwidth]{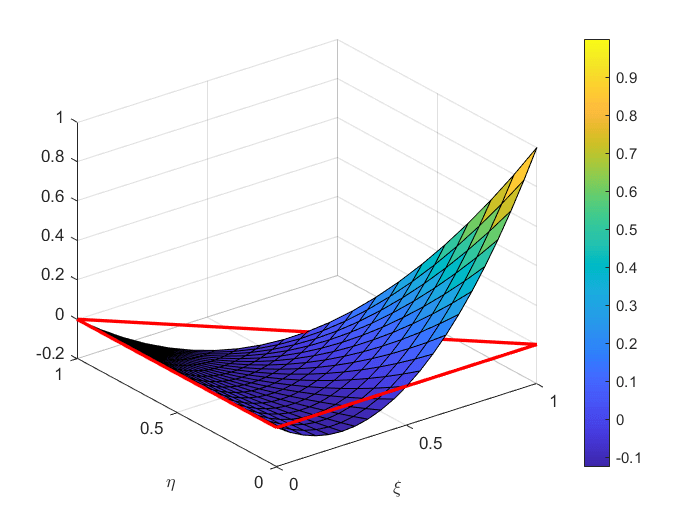}}
		\caption{The $ N=2 $ basis functions of the primary grid.}
		\label{p2-phi}
	\end{center}
\end{figure}

For the \textit{primary} finite element space of polynomial approximation degree $N$, let us define the  triangular reference element  $$T_{std}=\{(\xi,\eta) \in \mathbb{R}^{2,+} : \eta\leq1-\xi \wedge 0 \leq \xi \leq 1 \} ,$$ 
and let us impose the classical Lagrange interpolation condition $\phi_k(\vec \xi_l) = \delta_{kl}$ over the nodal points:
\begin{equation}
\vec \xi_k = (\xi_{k_1}, \eta_{k_2}) = \left( \frac{k_1}{N}, \frac{k_2}{N} \right),
\label{eqn.nodes}
\end{equation}
where $\delta$ indicates the Kronecker delta, $k=(k_1,k_2)$ is a multi-index and the ranges for indices $k_1$, $k_2$ are $0 \leq k_1 \leq N$ and $0 \leq k_2 \leq N - k_1$.
In total, $N_\phi=\frac{(N+1)(N+2)}{2}$ basis
functions are obtained for the reference element of the primary grid and they are labeled as $\{\phi_k \}_{k \in [1,N_\phi]}$. A visualization of them can be found in Figure \ref{p1-phi} for $N=1$ and in Figure \ref{p2-phi} for $N=2$.

\begin{figure}[htbp] 
	\begin{center}
		\subfloat[$\psi_{1}(\xi,\eta)$]{\includegraphics[width=0.35\textwidth]{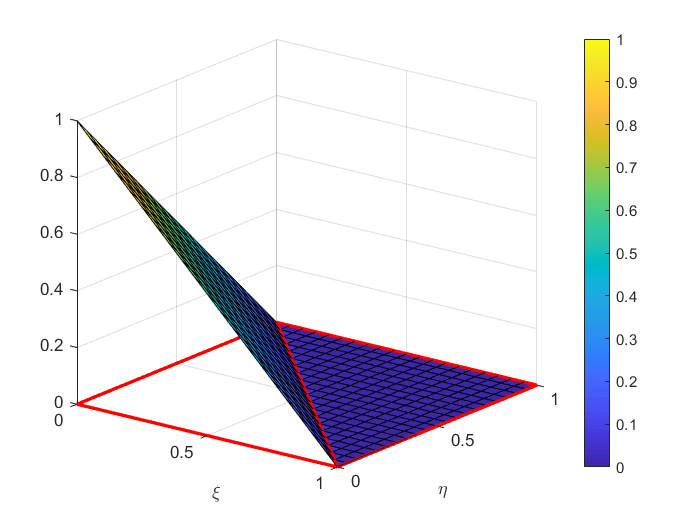}}
		\subfloat[$\psi_{2}(\xi,\eta)$]{\includegraphics[width=0.35\textwidth]{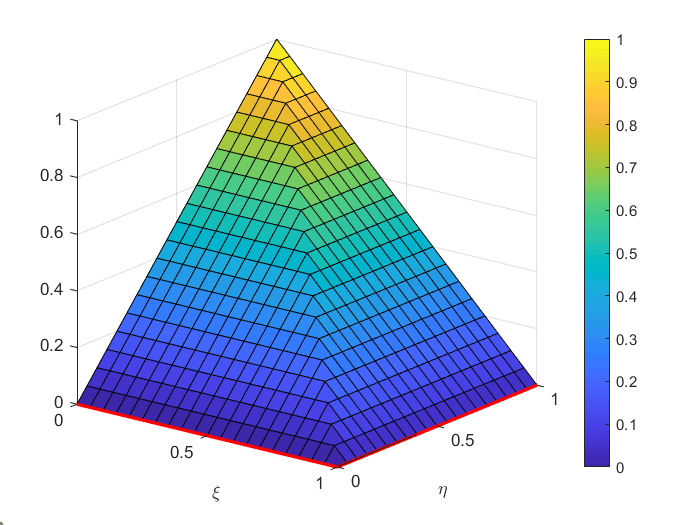}}
		\\
		\subfloat[$\psi_{3}(\xi,\eta)$]{\includegraphics[width=0.35\textwidth]{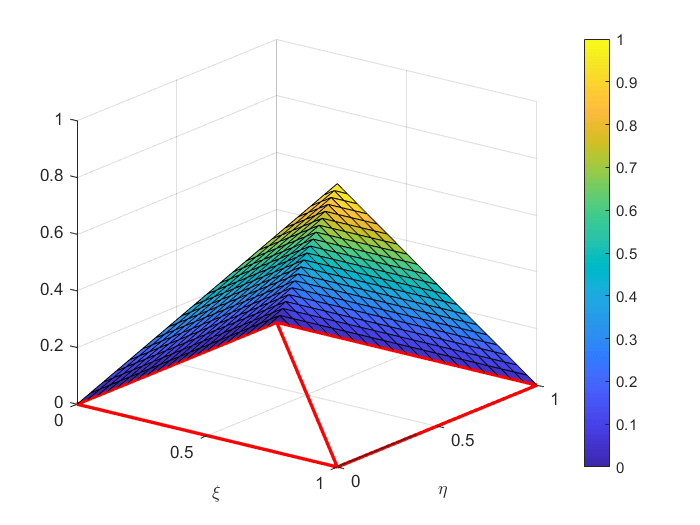}}
		\subfloat[$\psi_{4}(\xi,\eta)$]{\includegraphics[width=0.35\textwidth]{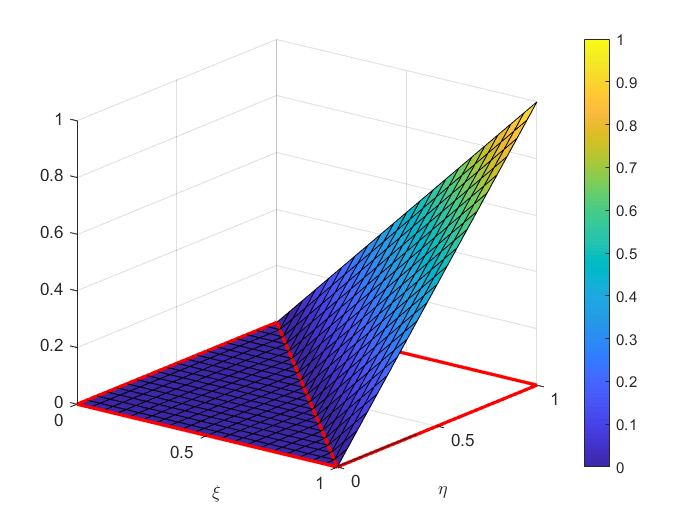}}
		\caption{The new $ N=1 $ basis functions of the dual grid proposed in this paper.}
		\label{p1-isop2}
	\end{center}
\end{figure}

\begin{figure}[htbp] 
	\begin{center}
		\subfloat[$\psi_{1}(\xi,\eta)$]{\includegraphics[width=0.35\textwidth]{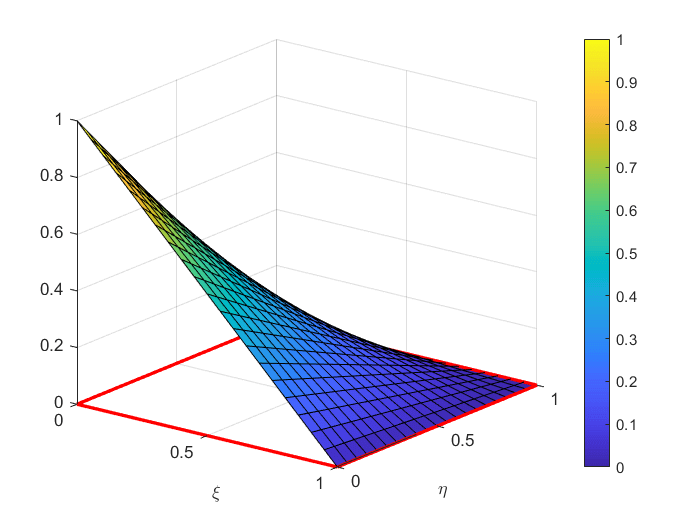}}
		\subfloat[$\psi_{2}(\xi,\eta)$]{\includegraphics[width=0.35\textwidth]{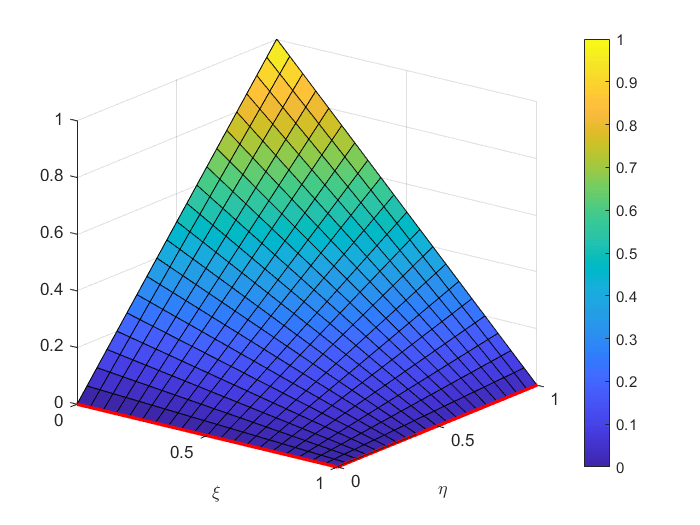}}
		\\
		\subfloat[$\psi_{3}(\xi,\eta)$]{\includegraphics[width=0.35\textwidth]{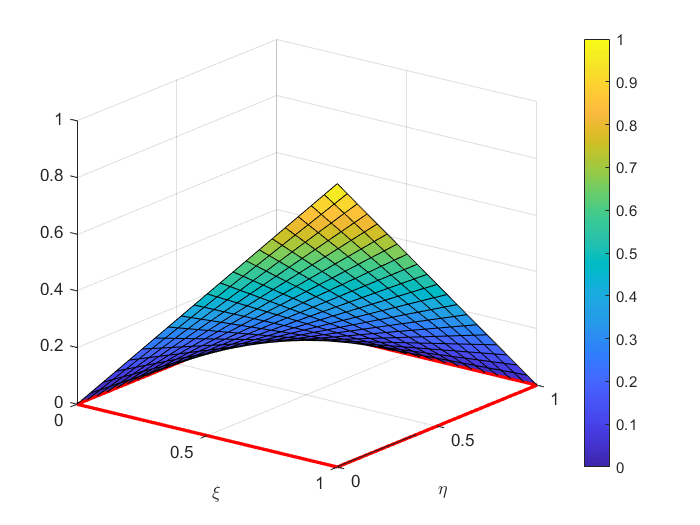}}
		\subfloat[$\psi_{4}(\xi,\eta)$]{\includegraphics[width=0.35\textwidth]{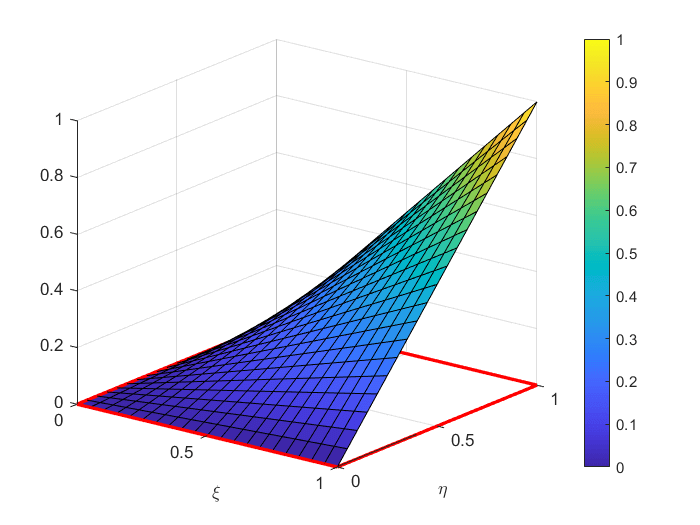}}
		\caption{The previous $ N=1 $ basis functions on the dual grid used in \cite{STINS2D,STINS2DTri}.}
		\label{p1-psi-old}
	\end{center}
\end{figure}

\begin{figure}[htbp] 
	\begin{center}
		\subfloat[$\psi_{1}(\xi,\eta)$]{\includegraphics[width=0.22\textwidth]{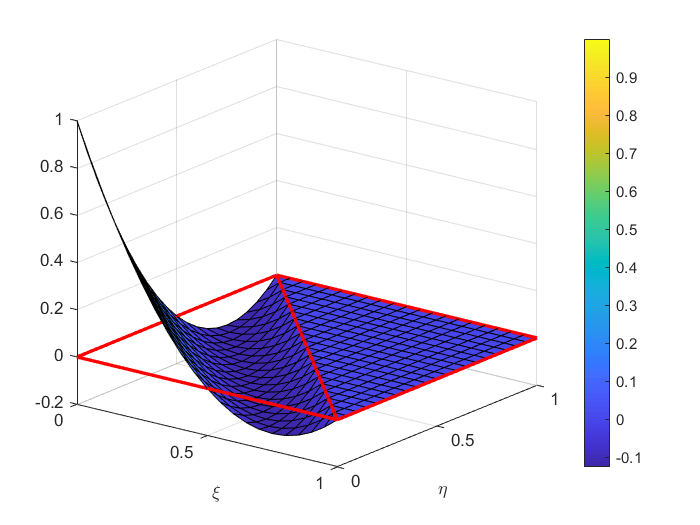}}
		\subfloat[$\psi_{2}(\xi,\eta)$]{\includegraphics[width=0.22\textwidth]{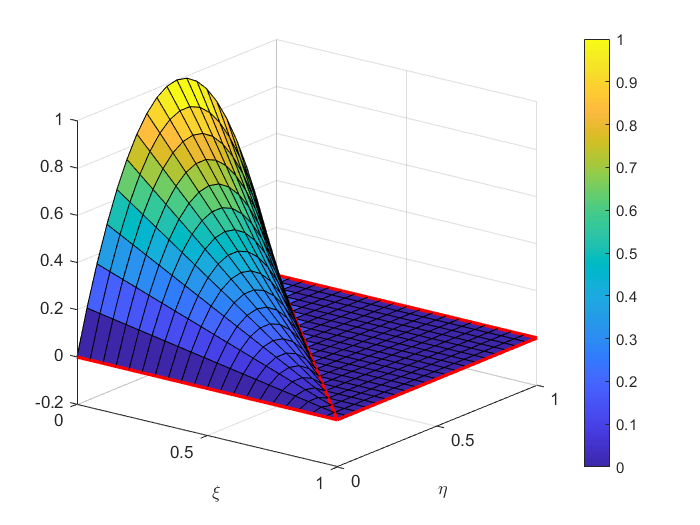}}
		\subfloat[$\psi_{3}(\xi,\eta)$]{\includegraphics[width=0.22\textwidth]{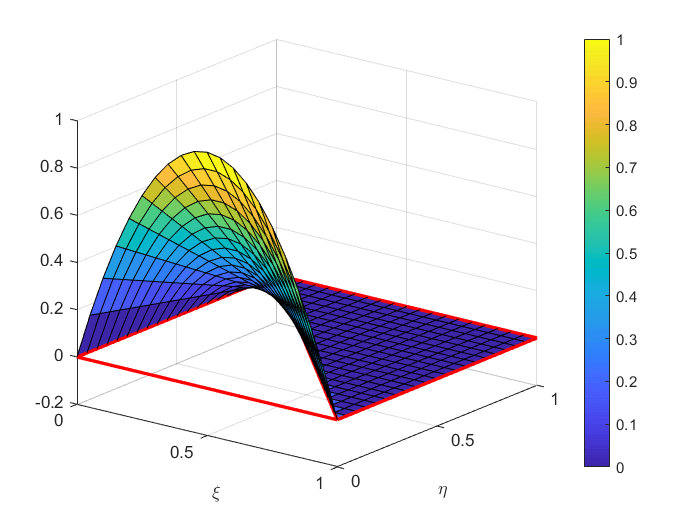}}
		\\
		\subfloat[$\psi_{4}(\xi,\eta)$]{\includegraphics[width=0.22\textwidth]{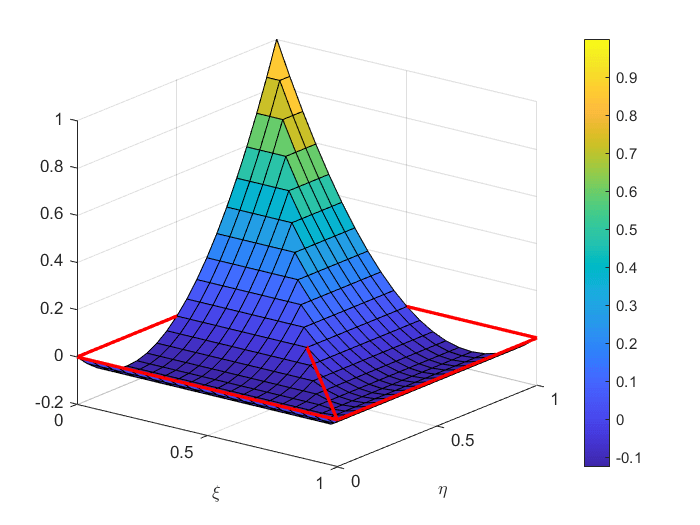}}
		\subfloat[$\psi_{5}(\xi,\eta)$]{\includegraphics[width=0.22\textwidth]{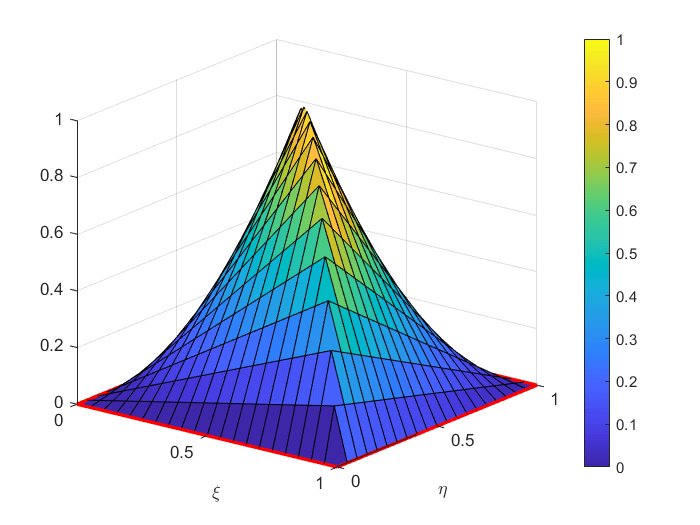}}
		\subfloat[$\psi_{6}(\xi,\eta)$]{\includegraphics[width=0.22\textwidth]{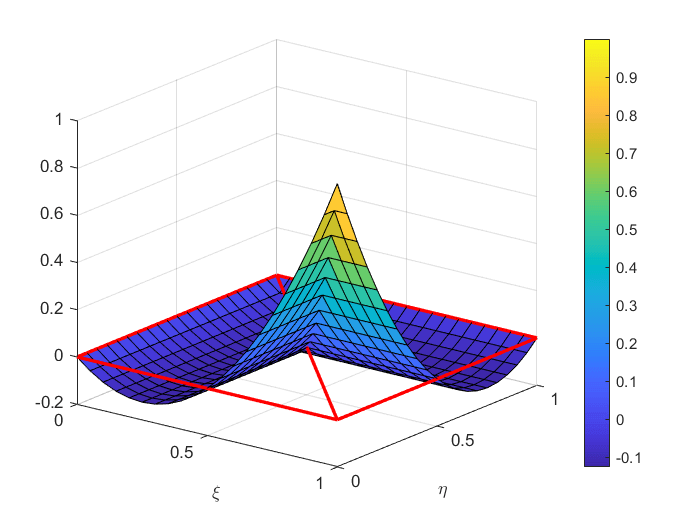}}
		\\
		\subfloat[$\psi_{7}(\xi,\eta)$]{\includegraphics[width=0.22\textwidth]{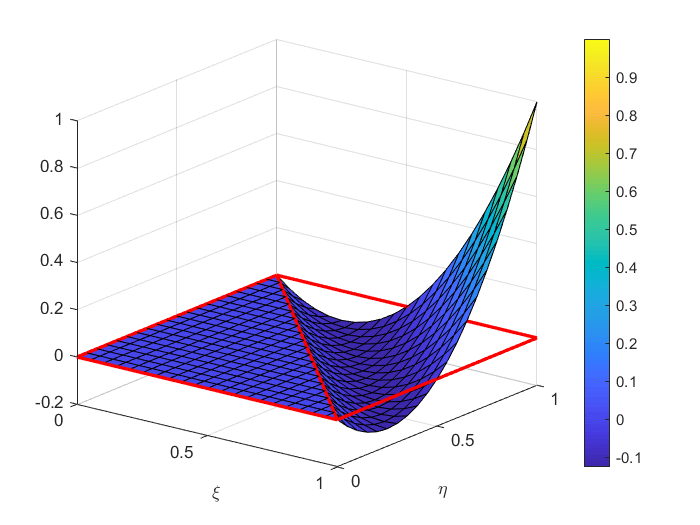}}
		\subfloat[$\psi_{8}(\xi,\eta)$]{\includegraphics[width=0.22\textwidth]{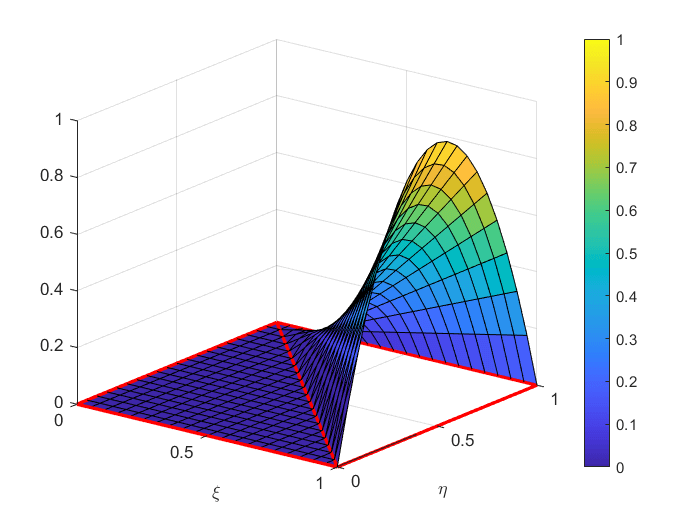}}
		\subfloat[$\psi_{9}(\xi,\eta)$]{\includegraphics[width=0.22\textwidth]{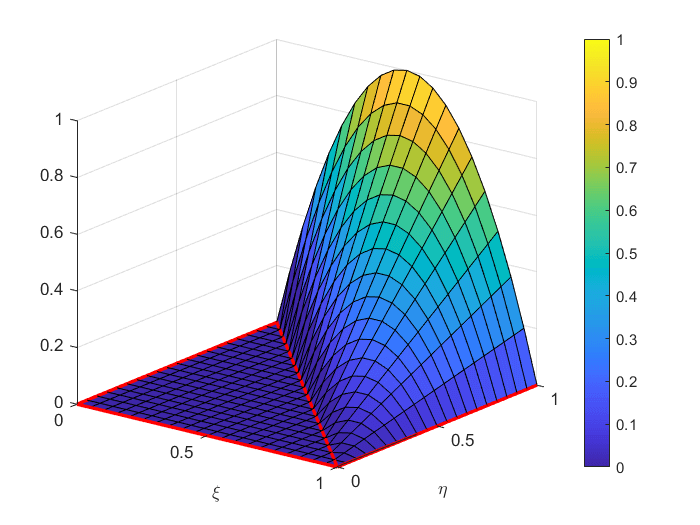}}
		\caption{The new $ N=2 $ basis functions of the dual grid proposed in this paper.}
		\label{p2-isop3}
	\end{center}
\end{figure}

\begin{figure}[htbp] 
	\begin{center}
		\subfloat[$\psi_{1}(\xi,\eta)$]{\includegraphics[width=0.22\textwidth]{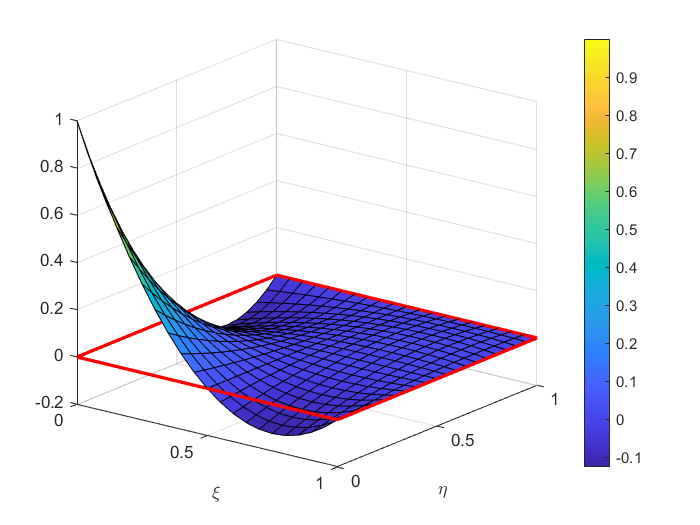}}
		\subfloat[$\psi_{2}(\xi,\eta)$]{\includegraphics[width=0.22\textwidth]{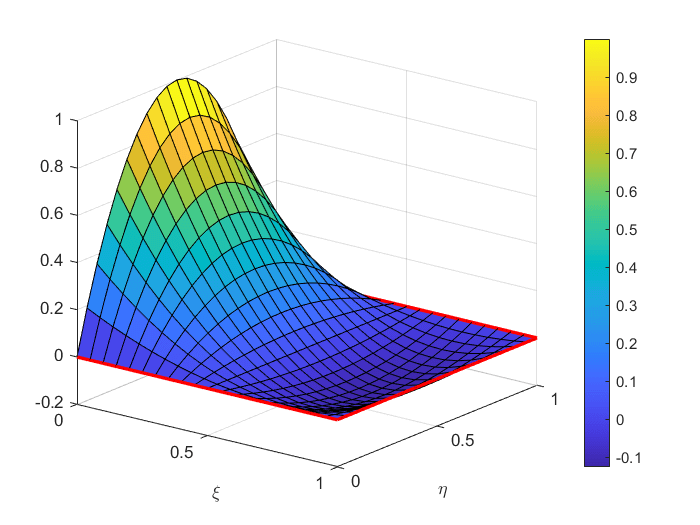}}
		\subfloat[$\psi_{3}(\xi,\eta)$]{\includegraphics[width=0.22\textwidth]{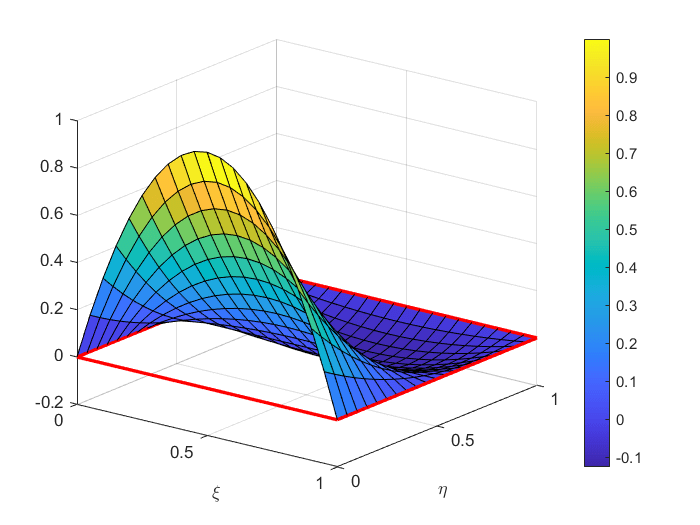}}
		\\
		\subfloat[$\psi_{4}(\xi,\eta)$]{\includegraphics[width=0.22\textwidth]{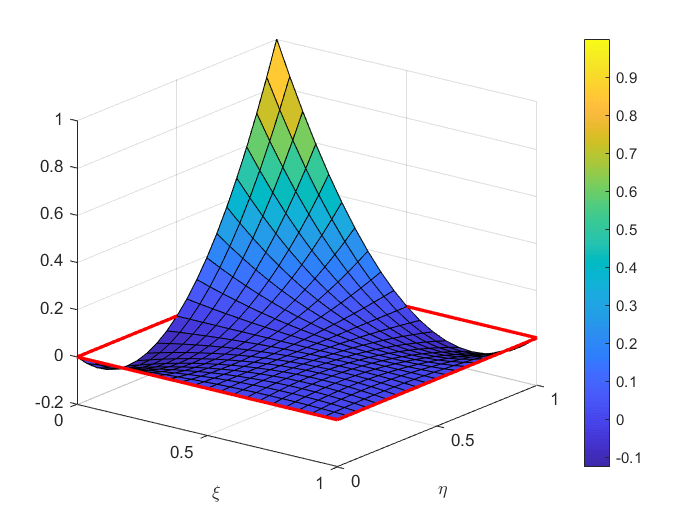}}
		\subfloat[$\psi_{5}(\xi,\eta)$]{\includegraphics[width=0.22\textwidth]{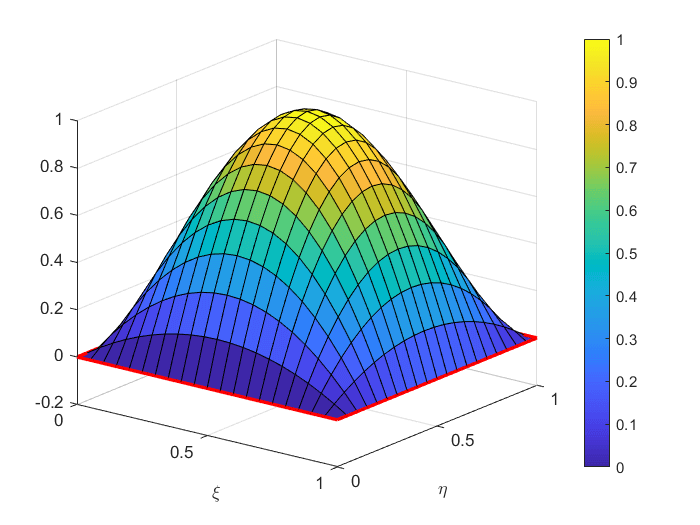}}
		\subfloat[$\psi_{6}(\xi,\eta)$]{\includegraphics[width=0.22\textwidth]{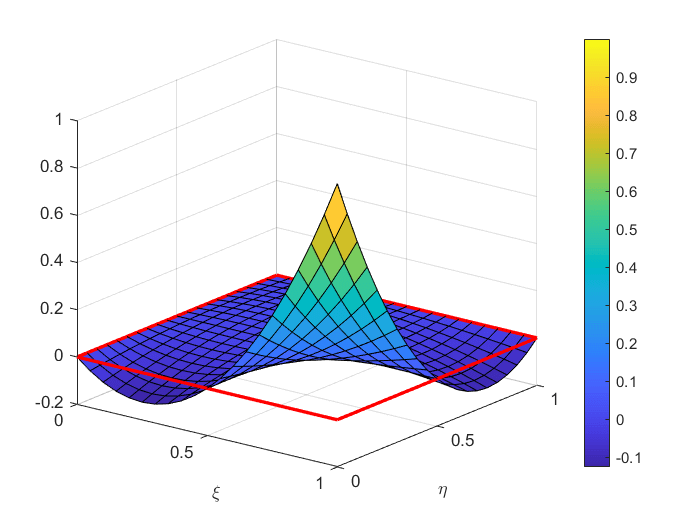}}
		\\
		\subfloat[$\psi_{7}(\xi,\eta)$]{\includegraphics[width=0.22\textwidth]{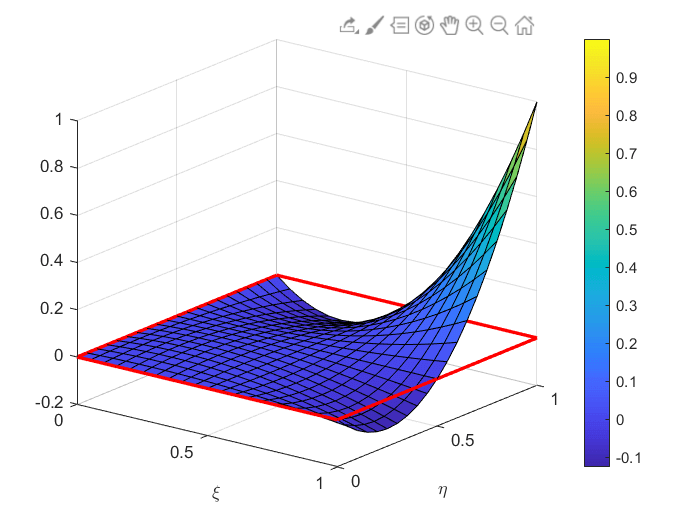}}
		\subfloat[$\psi_{8}(\xi,\eta)$]{\includegraphics[width=0.22\textwidth]{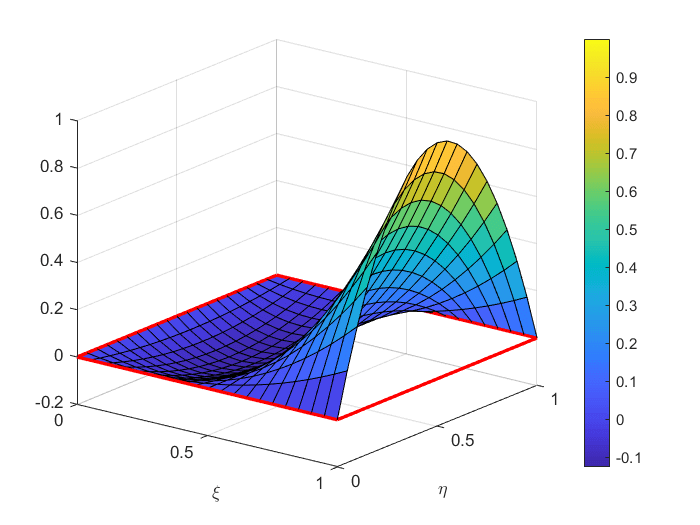}}
		\subfloat[$\psi_{9}(\xi,\eta)$]{\includegraphics[width=0.22\textwidth]{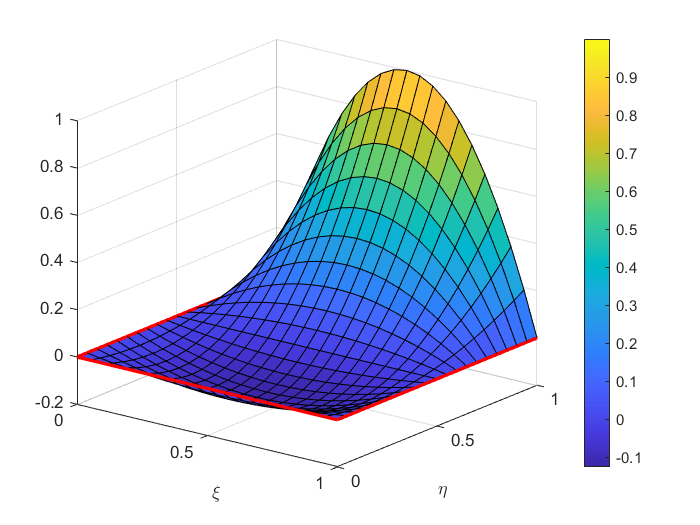}}
		\caption{The previous $ N=2 $ basis functions on the dual grid used in \cite{STINS2D,STINS2DTri}.}
		\label{p2-psi-old}
	\end{center}
\end{figure}

Similarly, for the \textit{dual} finite element space of polynomial approximation degree $N$, one could obtain another set of nodal basis functions on the reference square element $$R_{std}=\{(\xi,\eta) \in \mathbb{R}^{2,+} : 0 \leq \xi \leq 1 \wedge 0 \leq \eta \leq 1 \}$$
with the nodes given by the tensor product of $N+1$ one-dimensional quadrature points, for example the Newton-Cotes ones, see \cite{stroud}, as it was done in \cite{STINS2D,STINS2DTri}.
Instead, in this paper we define the basis functions on $R_{std}$ in a \textit{different} way, using 
the continuous union of piecewise polynomials of degree $N$ defined on each of the two subtriangles which constitute the reference square. Let us consider the square $R_{std}$ as the union of two sub-triangles $T_{I}$ and $T_{II}$, with $T_{I}=T_{std}$, and let us define the new basis functions following the standard nodal approach of conforming continuous finite elements with order $N$ inside the square ($T_{I}$ and $T_{II}$ are the so called "mini-elements"). This is an approach very similar to the $ \mathbb{P}^{k} - \text{iso} \mathbb{P}^{k+1} $ finite elements, described in \cite{ern2004}, which are applied for the velocity approximation in some mixed problems (for example the Stokes problem, where mixed finite element approximations are employed in order to numerically satisfy the \textit{inf-sup} or Ladyzhenskaya-Babu\v{s}ka-Brezzi compatibility condition). But, instead of creating new elements inside a triangular macro-element, in our staggered approach the macro-element is given by the square, and the two mini-elements naturally arise from the \textit{edge-staggered} geometric configuration. \\
In other words, we build the basis functions over the nodes which lie on sub-triangle $T_{I}$ exactly like we did for $T_{std}$, and then we extend them to $T_{II}$ with zero, in order to mantain the continuity inside the square. Viceversa, we define the basis functions over the nodes which lie on $T_{II}$ via the transformation $\sigma_{(\frac{1}{2},\frac{1}{2})}$ between $T_{II}$ and $T_{I}$, and then we extend them continuously to $T_{I}$ with zero. Here, $\sigma_{(\frac{1}{2},\frac{1}{2})}$ indicates the central symmetry with center $\left(\frac{1}{2},\frac{1}{2}\right)$ in the $\xi - \eta$ plane. In conclusion, $N_\psi=(N+1)^{2}$ basis functions are generated for the dual grid, i.e. for the spatial discretization of the velocity components, and they are labelled as $\{\psi_{k} \}_{k \in [1,N_\psi]}$. Their formal piecewise definition is as follows.
\begin{enumerate}
	\item If $1\leq k \leq \frac{N (N+1)}{2}$ ,
	\begin{equation}
	\psi_{k}(\vec \xi) = 
	\begin{cases}
	\phi_k(\vec \xi) \quad & \text{if } \vec \xi \in T_{I}\\
	0 \quad & \text{if } \vec \xi \in T_{II}
	\end{cases}
	\: ;
	\label{psi_k_1}
	\end{equation}
	\item if $\frac{N (N+1)}{2} < k \leq N_{\phi}$ ,
	\begin{equation}
	\psi_{k}(\vec \xi) = 
	\begin{cases}
	\phi_k(\vec \xi) \quad & \text{if } \vec \xi \in T_{I}\\
	\phi_{(N+1)^{2}-(k-1)}\left(\sigma_{(\frac{1}{2},\frac{1}{2})} (\vec \xi)\right) \quad & \text{if } \vec \xi \in T_{II}
	\end{cases}
	\: ;
	\label{psi_k_2}
	\end{equation}
	\item if $N_{\phi} < k \leq N_{\psi}$ ,
	\begin{equation}
	\psi_{k}(\vec \xi) =
	\begin{cases}
	0 \quad & \text{if } \vec \xi \in T_{I}\\ \phi_{k-N_{\phi}}\left(\sigma_{(\frac{1}{2},\frac{1}{2})}(\vec \xi)\right) \quad & \text{if } \vec \xi \in T_{II}
	\end{cases}
	\: .
	\label{psi_k_3}
	\end{equation}
\end{enumerate} 
A plot of them is available in Figure \ref{p1-isop2} for the case $ N=1 $ and in Figure \ref{p2-isop3} for the case $ N=2 $. A comparison with the old basis functions used in \cite{STINS2D,STINS2DTri} can be found in Figures \ref{p1-psi-old} and \ref{p2-psi-old}.

One clear advantage of this choice, from the programming point of view, is that the basis functions for the dual grid can be obtained by a simple manipulation of the basis functions of the primary grid. Indeed, notice that the primary grid and the dual grid have got the same degree of interpolation $N$. For this reason, we call this new pair of finite elements \textit{equal-order-interpolation} velocity-pressure elements.\\
Another crucial point is that the new velocity-pressure elements are in principle also suitable for an ALE implementation of a numerical solver, in which two important issues are the RAM memory required by the simulations and the computational time spent by the updating of the matrices in the overall algorithm. 
One of the purposes of this paper is to show that the ALE (moving meshes) implementation of our high-order accurate staggered numerical method with the new pair of velocity-pressure elements yields even better performance compared with the Eulerian (fixed mesh) implementation of the original method \cite{STINS2DTri}, concerning the two above-mentioned aspects. The performance comparisons of the two versions of the numerical solver are reported in detail in Section \ref{N_tests}.\\
The improved performace of the new method are indeed possible thanks to the new piecewise definition of the dual basis functions for the definition of the discrete velocity. As it will be clear when looking at the algebraic formulation of the problem (Paragraph \ref{algebr_f}), all the matrices of the final system to be solved are related to one and only one subset of the dual basis functions, hence to one and only one of the definitions in Formulae \eqref{psi_k_1}, \eqref{psi_k_2}, \eqref{psi_k_3}, namely: (i) the basis functions associated to the nodes of the left sub-triangle $T_{I}$, or (ii) the basis functions associated to the nodes belonging to the diagonal of $R_{std}$ connecting the points $(0,1)$ and $(1,0)$, or (iii) the basis functions associated to the nodes of the right sub-triangle $T_{II}$. It turns out that two great advantages are achieved.\\
First, given the approximation degree $N$, only a relatively small, constant number of \textit{universal} matrices computed over the $T_{std}$ reference element (for the volume space--time integral matrices), or over the unit interval $I_{std}=[0,1]$ (for the surface space--time integral matrices) are needed to be stored by the program. This means that only a universal version of the matrices and not all the local matrices  must be stored, therefore the new velocity-pressure elements allow to save a great amount of RAM memory, especially for simulations employing very refined meshes.\\
Secondly, the task of updating all the local matrices at each time step in the ALE algorithm is carried out through a simple tensor-reduction operation, where the local matrices are computed \textit{on the fly} from the universal matrices and from the geometric information of the moving meshes. Of course, this operation is much faster than recomputing from scratch all the integrals of the matrices, and the advantage in terms of computational time is significant when employing high-order finite elements, which would require very expensive numerical quadrature formulae in the recalculation of the matrices for every time step of the ALE simulation.

\subsection{Extension to the space--time basis functions}
\label{spcetime_bf}
\begin{figure}[htbp] 
	\begin{center}
		\subfloat[$\gamma_{1}(\tau)$]{\includegraphics[width=0.5\textwidth]{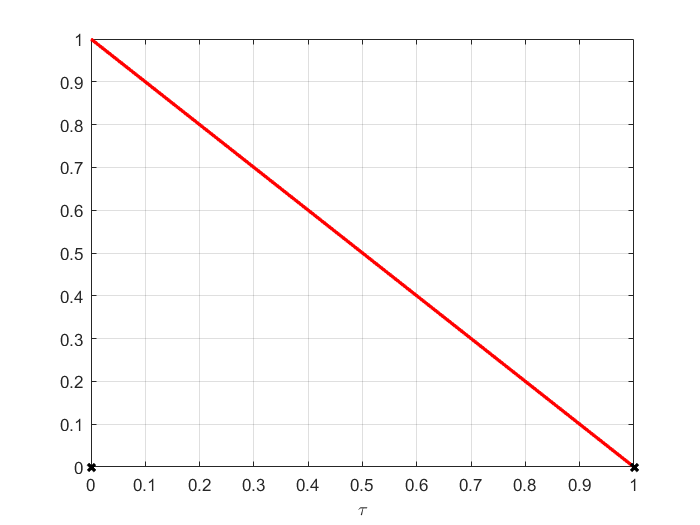}}
		\subfloat[$\gamma_{2}(\tau)$]{\includegraphics[width=0.5\textwidth]{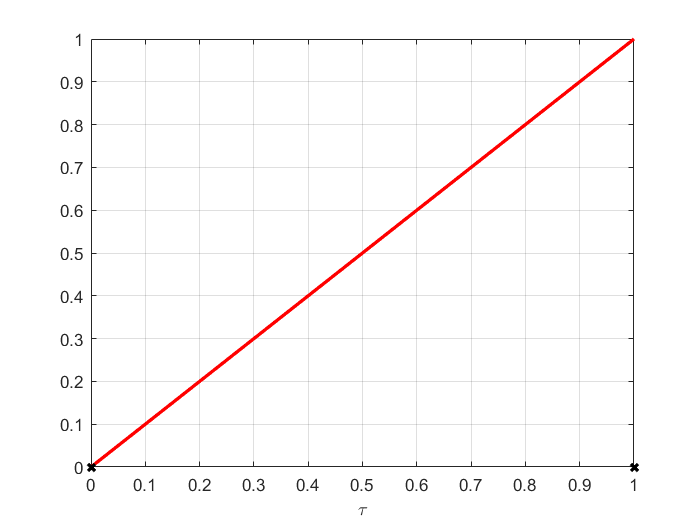}}
		\caption{The $ p_{\gamma}=1 $ basis functions on $I_{std}$.}
		\label{p1-gamma}
	\end{center}
\end{figure}

\begin{figure}[htbp] 
	\begin{center}
		\subfloat[$\gamma_{1}(\tau)$]{\includegraphics[width=0.5\textwidth]{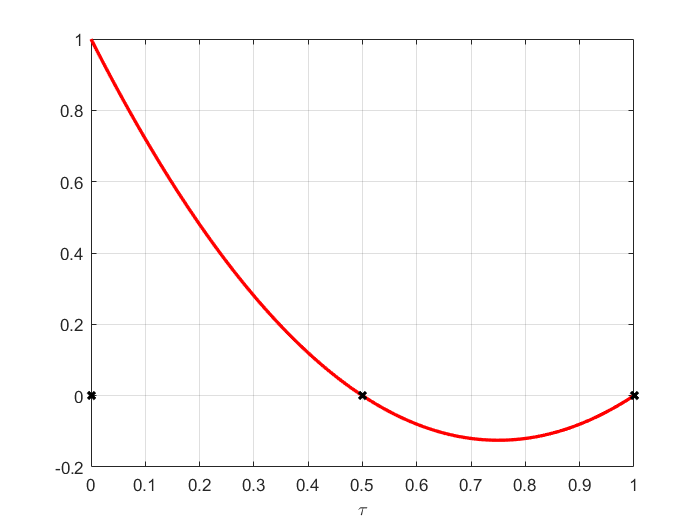}}
		\subfloat[$\gamma_{2}(\tau)$]{\includegraphics[width=0.5\textwidth]{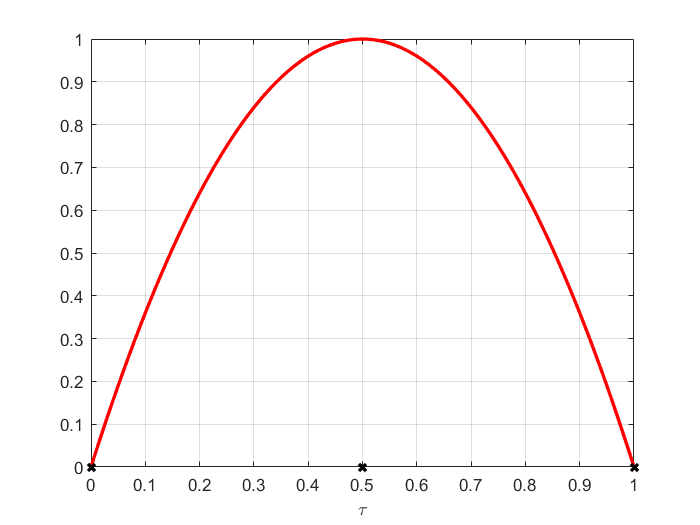}}
		\\
		\subfloat[$\gamma_{3}(\tau)$]{\includegraphics[width=0.5\textwidth]{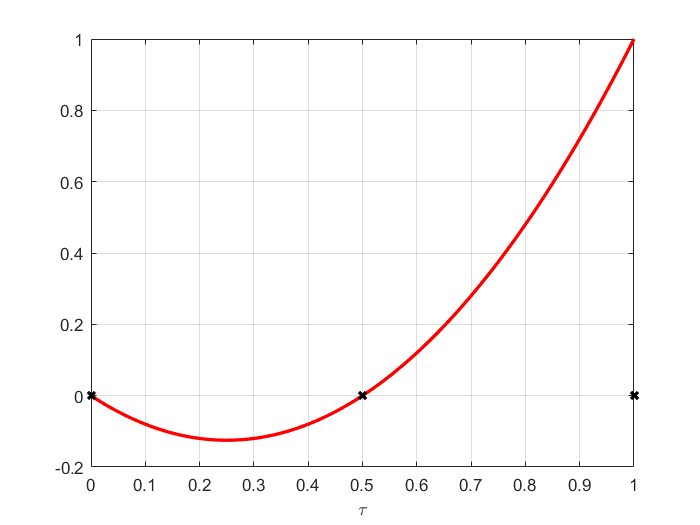}}
		\caption{The $ p_{\gamma}=2 $ basis functions on $I_{std}$.}
		\label{p2-gamma}
	\end{center}
\end{figure}

Since we apply high-order accurate methods also in time, our numerical scheme described in Section \ref{SIDG_method}, is of the kind of the \textit{space--time} DG schemes, see e.g.   \cite{KlaijVanDerVegt,Rhebergen2013,spacetimedg1,spacetimedg2}. Following \cite{STINS2DTri,STCNS,STDGLE}, we have to define the space--time basis functions as an extension of the spatial basis functions introduced in the previous paragraph. \\
If the time interval of our problem starts at instant $0$ and ends at instant $T$, we define a number of $N_{t}$ timesteps over the sequence of times:
$$0=t^0<t^1<t^2< \ldots <t^{N_{t}}=T ,$$
and we use the notation $T^{n+1}=[t^{n}, t^{n+1}]$ for the time interval and $\Delta t^{n+1} = t^{n+1}-t^{n} $ for the timestep length, for each $n=0, \ldots, N_{t}-1$.\\
For polynomials of degree $p_\gamma$, let us define the temporal basis functions on the reference interval $I_{std} = [0,1]$ as the Lagrange interpolation polynomials passing through the equidistant 1D Newton-Cotes quadrature nodes \cite{stroud}. The resulting $N_\gamma=p_\gamma+1$ basis functions are labelled as $\{\gamma_k\}_{k \in [1, N_\gamma]}$. A plot of them can be found in Figures \ref{p1-gamma} and \ref{p2-gamma}, for the cases $p_{\gamma}=1$ and $p_{\gamma}=2$, respectively. \\
Using the tensor product between the new spatial basis functions introduced in the previous paragraph and the temporal basis functions, we can define the basis functions on the space--time reference elements $T_{std}^{st}=T_{std}\times I_{std}$ and $R_{std}^{st}=R_{std}\times I_{std}$ as: 
\begin{equation}
\tilde{\phi}_{k}(\xi,\eta,\tau)=\phi_{m_1(k)}(\xi, \eta) \cdot \gamma_{m_2(k)}(\tau) \quad \forall k =1,\ldots,N_{\phi}^{st}
\label{bf_phi_st}
\end{equation} and
\begin{equation}
\tilde{\psi}_{k}(\xi, \eta,\tau)=\psi_{d_1(k)}(\xi, \eta) \cdot \gamma_{d_2(k)}(\tau) \quad \forall k =1,\ldots,N_{\psi}^{st},
\label{bf_psi_st}
\end{equation}
having defined $N_\phi^{st}=N_\phi \cdot N_\gamma$, $N_\psi^{st}=N_\psi \cdot N_\gamma$ and the sub-indices $m_1(k)=1+(k-1)\mod N_{\phi}$,
$m_2(k)=1+\lfloor k/N_{\phi}\rfloor$,
$d_1(k)=1+(k-1)\mod N_{\psi}$ and 
$d_2(k)=1+\lfloor k/N_{\psi}\rfloor$.
\subsection{The space--time Galerkin finite elements}
\label{subparam}
The aim of this paragraph is to define the space--time Galerkin approximation spaces, using the space--time basis functions \eqref{bf_phi_st}-\eqref{bf_psi_st} defined on the reference elements.\\
Let us employ the superscript $ st $ for indicating the space--time extension, within the time interval $T^{n+1}$, of all the geometric entities introduced in paragraph \ref{stagg_grid}. For instance, we have ${T}_i^{st}$ and ${R}_j^{st}$ for the generic elements of the main and of the dual grid, respectively. Notice that, in the perspective of the new ALE method, these control volumes are stretched in time according to the local mesh velocity, but in the Eulerian limit, i.e. when the mesh velocity is equal to zero, they actually are the perpendicular space--time prisms ${T}_i^{st} = {T}_i \times T^{n+1}$ and ${R}_j^{st}={R}_j \times T^{n+1}$.\\
\begin{figure}[htbp] 
	\begin{center}
		\includegraphics[width=0.95\textwidth]{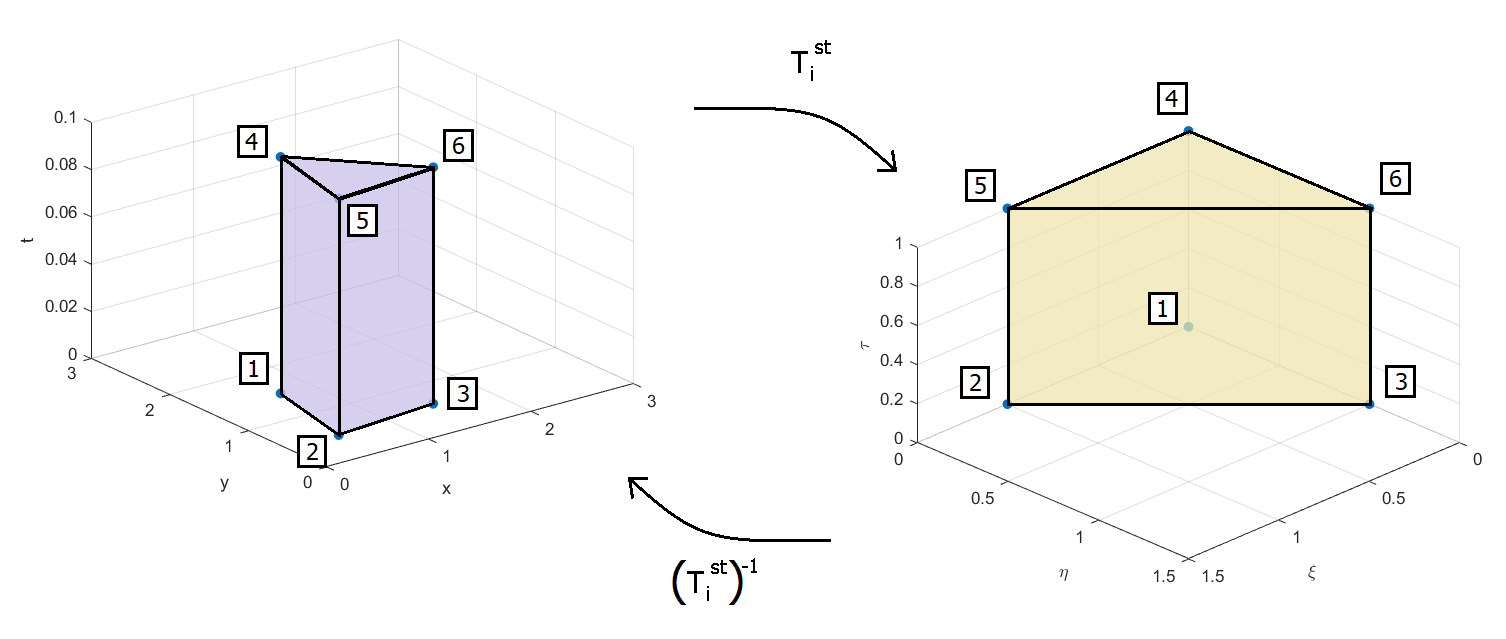}
		\caption{The transformations between the physical and the reference spaces,  for the primary space--time control volumes.}
		\label{transf_map}
	\end{center}
\end{figure}
Let $\{ (\hat{x}_k,\hat{y}_k,\hat{t}_k) \}_{k=1,\ldots,6}$ denote the coordinates of the vertices of $T_i^{st}$, as listed in Figure \ref{transf_map}; then the "subparametric" transformation from the reference space $\vec{\xi} = (\xi, \eta, \tau)$ to the physical space $\vec{x} = (x,y,t)$ is denoted by  $\left(\text{T}_i^{st}\right)^{-1}:T_{std}^{st} \rightarrow T_i^{st}$ and is defined by the action:
\begin{equation}
\left(\text{T}_i^{st}\right)^{-1}: \vec{\xi} \mapsto
\left\{
\begin{array}{l}
x=\sum\limits_{k=1}^6{\alpha_k(\vec{\xi}) \hat{x}_k} \\
y=\sum\limits_{k=1}^6{\alpha_k(\vec{\xi}) \hat{y}_k} \\
t=\sum\limits_{k=1}^6{\alpha_k(\vec{\xi}) \hat{t}_k}
\end{array}
\right.
,
\label{Tist_-1}
\end{equation}
where 
$\alpha_1(\vec{\xi}) = (1 - \xi - \eta)(1-\tau)$, $\alpha_2(\vec{\xi}) = \xi(1-\tau)$,
$\alpha_3(\vec{\xi}) = \eta(1-\tau)$,
$\alpha_4(\vec{\xi}) = (1 - \xi - \eta)\tau$,
$\alpha_5(\vec{\xi}) = \xi\tau$ and 
$\alpha_6(\vec{\xi}) = \eta\tau$.\\
Since $\hat{t}_1=\hat{t}_2=\hat{t}_3=t^n$ and $\hat{t}_4=\hat{t}_5=\hat{t}_6=t^{n+1}$, one can easily see that the third equation in \eqref{Tist_-1} reduces to:
$$t=(1-\tau)t^n+\tau t^{n+1},$$
therefore it is independent from the other two equations. After some calculations, it is possible to define also the inverse transformation $\text{T}_i^{st}:T_i^{st} \rightarrow T_{std}^{st}$ as:
\begin{equation}
\text{T}_i^{st}:\vec{x}\mapsto 
\left\{
\begin{array}{cl}
\xi   = & \frac{(y-C_y)B_x+(C_x-x)B_y}{A_yB_x - A_xB_y}\\
\eta  = & \frac{(x-C_y)A_y+(C_y-y)A_x}{A_yB_x - A_xB_y}\\
\tau  = & \frac{t-t^n}{t^{n+1}-t^{n}}
\end{array}
\right.
,
\label{Tist}
\end{equation}
where the coefficients in capital letters are calculated by the following expressions: 
\begin{equation}
\left\{
\begin{array}{cl}
A_x = & \tau(\hat{x}_5-\hat{x}_4)+(1-\tau)(\hat{x}_2-\hat{x}_1)\\
B_x = & \tau(\hat{x}_6-\hat{x}_4)+(1-\tau)(\hat{x}_3-\hat{x}_1)\\
C_x = & \tau \hat{x}_4 + (1-\tau) \hat{x}_1\\
A_y = & \tau(\hat{y}_5-\hat{y}_4)+(1-\tau)(\hat{y}_2-\hat{y}_1)\\
B_y = & \tau(\hat{y}_6-\hat{y}_4)+(1-\tau)(\hat{y}_3-\hat{y}_1)\\
C_y = & \tau \hat{y}_4 + (1-\tau) \hat{y}_1 
\end{array}
\right.
\label{cap_coeff}
\end{equation}
and by substituting in them the expression of $\tau$ given by the third equation of \eqref{Tist}.\\
For the staggered dual grid, the transformation between the reference space and the physical space is separately given  for the two sub-prisms that compose $R_j^{st}=T_{l(j),j}^{st}\cup T_{r(j),j}^{st} \: $. The transformations $\text{T}_{l(j),j}^{st}$ and $\text{T}_{r(j),j}^{st}$ are defined in the same way as in \eqref{Tist_-1}-\eqref{Tist}, with the only difference that the geometric coefficients refer to the vertices of the left and of the right sub-prisms, respectively. Of course, if $ j \in \mathcal{B}(\Omega_{h}) $, only one transformation is defined; by convention, we choose the left one: $ R_j^{st} = T_{l(j),j}^{st} \: $. \\
At this point, we are ready to define the proper approximation spaces for our problem. Let us define:
\begin{equation}
\mathcal{P}_{N,p_{\gamma}} = \{ \tilde{\varphi} : T_{std}^{st} \rightarrow \mathbb{R} \text{ s.t. } \tilde{\varphi}(\vec\xi) = \varphi(\xi,\eta)\cdot\theta(\tau), \varphi \in \mathbb{P}^N(T_{std}), \theta \in \mathbb{P}^{p_{\gamma}}(I_{std}) \} ,
\label{P_N,p_gamma}
\end{equation}
and, for every $ j \in [1, N_d]\setminus \mathcal{B}(\Omega_{h}) \: : $
\begin{equation}
\begin{array}{c}
\mathcal{D}_{N,p_{\gamma}}(R_j^{st}) = \{ w : R_{j}^{st} \rightarrow \mathbb{R} \text{ s.t. } 
w|_{T_{l(j),j}^{st}} = \tilde{\omega}_{1} \circ \text{T}_{l(j),j}^{st}, w|_{T_{r(j),j}^{st}} = \tilde{\omega}_{2} \circ \text{T}_{r(j),j}^{st}, \\ 
\tilde{\omega}_{1},\tilde{\omega}_{2} \in \mathcal{P}_{N,p_{\gamma}} \} .
\label{D_N,p_gamma_j}
\end{array}
\end{equation}
In formula \eqref{P_N,p_gamma}, $\mathbb{P}^N$ denotes the space of polynomials of degree less or equal to $N$, either in one or in two dimensions. Then our approximation spaces for the main and for the dual grid are:
\begin{equation}
V_{h}^{m} = \{ q_h : q_h|_{T_i^{st}} = \tilde{\varphi} \circ \text{T}_i^{st}, \tilde{\varphi} \in \mathcal{P}_{N,p_{\gamma}}, \forall i = 1,\ldots, N_e \} ,
\label{Vh_main}
\end{equation}
\begin{equation}
\begin{array}{clrc}
V_{h}^{d} = \{ w_h : & w_h|_{R_j^{st}}  \in \mathcal{D}_{N,p_{\gamma}}(R_j^{st}), & \forall j \in [1, N_d]\setminus \mathcal{B}(\Omega_{h}) , & \\
 & w_h|_{R_j^{st}} = \tilde{\varphi} \circ \text{T}_{l(j),j}^{st}, \tilde{\varphi} \in \mathcal{P}_{N,p_{\gamma}}, & \forall j \in \mathcal{B}(\Omega_{h})  \} & .
\end{array}
\label{Vh_dual}
\end{equation}
\\
For the discrete pressure $p_h$, we set $p_i=p_{h}|_{T_i^{st}}$ and we introduce the vector of the degrees of freedom of the numerical solution $\etah_i^{n+1}$, for every $i=1,\ldots,N_{e}$ and for every $n=0,\ldots,N_{t}-1 \: :$
\begin{equation}
p_i(\vec x)=\sum\limits_{l=1}^{N_{\phi}^{st}} \tilde{\phi}_l^{(i)}(\vec x)\hat{p}_{l,i}^{n+1}
=\tilde{\bphi}^{(i)}(\vec x) \cdot \etah_i^{n+1}
,
\label{eq:D_1}
\end{equation}
where of course $\tilde{\phi}_{l}^{(i)}=\tilde{\phi}_{l}\circ \text{T}_i^{st}$, for every $l=1,\ldots,N_{\phi}^{st}$.
\\
Similarly, for the discrete velocity vector field $\mathbf{v}_h$, we set $\mathbf{v}_j=\mathbf{v}_{h}|_{R_j^{st}}$ and we define the vectors $\tilde{\bpsi}^{(j)}$ and $\vvh_j^{n+1}$ such that:
\begin{equation}
\mathbf{v}_j(\vec x)=
\begin{cases}
\sum\limits_{l=1}^{N_\psi^{st}} \tilde{\psi}_l^{(j)}(\vec x) \hat{\mathbf{v}}_{l,j}^{n+1}
=\tilde{\bpsi}^{(j)}(\vec x)\cdot \vvh_j^{n+1} \quad & \forall j \in [1, N_d]\setminus \mathcal{B}(\Omega_{h})  \\
\sum\limits_{l=1}^{N_\phi^{st}} \tilde{\psi}_l^{(j)}(\vec x) \hat{\mathbf{v}}_{l,j}^{n+1}
=\tilde{\bpsi}^{(j)}(\vec x)\cdot \vvh_j^{n+1} \quad & \forall j \in \mathcal{B}(\Omega_{h})
\end{cases}
,
\label{eq:D_3}
\end{equation}
where we have distiguished the two cases for the internal and for the boundary elements.

\section{The semi-implicit Discontinuous Galerkin method}
\label{SIDG_method}
In this Section, a description with the key points of our staggered space--time DG method \cite{STINS2D,STINS2DTri,STINS3D} is outlined. First, the problem is arranged in the setting of the variational or weak formulation (Paragraph \ref{weak_f}); then it is written in the algebraic formulation (Paragraph \ref{algebr_f}); and, last, the solution algorithm is summarized in Paragraph \ref{p_corr}.
\subsection{The weak formulation}
\label{weak_f}
Let us consider equations \eqref{eq:CS_2_div}-\eqref{eq:CS_2_2_0_div}.
The Galerkin finite element problem is obtained by multiplying equation \eqref{eq:CS_2_div} by a generic function of the primary finite element space and then integrating over the space--time domain, and - in parallel - by multiplying equation \eqref{eq:CS_2_2_0_div} by a generic vectorial function of the dual finite element space and then integrating over the space--time domain. The problem reads as follows:\\
To find $\tilde{\mathbf{v}}_h \in  \left[V_{h}^{d}\right]^2$ and $p_h \in V_{h}^{m}$ such that:
\begin{equation}
\int\limits_{0}^{T} \int\limits_{\Omega} {q_h \left( \tilde{\nabla} \cdot \tilde{\mathbf{v}}_h \right)d\vec{x}}=0 \quad \forall q_h \in V_{h}^{m} ,
\label{eq:CS_4_FLRn_h}
\end{equation}
\begin{equation}
\int\limits_{0}^{T} \int\limits_{\Omega}{\mathbf{w}_h \cdot \left( \tilde{\nabla} \cdot \tilde{\F}_{cv,h} \right) d\vec{x}}+\int\limits_{0}^{T} \int\limits_{\Omega}{\mathbf{w}_h \cdot \left( \tilde{\nabla} \cdot \tilde{\F}_{p,h} \right) d\vec{x}}=\int\limits_{0}^{T} \int\limits_{\Omega}{\mathbf{w}_h \cdot \mathbf{S}_h d\vec{x}} \quad \forall \mathbf{w}_h \in  \left[V_{h}^{d}\right]^2.
\label{eq:CS_5_FLRn_h}
\end{equation}
Introducing the spatial and temporal discretizations from Section \ref{new_stag}, the problem given by equations \eqref{eq:CS_4_FLRn_h}-\eqref{eq:CS_5_FLRn_h} is equivalent to solve, for every $n=0, \ldots, N_{t}-1 $, the following systems of equations:
\begin{equation}
\int\limits_{T_i^{st}}{\tilde{\phi}_k^{(i)} \left( \tilde{\nabla} \cdot \tilde{\mathbf{v}}_h \right) d\vec{x}}=0,
\label{eq:CS_4n}
\end{equation}
for every $i=1,\ldots,N_e $, for every $k = 1, \ldots, N_\phi^{st}$, and:
\begin{equation}
\int\limits_{R_j^{st}}{\tilde{\psi}_k^{(j)}\left( \tilde{\nabla} \cdot \tilde{\F}_{cv,h} \right) d\vec{x}}+
\int\limits_{R_j^{st}}{\tilde{\psi}_k^{(j)} \left( \tilde{\nabla} \cdot \tilde{\F}_{p,h} \right) d\vec{x}}=
\int\limits_{R_j^{st}}{\tilde{\psi}_k^{(j)} \mathbf{S}_h d\vec{x}},
\label{eq:CS_5n}
\end{equation}
for every $j=1,\ldots,N_d$, for every $k = 1, \ldots, N_\psi^{st}$.\\
The integration by parts of equations \eqref{eq:CS_4n} leads to:
\begin{equation}
\oint\limits_{\partial T_i^{st}}{\tilde{\phi}_k^{(i)} \tilde{\mathbf{v}}_h \cdot \tilde{\mathbf{n}}_{i} \, ds dt}-\int\limits_{T_i^{st}}{\tilde{\nabla} \tilde{\phi}_k^{(i)} \cdot \tilde{\mathbf{v}}_h d\vec{x}}  =0,
\label{eq:CS_6n}
\end{equation}
where $\tilde{\mathbf{n}}_{i} = \left( \tilde n_t, \tilde n_x,\tilde n_y \right)^T $ is the outward pointing space--time unit normal vector referred to the control volume $ T_i^{st} $. Furthermore, by highlighting the contributions of the three sub-elements, we get:
\begin{equation}
\sum\limits_{j \in S_i}\left( \int\limits_{\Gamma_j^{st}}{\tilde{\phi}_k^{(i)} \tilde{\mathbf{v}}_j \cdot \tilde{\mathbf{n}}_{i,j}  ds dt}-
\int\limits_{T_{i,j}^{st}}{\tilde{\nabla} \tilde{\phi}_k^{(i)} \cdot \tilde{\mathbf{v}}_j d\vec{x}}  \right)=0,
\label{eq:CS_8n}
\end{equation}
where $\tilde{\mathbf{n}}_{i,j}=\tilde{\mathbf{n}}_{i}|_{\Gamma_j^{st}} \: $. \\
Now, let us focus on equations \eqref{eq:CS_5n}. The second term can be integrated directly, by interpreting the gradient of the pressure in the sense of distributions; thus:
\begin{equation}
\begin{array}{l}
\int\limits_{R_j^{st}}{\tilde{\psi}_k^{(j)}\left( \tilde{\nabla} \cdot \tilde{\F}_{cv,j} \right)  d\vec{x}}
+\int\limits_{T_{\ell(j),j}^{st}}{\tilde{\psi}_k^{(j)} \nabla p_{\ell(j)} d\vec{x}}
+\int\limits_{T_{r(j),j}^{st}}{\tilde{\psi}_k^{(j)} \nabla p_{r(j)} d\vec{x}} +\\
+\int\limits_{\Gamma_j^{st}}{\tilde{\psi}_k^{(j)} \left(p_{r(j)}-p_{\ell(j)}\right) \vec{n}_{l(j),j} ds dt}=
\int\limits_{R_j^{st}}{\tilde{\psi}_k^{(j)} \mathbf{S}_j d\vec{x}},
\end{array}
\label{eq:CS_9n}
\end{equation}
where $\tilde{\F}_{cv,j} = \tilde{\F}_{cv,h}|_{R_{j}^{st}} \: $; $\mathbf{S}_j = \mathbf{S}_h|_{R_j^{st}} \: $; and $ \vec{n}_{l(j),j} = (\tilde{n}_{x}, \tilde{n}_{y})^{T} $ is the spatial normal vector defined on the surface $\Gamma_j^{st}$ and directed from the left $ T_{\ell(j),j}^{st}$ to the right $ T_{r(j),j}^{st} \: $. Notice that the same result could also have been derived by forward and backward integration by parts, following \cite{BassiRebay}. Then, let us investigate the integration by parts of the first term in \eqref{eq:CS_9n}. This becomes:
\begin{equation}
\int \limits_{\partial R_j^{st}} \tilde{\psi}_k^{(j)} \tilde{\F}_{cv,j} \cdot \tilde{\mathbf{n}}_{j} ds dt 
- \int \limits_{R_j^{st}} \tilde{\F}_{cv,j} \cdot \tilde \nabla \tilde{\psi}_k^{(j)} d\vec{x} ,
\label{eqn.GaussPDE}
\end{equation}
where $\tilde{\mathbf{n}}_{j} = \left( \tilde n_t, \tilde n_x,\tilde n_y \right)^T $ is the outward pointing space--time unit normal vector referred to $ R_j^{st} $. Specifically, the boundary of $ R_j^{st} $ is:
\begin{equation}
\partial R_j^{st} = \left( \bigcup \limits_{ll \in E_{j}} \Lambda_{ll}^{st} \right) 
\cup R_j^{n} \cup R_j^{n+1},
\label{dCi}
\end{equation}
where $ R_j^{n} $ and $ R_j^{n+1} $ are the cell configurations at times $t^n$ and $t^{n+1}$, respectively. Using $\tilde{\mathbf{n}}_{j}|_{R_{j}^{n+1}} = \left(1,0,0\right)^{T}$ and $\tilde{\mathbf{n}}_{j}|_{R_{j}^{n}} = \left(-1,0,0\right)^{T} $, expression \eqref{eqn.GaussPDE} is modified into:
\begin{equation}
\begin{array}{l}
\int \limits_{R_j^{n+1}} \tilde{\psi}_k^{(j)} \mathbf{v}_j(x,y,t^{n+1}) dxdy 
- \int \limits_{R_j^{n}} \tilde{\psi}_k^{(j)} \mathbf{v}_j(x,y,t^{n}) dxdy + \\
+ \sum_{ll \in E_{j}} \left(\int \limits_{\Lambda_{ll}^{st}}{\tilde{\psi}_k^{(j)} \tilde{\F}_{cv,j} \cdot \tilde{\mathbf{n}}_{j,ll} ds dt} \right)
- \int \limits_{R_j^{st}} \tilde{\F}_{cv,j} \cdot \tilde \nabla \tilde{\psi}_k^{(j)} d\vec{x},
\end{array}
\label{eqn.GaussPDE1}
\end{equation}
where $\tilde{\mathbf{n}}_{j,ll}=\tilde{\mathbf{n}}_{j}|_{\Lambda_{ll}^{st}} \:$. The third term in \eqref{eqn.GaussPDE1} is replaced with a numerical flux function in space--time normal direction, denoted by $\tilde{\vec{\mathcal{G}}} \cdot \tilde{\mathbf{n}}_{j,ll} \: $, which takes into account the information of the two adjacent elements $(\mathbf{v}_j^-,\nabla \mathbf{v}_j^-)$ and $(\mathbf{v}_j^+,\nabla \mathbf{v}_j^+)$.\footnote{Or, the values of the boundary conditions, if $j \in \mathcal{B}(\Omega_{h})$ and $\Lambda_{ll} \subset \partial{\Omega}$.} As an approximate Riemann solver, we rely on a very robust Rusanov flux \cite{Rusanov:1962}, and, following \cite{DumbserNSE}, it includes both the convective and the viscous terms:
\begin{equation}
\tilde{\vec{\mathcal{G}}} \cdot \tilde{\mathbf{n}}_{j,ll}  =   
\frac{1}{2} \left( \tilde{\F}_{cv,j}^{+} + \tilde{\F}_{cv,j}^{-}  \right) \cdot \tilde{\mathbf{n}}_{j,ll} - 
\frac{1}{2} \left(|s_{\max}| + 2\eta |s_{\max}^{\nu}| \right) \left( \mathbf{v}_j^+ - \mathbf{v}_j^- \right).  
\label{eqn.rusanov_} 
\end{equation} 
We propose this Riemann solver because it has already been applied in the ALE context \cite{Lagrange2D,Lagrange3D,LagrangeDG}, however the following definitions do not involve any mesh velocity, because we are only interested in the Eulerian limit. Thus, $s_{\max}$ represents the maximum eigenvalue of the Jacobian matrix, in normal direction, coming from the purely convective transport operator:
\begin{equation} 
\mathbf{A}_{\mathbf{n}}
=\left(\sqrt{\tilde n_x^2 + \tilde n_y^2 
}\right)\mathbf{A} \cdot \mathbf{n}_{j,ll} , \qquad    
\mathbf{A}=\frac{\partial \mathbf{F}_{cv}(\mathbf{v},\nabla \mathbf{v}
	)}{\partial \mathbf{v}} , \qquad
\mathbf{n}_{j,ll} = \frac{(\tilde n_x, \tilde n_y
	)^T}{\sqrt{\tilde n_x^2 + \tilde n_y^2
}}\: ;
\label{eqn.An}
\end{equation} 
and $s_{\max}^{\nu}$ is the maximum eigenvalue of the Jacobian matrix of the viscous operator, in normal direction, given by:
\begin{equation}
\mathbf{D}_{\mathbf{n}}=\frac{\partial \mathbf{F}_{cv}(\mathbf{v},\nabla \mathbf{v})}{\partial (\nabla \mathbf{v}
	)} \cdot \mathbf{n}_{j,ll} \: .
\end{equation}
Finally, according to \cite{DumbserNSE,MunzDiffusionFlux,HidalgoDumbser}, the factor $\eta$ is estimated as:
\begin{equation}
\eta = \frac{2N+1}{h_\nu \sqrt{\frac{\pi}{2}}},
\label{eqn.eta}
\end{equation}
where the characteristic size $h_\nu$ is given by $h^{+} + h^{-}$, being $h^-$ and $h^+$ the radii of the inscribed circles in $R_j^{n}$ and in the neighbouring element, respectively.
\subsection{The algebraic formulation}
\label{algebr_f}
Let us now derive the algebraic formulation of the Galerkin problem. If we substitute the expressions of the unknowns $p_{i}$ and $\mathbf{v}_{j}$ in terms of the basis functions, see Formulae \eqref{eq:D_1}-\eqref{eq:D_3}, we can write our problem, for every $n=0,\ldots,N_t -1$, as:
\begin{equation}
\sum\limits_{j \in S_i}\left(\int\limits_{\Gamma_j^{st}}{\tilde{\phi}_k^{(i)}\tilde{\psi}_l^{(j)} \tilde{\mathbf{n}}_{i,j} ds dt} \cdot \hat{\tilde{\mathbf{v}}}_{l,j}^{n+1}-
\int\limits_{T_{i,j}^{st}}{\tilde{\nabla} \tilde{\phi}_k^{(i)}\tilde{\psi}_l^{(j)} d\vec{x}} \cdot \hat{\tilde{\mathbf{v}}}_{l,j}^{n+1} \right)=0,
\label{eq:CS_10}
\end{equation}
for every $i=1,\ldots,N_e$, for every $k = 1, \ldots, N_\phi^{st}$, and:
\begin{equation}
\begin{array}{cl}
\int \limits_{R_j^{n+1}}&\tilde{\psi}_k^{(j)} \tilde{\psi}_l^{(j)} \vvh_{l,j}^{n+1} ds 
- \int \limits_{R_j^{n}} \tilde{\psi}_k^{(j)} \tilde{\psi}_l^{(j)} \vvh_{l,j}^{n} ds
+ \sum_{ll \in E_{j}} \left(\int \limits_{\Lambda_{ll}^{st}} \tilde{\psi}_k^{(j)} \tilde{\vec{\mathcal{G}}} \cdot \tilde{\mathbf{n}}_{j,ll} \, ds dt \right) + \\
&- \int \limits_{R_j^{st}} \tilde{\F}_{cv,j} \cdot \tilde \nabla \tilde{\psi}_k^{(j)} d\vec{x}
+\int\limits_{T_{\ell(j),j}^{st}}{\tilde{\psi}_k^{(j)} \nabla  \tilde{\phi}_{l}^{(\ell(j))} d\vec{x}}  \cdot \etah_{l,\ell(j)}^{n+1} + \\
&+\int\limits_{T_{r(j),j}^{st}}{\tilde{\psi}_k^{(j)} \nabla \tilde{\phi}_{l}^{(r(j))} d\vec{x}} \cdot \etah_{l,r(j)}^{n+1}  
+\int\limits_{\Gamma_j^{st}}{\tilde{\psi}_k^{(j)} \tilde{\phi}_{l}^{(r(j))} \vec{n}_{l(j),j} ds dt} \cdot \etah_{l,r(j)}^{n+1} +  \\
&-\int\limits_{\Gamma_j^{st}}{\tilde{\psi}_k^{(j)} \tilde{\phi}_{l}^{(\ell(j))} \vec{n}_{l(j),j} ds dt} \cdot \etah_{l,\ell(j)}^{n+1}
=\int\limits_{R_j^{st}}{\tilde{\psi}_k^{(j)} \mathbf{S}_{j} d\vec{x}},
\end{array}
\label{eq:CS_11_1}
\end{equation}
for every $j=1,\ldots,N_d$, for every $k = 1, \ldots, N_\psi^{st}$, where we have used the standard summation convention for the repeated index $l$.\\
Equations \eqref{eq:CS_10}-\eqref{eq:CS_11_1} are written in a compact matrix form as:
\begin{equation}
\sum\limits_{j \in S_{i}} \vec{\D}_{i,j} \vvh_j^{n+1} = \mathbf{0} 
\label{eq:CS_12},
\end{equation}
\begin{equation}
\left(\Mpsi_j^+ - \Mpsi_j^\circ \right) \vvh_j^{n+1} - \Mpsi_j^-\vvh_j^{n} + \sum\limits_{ll \in E_{j}} \left( \tilde{\vec{\Nn}}_{ll,j} \right) - \tilde{\vec{\Ff}}_{j} - \vec{\LM}_j \etah_{\ell(j)}^{n+1}+ \vec{\RM}_j \etah_{r(j)}^{n+1} =\tilde{\vec{\Ss}}_j \: , 
\label{eq:CS_12_1}
\end{equation}
where:
\begin{equation}
\vec{\D}_{i,j} =\int\limits_{\Gamma_j^{st}}{\tilde{\phi}_k^{(i)}\tilde{\psi}_l^{(j)} \vec{n}_{i,j} ds dt}-\int\limits_{T_{i,j}^{st}}{\nabla \tilde{\phi}_k^{(i)}\tilde{\psi}_l^{(j)}d\vec{x}}, 
\label{eq:MD_3} 
\end{equation}
\begin{eqnarray}
\Mpsi_j^+ &=& \int\limits_{R_j^{n+1}}{\tilde{\psi}_k^{(j)}(x,y,t^{n+1})\tilde{\psi}_l^{(j)}(x,y,t^{n+1})  dx dy}, \label{eq:MD_2} \\
\Mpsi_j^- &=& \int\limits_{R_j^{n}}{\tilde{\psi}_k^{(j)}(x,y,t^{n})\tilde{\psi}_l^{(j)}(x,y,t^{n})  dx dy}, \label{eq:MD_2_1} \\
\Mpsi_j^\circ &=& \int\limits_{R_j^{st}}{\diff{\tilde{\psi}_k^{(j)}}{t} \tilde{\psi}_l^{(j)} d\vec{x}}, \label{eq:MD_2_2} 
\end{eqnarray}
\begin{eqnarray}
\tilde{\vec{\Nn}}_{ll,j}&=&
\int \limits_{\Lambda_{ll}^{st}} \tilde{\psi}_k^{(j)} \tilde{\vec{\mathcal{G}}} \cdot \tilde{\mathbf{n}}_{j,ll} ds dt, \label{eq:MD_34} \\
\tilde{\vec{\Ff}}_{j} &=& \int \limits_{R_j^{st}} \tilde{\F}_{cv,j} \cdot \tilde{\nabla} \tilde{\psi}_k^{(j)} d\vec{x}, 
\label{eq:MD_35} 
\end{eqnarray}
\begin{eqnarray}
\vec{\LM}_{j} &=& \int\limits_{\Gamma_j^{st}}{\tilde{\psi}_k^{(j)} \tilde{\phi}_{l}^{(\ell(j))}\vec{n}_{l(j),j} ds dt}-\int\limits_{T_{\ell(j),j}^{st}}{\tilde{\psi}_k^{(j)} \nabla \tilde{\phi}_{l}^{(\ell(j))} d\vec{x}},
\label{eq:MD_5} \\
\vec{\RM}_{j} &=& \int\limits_{\Gamma_j^{st}}{\tilde{\psi}_k^{(j)} \tilde{\phi}_{l}^{(r(j))}\vec{n}_{l(j),j} ds dt}+\int\limits_{T_{r(j),j}^{st}}{\tilde{\psi}_k^{(j)} \nabla \tilde{\phi}_{l}^{(r(j))} d\vec{x}},
\label{eq:MD_4}
\end{eqnarray}
\begin{equation}
\tilde{\vec{\Ss}}_j = \int\limits_{R_j^{st}}{\tilde{\psi}_k^{(j)} \mathbf{S}_{j} d\vec{x}}. 
\label{eq:MD_5_2}
\end{equation}
In order to simplify the notation, we introduce the new matrix $\Q_{i,j} \:$:
\begin{equation}
\Q_{i,j}=\int\limits_{T_{i,j}^{st}}{\tilde{\psi}_k^{(j)} \nabla \tilde{\phi}_{l}^{(i)} d\vec{x}}-\int\limits_{\Gamma_j^{st}}{\tilde{\psi}_k^{(j)} \tilde{\phi}_{l}^{(i)}\sigma_{i,j} \vec{n}_{l(j),j} ds dt},
\label{eq:MD_6}
\end{equation}
where $\sigma_{i,j}$ is defined by:
\begin{equation}
\sigma_{i,j}=\frac{r(j)-2i+\ell(j)}{r(j)-\ell(j)}.
\label{eq:SD_1}
\end{equation}
In this way, $\Q_{\ell(j),j}=-\vec{\LM}_j$ and $\Q_{r(j),j}=\vec{\RM}_j$. \\
In order to further ease the notation, we use the abbreviation $\Mpsi_j = \Mpsi_j^+ - \Mpsi_j^\circ$, therefore 
equations \eqref{eq:CS_12}-\eqref{eq:CS_12_1} are rewritten as: 
\begin{equation}
\sum\limits_{j \in S_i}\vec{\D}_{i,j}\vvh_j ^{{n+1}}=\mathbf{0}, 
\label{eq:CS_15}
\end{equation}
\begin{equation}
\Mpsi_j\vvh_j^{n+1} - \Mpsi_j^-\vvh_j^{n} + \sum\limits_{ll \in E_{j}} \left( \tilde{\vec{\Nn}}_{ll,j} \right) - \tilde{\vec{\Ff}}_{j} + \Q_{\ell(j),j} \etah_{\ell(j)} ^{{ n+1}} + \Q_{r(j),j} \etah_{r(j)}^{{n+1}} = \tilde{\vec{\Ss}}_j \: 
. 
\label{eq:CS_16}
\end{equation}
\begin{rem}
	\label{rem_ALE}
	As it was discussed before in Paragraph \ref{spatial_bf}, we have adopted an efficient memory-saving, quadrature-free implementation of the matrices of our problem. We remark that, in the perspective of an ALE algorithm, each of the matrices defined above must be recomputed at every time step, in view of different mesh configurations between times $t^n$ and $t^{n+1}$. In order to reduce the cost of this operation, we carry out this task with the aid of a \textit{universal} matrix in the following way.\\
	Let $I=a^{K}_{k,l}$ represent one element in the generic local matrix $\left[\vec{\mathcal{A}}_K\right]_{k,l} \: $. Let us suppose that, after the appropriate transformation $\text{T}: K \rightarrow K_{std}$ from the physical space $\vec{x}$ to the reference space $\vec{\xi}$, the integrand function has got a factor $g_{K}^{n+1}(\vec{\xi})$ which depends on the specific control volume $K$ (e.g. the outward normal vector, or the determinant of the Jacobian matrix of the inverse transformation $\mathbf{J}(\text{T}^{-1})$). Notice that, because of the theorem of the change of variables, the determinant of $\mathbf{J}(\text{T}^{-1})$ is always required in the computation of $I$. Let us provide the approximation of $g_{K}^{n+1}$ as a function in $\mathcal{P}_{N,p_{\gamma}}$:
	$$
	g_{K}^{n+1}(\vec{\xi})=\sum_{m=1}^{N_{\phi}^{st}}{\tilde{\phi}_{m}(\vec{\xi})\cdot \hat{g}_{K,m}^{n+1}} \: . 
	$$
	Thus, we are able to split the computation of $I$ in the scalar product of a \textit{universal} factor and a \textit{mesh-dependent} factor: $I = \tilde{I}_{k,l,m}\cdot\hat{g}_{K,m}^{n+1}$. Looking that from a broader view, this means that we are updating the generic tensor $\vec{\mathcal{A}}$ by means of a cheap tensor-reduction operation, which allows us to save a lot of time in the simulation's process.\\
	At the same time, for each tensor $\vec{\mathcal{A}}$, only its universal version $\tilde{\vec{I}}$ must be stored in the RAM memory. Since its size does not depend on the mesh $\mathcal{T}_{h}$, then, given $N$ and $p_{\gamma}$, only a relatively small, constant piece of RAM memory is occupied by this tensor. Incidentally, we also remark that the universal tensors associated to all the matrices of our problem are ready-to-use for every new simulation, because they are saved in binary files, for the most common combinations of $\left(N,p_{\gamma}\right)$. For their computation, we have used a quadrature-free approach in a program employing symbolic computation. Later, in Section \ref{N_tests} concerning numerical test problems we will clearly show, e.g. via Table \ref{TG_RAM_memory} and Figure \ref{TG_Comp_Times}, that the new algorithm proposed in this paper requires much less RAM memory and less CPU time than the original staggered semi-implicit DG scheme proposed in \cite{STINS2D,STINS2DTri,STINS3D}.  
\end{rem}
Let us now specify the treatment of the convective and viscous fluxes\footnote{For simplicity, the Lagrangian contribute due to the temporal derivative is neglected in this article, since we only operate in the Eulerian limit of the ALE method.} in vectors $\tilde{\vec{\Nn}}_{ll,j}$ and $\tilde{\vec{\Ff}}_{j}$. We decide not to take both terms implicitly, otherwise we would have to solve a nonlinearity in the final system. Following \cite{STINS2D,STINS2DTri}, we discretize the velocity $\mathbf{v}_h$ explicitly and its gradient $\nabla \mathbf{v}_h$ implicitly. In this way, we avoid additional restrictions on the maximum time step given by the viscous terms. Indeed, the above numerical scheme requires that the time step satisfies the CFL-type stability condition for DG schemes:
\begin{equation}
\Delta t^{n+1}= \frac{1}{2} \frac{\textnormal{CFL}}{2N+1}\cdot \frac{h_{min}}{|s_{max}|},
\label{eq:CFLC}
\end{equation}
where $h_{min}$ is the smallest incircle radius among all the triangles $T_{i}^{n}$, with $i=1,\ldots,N_{e}$; $\textnormal{CFL}<1$; and $s_{max}$ is the maximum eigenvalue of the convective flux tensor, %that is $2|\mathbf{v}\cdot\mathbf{n}|$.
defined in Formula \eqref{eqn.An}.
\subsection{Pressure correction formulation and final algorithm}
\label{p_corr}
At this level, we are ready to outline the final solving algorithm of our problem. The basic idea is to introduce the classical SIMPLE algorithm \cite{patankarspalding} used for the resolution of the incompressible Navier-Stokes equations. According to \cite{STINS2DTri}, we call this method a coupled space--time pressure correction algorithm.\\
If we define an intermediate solution $\hat{\mathbf{v}}_{j}^{n+\frac{1}{2}}$, in place of $\hat{\mathbf{v}}_{j}^{n+1}$, for the system of equations \eqref{eq:CS_16}, by taking the pressure explicitly, rather than implicitly, then we can readily compute it by solving the system:
\begin{equation}
\Mpsi_j \hat{\mathbf{v}}_j^{n+\frac{1}{2}} = \Mpsi_j^{-}\vvh_j^{n} -  \sum\limits_{ll \in E_{j}} \left( \tilde{\vec{\Nn}}_{ll,j}  \right) + \tilde{\vec{\Ff}}_{j} - \Q_{\ell(j),j} \etah_{\ell(j)} ^{n} - \Q_{r(j),j} \etah_{r(j)}^{n} + \tilde{\vec{\Ss}}_j \: .
\label{eq:A1_mod}
\end{equation}
From the comparison between equations \eqref{eq:CS_16} and \eqref{eq:A1_mod}, we can restore the system of equations for the velocity degrees of freedom $\hat{\mathbf{v}}_{j}^{n+1} \: $:
\begin{equation}
\Mpsi_j	\vvh_j^{n+1}-\Mpsi_j \vvh_j^{n+\frac{1}{2}}
+\Q_{\ell(j),j}\Delta\etah_{\ell(j)}^{n+1}
+\Q_{r(j),j} \Delta\etah_{r(j)}^{n+1}
=\mathbf{0},
\label{press_corr}
\end{equation}
where we have defined the pressure variation $\Delta \hat{\mathbf{p}}_{i}^{n+1}=\hat{\mathbf{p}}_{i}^{n+1} - \hat{\mathbf{p}}_{i}^{n}$. But, before doing that, the pressure correction formulation provides for the computation of these vectors $\Delta \hat{\mathbf{p}}_{i}^{n+1}$ from the so called pressure Poisson equation. This comes from the formal substitution of the velocity degrees of freedom into \eqref{eq:CS_15}, see also \cite{CasulliCheng1992,DumbserCasulli}:
\begin{equation}
\sum\limits_{j\in S_i} \vec{\D}_{i,j}\Mpsi_j ^{-1}\Q_{\ell(j),j} \Delta\etah_{\ell(j)}^{n+1}
+\sum\limits_{j\in S_i}\vec{\D}_{i,j}\Mpsi_j^{-1}\Q_{r(j),j} \Delta\etah_{r(j)}^{n+1}
=\sum\limits_{j \in S_i} \vec{\D}_{i,j} \hat{\mathbf{v}}_{j}^{n+\frac{1}{2}}.
\label{eq:CS_19}
\end{equation}
We underline one of the main advantages of the \textit{staggered} approach: the linear system for the pressure correction is made of only four non-zero blocks, because for each primary element, there is the central contribution and only three other contributions of the neighbours; therefore this system is very sparse and quite inexpensive to be solved. \\
On the other hand, the fact that the above algorithm contains a sort of operator-splitting or \textit{fractional step} procedure is the cause of the low-order in time (indeed we would get only a first order method). In order to overcome this problem, we use a simple \textit{Picard iteration}. This approach is inspired by the local space--time Galerkin predictor method proposed for the high-order time discretization of $P_NP_M$ schemes in \cite{DumbserNSE,Dumbser2008}.
For every $n=0, \dots, N_{t}-1$, the final algorithm can be summarized as follows.
\begin{enumerate}
	\item Initialize $\hat{\mathbf{v}}_{j}^{n+1,0}$ copying the information from the previous time step (or using the initial conditions, if $n=0$); and initialize $\hat{\mathbf{p}}_{i}^{{n+1,0}}$ to a constant initial guess (for instance, $1$).\footnote{Alternatively, we could also choose the extrapolation of $p_{h}^{n}$ in $T^{n+1}$.}
	\item For every $k=0, \ldots, N_{pic}=p_{\gamma} \: $:
	\begin{enumerate}
		\item compute $\hat{\mathbf{v}}_{j}^{n+1,k+\frac{1}{2}}$:
		\begin{equation}
		\begin{array}{cl}
		\hat{\mathbf{v}}_j^{n+1,k+\frac{1}{2}} &= \Mpsi_j^{-1}\Mpsi_j^{-}\vvh_j^{n} 
		-  \sum\limits_{ll \in E_{j}} \Mpsi_j^{-1} \tilde{\vec{\Nn}}_{ll,j} \left(\mathbf{v}_h^{n+1,k},\nabla \mathbf{v}_h^{n+1,k+\frac{1}{2}} \right) + \\
		&+ \Mpsi_j^{-1} \tilde{\vec{\Ff}}_{j} \left(\mathbf{v}_h^{n+1,k},\nabla \mathbf{v}_h^{n+1,k+\frac{1}{2}} \right) + \\
		&- \Mpsi_j^{-1} \Q_{\ell(j),j} \etah_{\ell(j)}^{{n+1,k}}- \Mpsi_j^{-1} \Q_{r(j),j} \etah_{r(j)}^{{n+1,k}} + \Mpsi_j^{-1} \tilde{\vec{\Ss}}_j \: ;
		\label{eq:A1_mod_k}
		\end{array}
		\end{equation}
		\item compute $\Delta \etah_{i}^{n+1,k+1}$:
		\begin{equation}
		\begin{array}{l}
		\sum\limits_{j\in S_i} \vec{\D}_{i,j}\Mpsi_j ^{-1}\Q_{\ell(j),j}  \Delta\etah_{\ell(j)}^{n+1,k+1}
		+\sum\limits_{j\in S_i}\vec{\D}_{i,j}\Mpsi_j^{-1}\Q_{r(j),j} \Delta \etah_{r(j)}^{n+1,k+1} 
		 =  \\
		=\sum\limits_{j \in S_i} \vec{\D}_{i,j} \hat{\mathbf{v}}_j^{n+1,k+\frac{1}{2}};
		\end{array}
		\label{eq:A2_mod_k}
		\end{equation}
		\item update $\hat{\mathbf{v}}_{j}^{n+1,k+1}$ and $\etah_{i}^{n+1,k+1}$:
		\begin{equation}
		%\begin{array}{c}
		\vvh_j^{n+1,k+1}= \hat{\mathbf{v}}_j^{n+1,k+\frac{1}{2}}
		- \Mpsi_j^{-1} \Q_{\ell(j),j} \Delta \etah_{\ell(j)} ^{{ n+1, k+1}}
		- \Mpsi_j^{-1} \Q_{r(j),j} \Delta\etah_{r(j)}^{{n+1,k+1}} ;
		%\end{array}
		\label{eq:A3_mod}
		\end{equation}
		\begin{equation}
		\etah_{i}^{n+1,k+1}=\etah_{i}^{n+1,k}+\Delta\etah_{i}^{n+1,k+1}.
		\end{equation}
	\end{enumerate}	
	\item Set $\hat{\mathbf{v}}_{j}^{n+1}=\hat{\mathbf{v}}_{j}^{n+1,k+1}$	and  $\etah_{i}^{{n+1}}=\etah_{i}^{{n+1,k+1}}$. 
\end{enumerate}
In \eqref{eq:A2_mod_k} and \eqref{eq:A3_mod}, we have defined the pressure correction for the $k$-th Picard iteration $\Delta \hat{\mathbf{p}}_i^{{ n+1, k+1}} = \hat{\mathbf{p}}_i^{{ n+1, k+1}}-\hat{\mathbf{p}}_i^{{ n+1, k}}$.

\section{Numerical Tests}
\label{N_tests}
In this Section, we present the numerical tests that we have performed in order to validate our new numerical method. Several benchmarks are considered, from the 1D elliptic and parabolic equations to the full model of the incompressible Navier-Stokes equations. For each test, we have compared the performances obtained by the new ALE implementation of the method, in the limit $V=0$, together with the Eulerian implementation presented in \cite{STINS2DTri}, showing how the new velocity-pressure elements allow to build a very efficient ALE algorithm.
\subsection{The Couette flow}
\label{Couette}
\begin{figure}[htbp] 
	\begin{center}
		\includegraphics[width=0.75\textwidth]{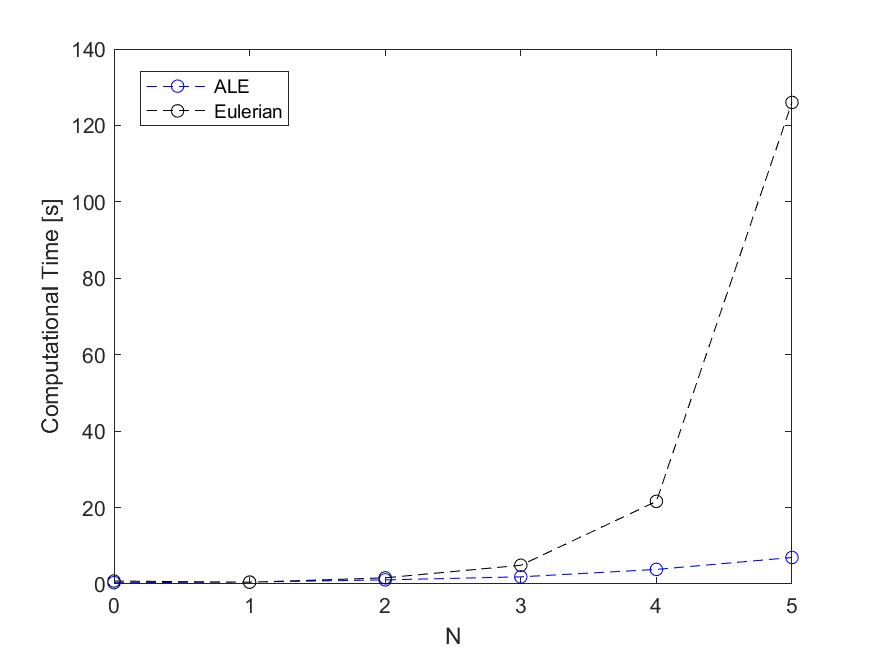}	
		\caption{Steady Couette flow test: comparison between the computational times of the simulations for the ALE and the Eulerian implementations.}
		\label{TEST26LV_Times}
	\end{center}
\end{figure}

\begin{figure}[htbp] 
	\begin{center}
		\includegraphics[width=0.75\textwidth]{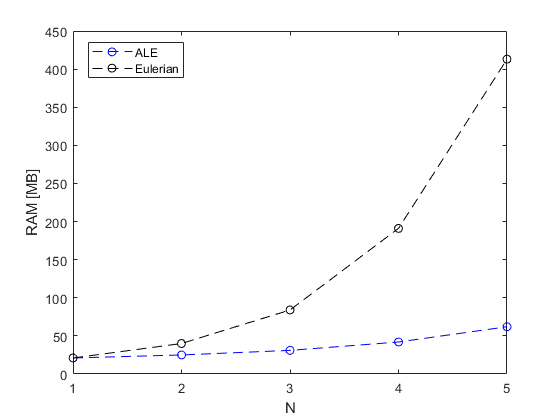}	
		\caption{Steady Couette flow test: comparison between the memory usage for the ALE and the Eulerian implementations.}
		\label{TEST26LV_Memory}
	\end{center}
\end{figure}

In fluid dynamics, the Couette flow is a classical example of shear-driven flow for viscous fluids. Let us consider two infinite, parallel plates, separated by a distance $L=1$, and, between them, a viscous fluid with kinematic viscosity coefficient $\nu=0.01$. Here, we assume that the flow is uni-directional, that is $v = 0$, and that the steady state is already reached. 
A schematic representation of this test can be found in Figure \ref{Couette_description}: 
\begin{figure}[htbp] 
	\begin{center}
		\includegraphics[width=0.40\textwidth]{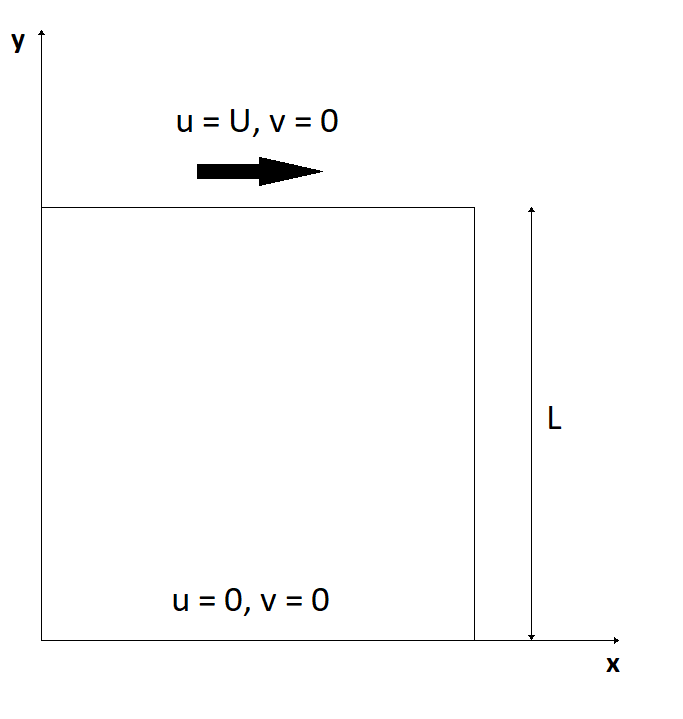}
		\caption{The Couette flow test: a schematic representation.}
		\label{Couette_description}
	\end{center}
\end{figure}
a no-slip (or "wall") boundary condition, expressed by $u=0$, is imposed in the lower plate, located at $y=0$; a non homogeneous Dirichlet's boundary condition on the velocity, expressed by $u=U$, with $U=1$, is imposed on the upper plate; finally, periodic boundary conditions are used at the left and at the right sides of the domain. In these conditions, the steady regime that originates is a laminar flow driven by the viscous drag force. From Equations \eqref{eq:CS_2} and \eqref{eq:CS_2_2_0}, if the gravity term is not considered, and since no pressure gradient is present, the only non-trivial equation remains the $x$-momentum equation:
\begin{equation}
- \nu \frac{\partial ^2 u}{\partial y ^2} = 0 .
\label{SteadyCouette}
\end{equation}
The exact solution of this problem is given by:
$$ u(y) = U\frac{y}{L} , $$
and, according to Newton's law, the tangential stress is a constant given by:
$$ \tau = \mu \frac{\partial u}{\partial y} = \mu \frac{U}{L} \quad \forall y \in (0,L) .$$
\begin{figure}[htbp] 
	\begin{center}
		\includegraphics[width=0.45\textwidth]{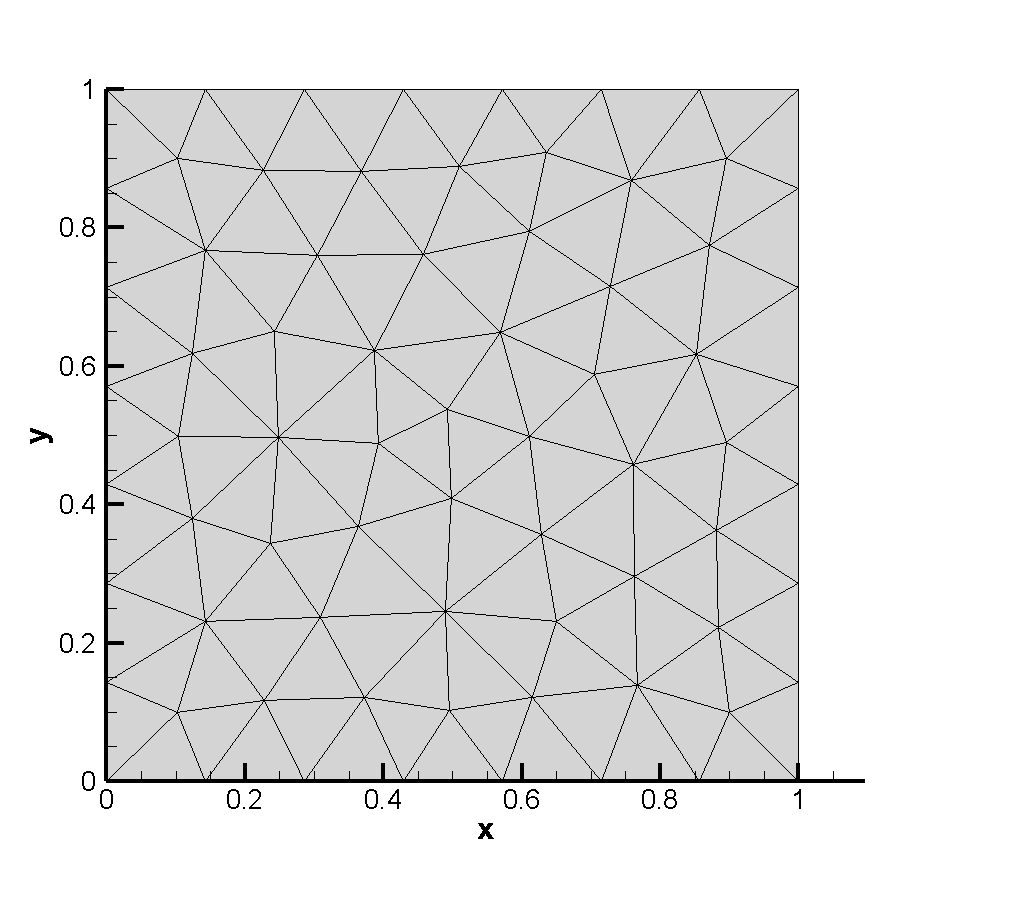}
		\caption{The unstructured mesh with $N_e=116$.}
		\label{mesh116}
	\end{center}
\end{figure}
In order to check that our numerical method is able to reproduce the exact solution, we set up a numerical test where a coarse unstructured mesh with $N_{e}=116$ elements (Figure \ref{mesh116}) was used and an increasing number $N$ of polynomial degrees from $0$ to $5$ was employed. The results in Table \ref{TEST26LV} confirm that the linear solution is correctly reproduced up to machine precision, for all degrees $N\geq1$.\\
\begin{table}[htbp] 
	\caption{Steady Couette flow test: errors in the $ L^2 $ norms for the velocity, after $20$ timesteps, obtained by $p_{\gamma}=0$ and increasing polynomial degree $N$. In column (b), the Eulerian implementation was used.} 
	\begin{center}
		\scriptsize 
		\subfloat[ALE, $V=0$]{
			\label{TEST26LV_moving_Ana}
			\begin{tabular}{|l|c|} 
				\hline 
				$ N $  &$E_{2}^{v}$ 	\\ \hline 
				0   & 2.35E-02 \\ \hline 
				1   & 8.51E-14 \\ \hline 
				2   & 2.55E-13  \\ \hline 
				3   & 3.05E-13 \\ \hline
				4   & 2.70E-13 \\ \hline
				5   & 2.07E-13 \\ \hline
			\end{tabular} 
		}
		\subfloat[Eulerian]{
			\label{TEST26LV_fixed}
			\begin{tabular}{|c|} 
				\hline 
				$E_{2}^{v}$ 	\\ \hline 
				2.35E-02 \\ \hline 
				1.97E-14 \\ \hline 
				5.82E-14 \\ \hline 
				7.02E-14 \\ \hline
				2.26E-13 \\ \hline 
				9.13E-12 \\ \hline
			\end{tabular} 
		}
		\label{TEST26LV} 
	\end{center}
\end{table} 
A comparison in terms of computational times between the new quadrature-free ALE implementation with the new discrete velocity space and the standard Eulerian approach using a simple tensor product basis on the unit square for the velocity space is reported in Figure \ref{TEST26LV_Times}. The advantage of using the new implementation is evident for the highest polynomial degree $N=5$, where we have $126$ seconds for the old 
and only $7$ seconds for the new implementation, which means that only $5.5 \%$ of the time is required by the new scheme. \\
For the same test, also a comparison in terms of the RAM memory was carried out: the results are reported in Figure \ref{TEST26LV_Memory}. We can observe that the memory required by the new algorithm is only slightly increasing with $N$ and it is significantly smaller than the old algorithm, even in a case of small number of elements, $N_{e}=116$. For example, in the simulations with $N=5$, the new algorithm required $62 MB$ of RAM, versus $413 MB$ of the old one, which means only $15\%$.

\begin{figure}[htbp] 
	\begin{center}
		\includegraphics[width=0.75\textwidth]{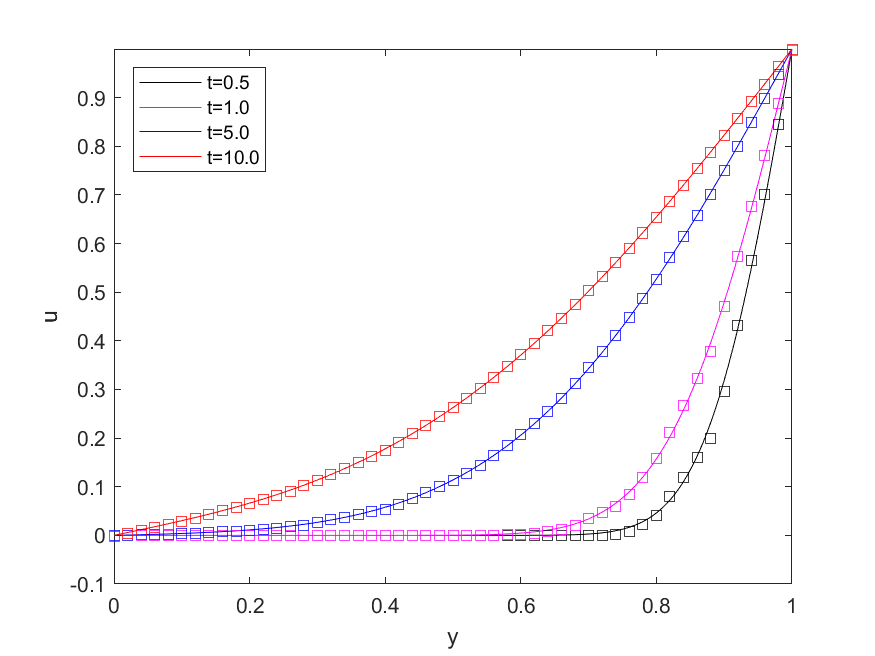}
		\caption{Unsteady Couette flow test: plots of the exact solution (continuous line) and the numerical solution (squares) at times $t=0.5$, $t=1.0$, $t=5.0$ and $t=10.0$.}
		\label{Test27LV_Graphs}
	\end{center}
\end{figure}

Now, let us study how the steady solution of the Couette flow is reached, starting from an initial condition of still fluid. Let us indicate the horizontal velocity of the fluid at height $y$ and time $t$ as $u(y,t)$. The problem is described by the parabolic equation:
\begin{equation}
\frac{\partial u}{\partial t} - \nu \frac{\partial ^2 u}{\partial y ^2} = 0 ,
\label{UnsteadyCouette}
\end{equation}
by the boundary conditions:
$$u(0,t)=0, u(L,t)=U, \quad \forall t \in \left[0,T\right],$$
and by the initial condition:
$$u(y,0)=0, \quad \forall y\in (0,L) .$$
The exact solution of this problem can be found, for example, by the method of separation of variables and is given by:
$$u(y,t) = U\frac{y}{L}-2\frac{U}{\pi}\sum_{n=1}^{\infty}\frac{1}{n} e^{-n^2\pi^2\nu\frac{t}{L^2}} \sin\left( n\pi\left(1-\frac{y}{L}\right) \right).$$
Also for this unsteady version of the Couette problem, we choose the following parameters: $\nu=0.01$, $L=1$,  $U=1$. In Figure \ref{Test27LV_Graphs} the solution obtained by our numerical method with $N=p_{\gamma}=1$ on the coarse unstructured mesh of Figure \ref{mesh116} is reported, together with the plot of the exact solution, for different times ($t=0.5$, $t=1.0$, $t=5.0$, $t=10.0$). The comparison beween the numerical and the analytical solution shows a very good agreement. We observe that at large times, the solution of this problem tends towards the steady solution $ u(y) = U\frac{y}{L} $ which was studied in the previous numerical test.
\\
In order to compare the two versions of our numerical method, we executed several simulations of the unsteady Couette flow test, up to the final time $T =10.0 $. The results for a set of increasing space and time polynomial degrees $N = p_{\gamma}$ are shown in Table \ref{TEST27LV_moving_vs_fixed}. We can conclude, also for this numerical test, that:
\begin{enumerate}
	\item the new ALE numerical method is consistent in the Eulerian limit, since the results obtained by the "ALE, $V=0$" and "Eulerian" implementations are very similar;
	\item the new algorithm is cheaper from the point of view of the computational time;
	\item the new algorithm is cheaper from the point of view of the computer memory usage.
\end{enumerate} 
\begin{table}[htbp] 
	\caption{Unsteady Couette flow test: errors in the $ L^2 $ norm for the velocity, CPU computational time (in seconds), and RAM's usage (in $MB$). In table (b), the Eulerian implementation was used.} 
	\begin{center}
		\scriptsize 
		\subfloat[ALE, $V=0$]{
			\label{TEST27LV_moving}
			\begin{tabular}{|l|c|c|c|} 
				\hline 
				$ N=p_{\gamma} $  	&$E_{2}^{v}$ 	&CPU time [$s$]  & RAM [$MB$] \\ \hline 
				1   & 9.93E-04 & 46 & 25	\\ \hline 
				2   & 2.67E-04 & 723 & 49	\\ \hline 
				3   & 1.07E-04 & 15579 & 163	\\ \hline 
			\end{tabular} 
		}
		\subfloat[Eulerian]{
			\label{TEST27LV_fixed}
			\begin{tabular}{|c|c|c|} 
				\hline 
				$E_{2}^{v}$ 	&CPU time [$s$]  & RAM [$MB$]  	\\ \hline 
				8.58E-04 & 50 & 30 \\ \hline 
				2.66E-04 & 2290 & 109	\\ \hline 
				1.07E-04 & 76302 & 487	\\ \hline 
			\end{tabular} 
		}
		\label{TEST27LV_moving_vs_fixed} 
	\end{center}
\end{table} 
\subsection{The Taylor-Green vortex}
\label{Tay_Gre}
\begin{figure}[htbp]
	\begin{center}
		\subfloat[$N=p_{\gamma}=0$]{
			\includegraphics[width=0.5\textwidth]{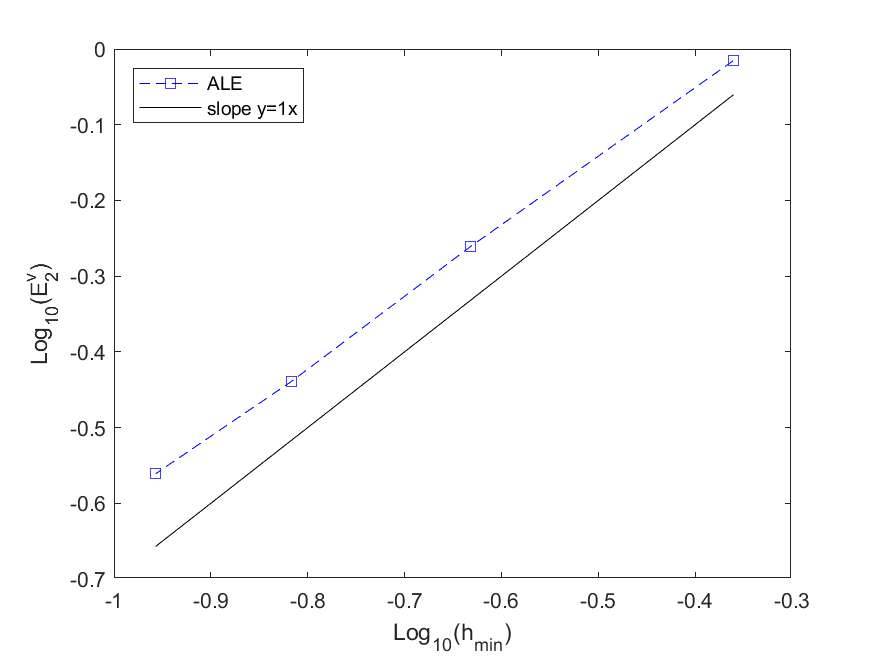}
		}
		\subfloat[$N=p_{\gamma}=1$]{
			\includegraphics[width=0.5\textwidth]{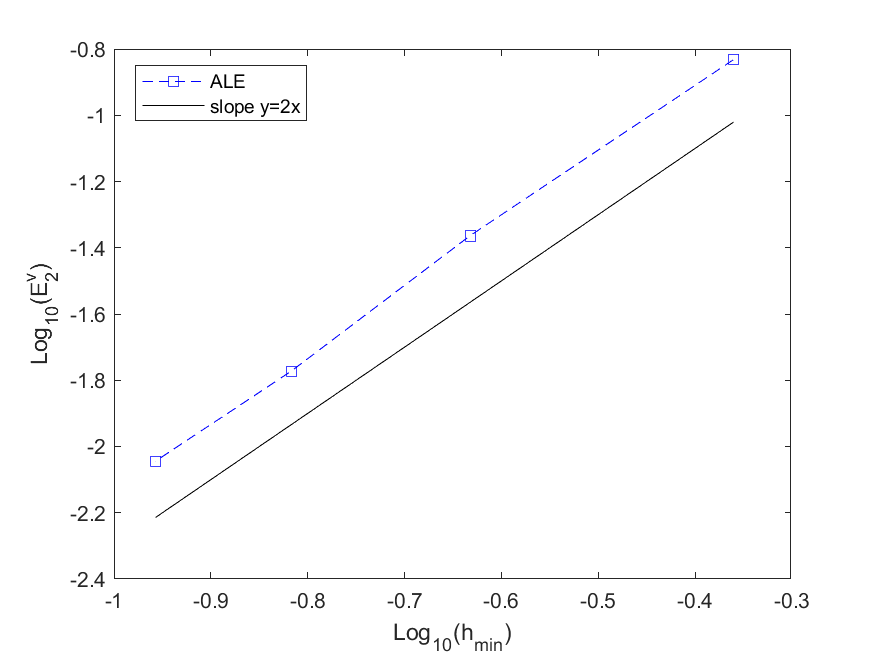}
		}
		\\
		\subfloat[$N=p_{\gamma}=2$]{
			\includegraphics[width=0.5\textwidth]{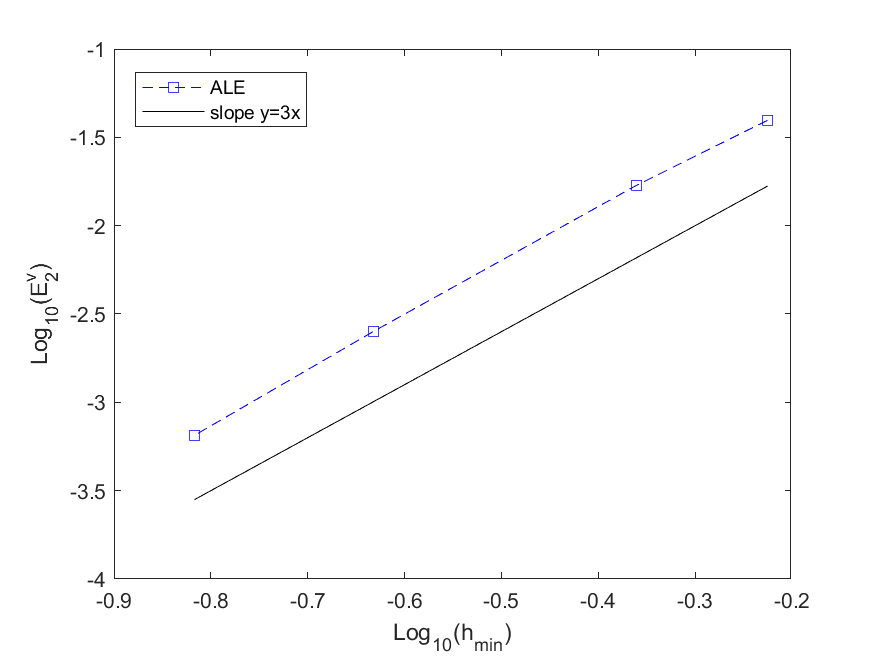}
		}
		\subfloat[$N=p_{\gamma}=3$]{
			\includegraphics[width=0.5\textwidth]{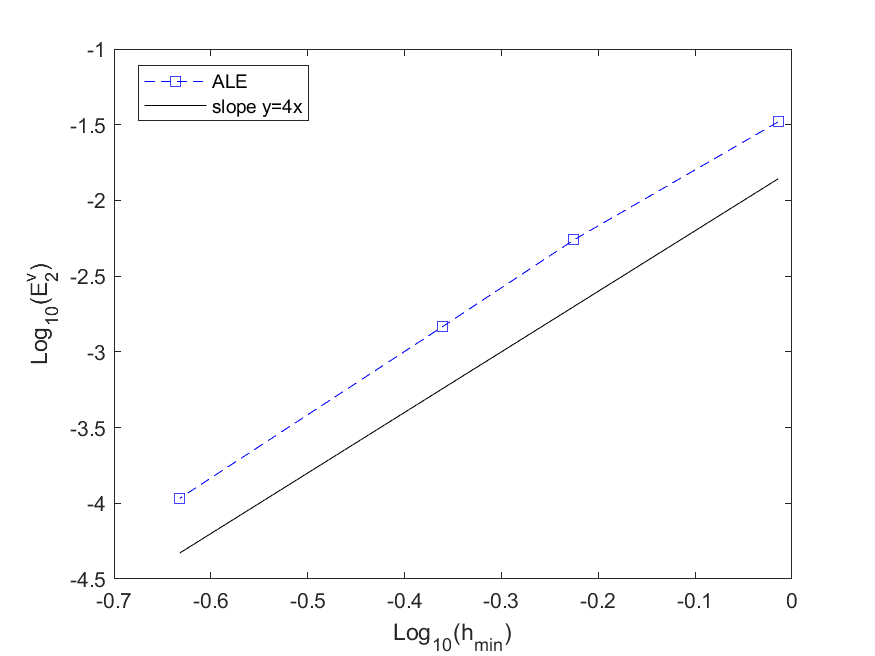}
		}
		\caption{Unsteady Taylor-Green vortex test: slope in the logarithmic plane of the $L^2$ norm of the error for the velocity field $(u,v)$, for different polynomial degrees.}
		\label{TG_Log_Graphs}
	\end{center}
\end{figure}

\begin{figure}[htbp]
	\begin{center}
		\subfloat[$N=p_{\gamma}=0$]{
			\includegraphics[height=0.33\textwidth]{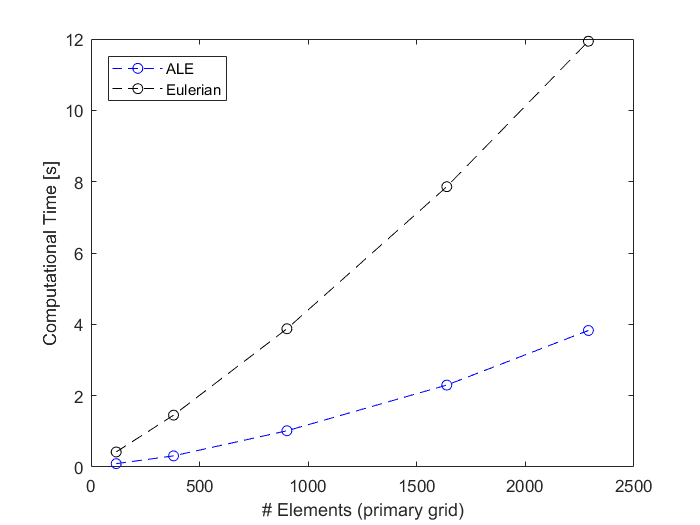}
		}
		\subfloat[$N=p_{\gamma}=1$]{
			\includegraphics[width=0.5\textwidth]{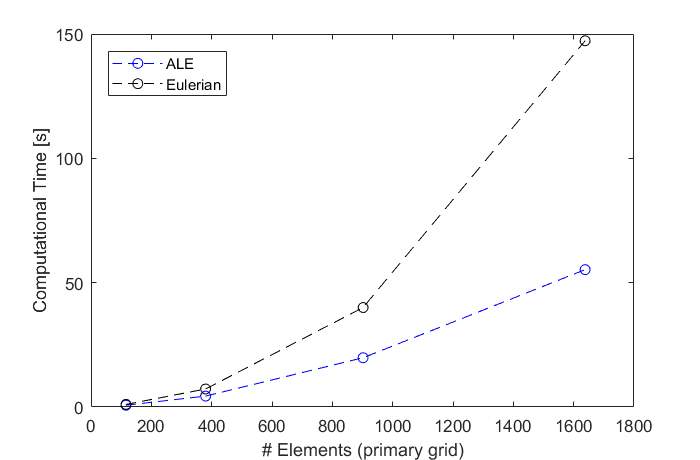}
		}
		\\
		\subfloat[$N=p_{\gamma}=2$]{
			\includegraphics[width=0.5\textwidth]{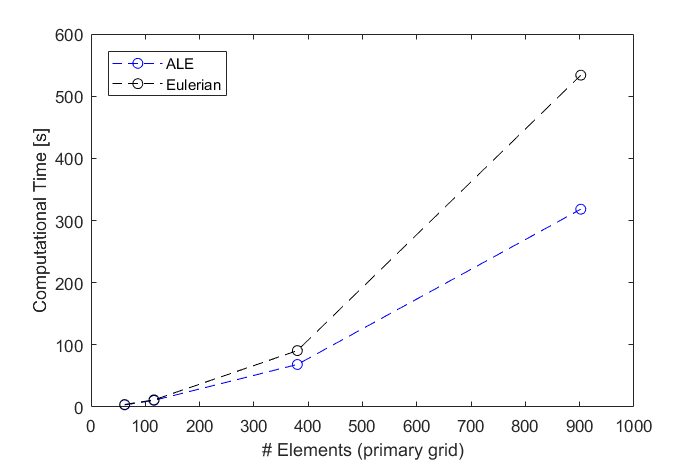}
		}
		\subfloat[$N=p_{\gamma}=3$]{
			\includegraphics[width=0.5\textwidth]{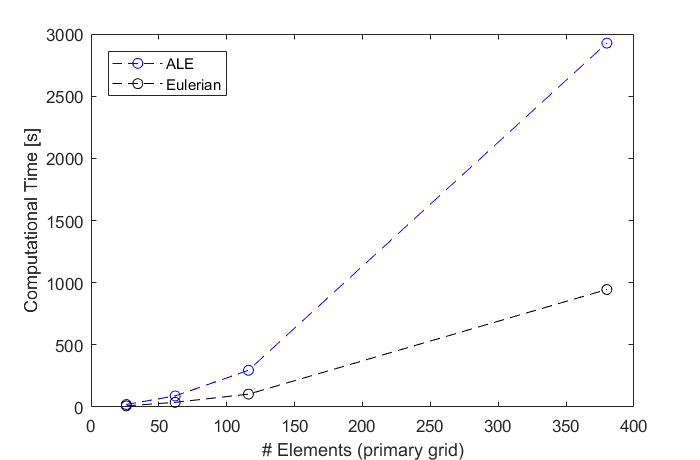}
		}
		\caption{Unsteady Taylor-Green vortex test: comparison between the computational times of the simulations for the ALE and the Eulerian implementations.}
		\label{TG_Comp_Times}
	\end{center}
\end{figure}

As next numerical test for the incompressible Navier-Stokes equations we consider the two-dimensional version of the Taylor-Green vortex problem \cite{TaylorGreen:1937}. We consider periodic boundary conditions for the square $ \Omega = [-\pi,\pi]^2 $ and a final time given by $T = 0.1$. Because of the presence of the viscous stresses, the vortices decay in time, losing energy. The analytical solution of the problem is given by:
\begin{eqnarray}
u(x,y,t) &=& \sin(x) \cos(y) e^{-2\nu t}; \\
v(x,y,t) &=& -\cos(x) \sin(y) e^{-2\nu t}; \\
p(x,y,t) &=& \frac{1}{4}\left( \cos(2x)+ \cos(2y)\right) e^{-4\nu t} .
\label{unst_TG}
\end{eqnarray}
The assessment of the order of convergence, for multiple simulations, for different polynomial degrees $N$, is reported in Table \ref{TEST19_P0P3} and in Figure \ref{TG_Log_Graphs}.
We observe that the errors obtained by the two implementations are very similar, which is a confirmation of the consistency of our numerical method in the Eulerian limit.
Moreover, regarding the estimated order of convergence, it is definitively close to the value $N+1$ for the $L^2$ norm of the velocity error, while for the pressure the value $N+1$ is clearly achieved only for the cases $N=0$ and $N=2$.\\
\begin{table}[htbp] 
	\caption{Unsteady Taylor-Green vortex test: errors in the $ L^2 $ norms for the pressure $p$ and the velocity field $(u,v)$, and convergence ratio, obtained by $N = p_{\gamma} =0,1,2,3$ on unstructured meshes with increasing number of triangles $ N_{e} $.}
	\begin{center}
		\scriptsize 
		\subfloat[$N=0$, ALE]{
			\begin{tabular}{|l|c|c|c|c|} 
				\hline 
				$ N_{e} $   &$E_{2}^{p}$ 	&$E_{2}^{v}$ 	&$\sigma_{2}^{p}$ 	&$\sigma_{2}^{v}$	 \\ \hline 
				116   & 8.66E-01 & 9.66E-01 &  -  &  -  \\ \hline 
				380   & 4.34E-01 & 5.48E-01 & 1.1 & 0.9 \\ \hline 
				902   & 2.38E-01 & 3.64E-01 & 1.4 & 1.0 \\ \hline 
				1638   & 1.95E-01 & 2.75E-01 & 0.6 & 0.9 \\ \hline 
				%2290   & 2.67E-01 & 2.34E-01 & -2.1 & 1.1 \\ \hline 
			\end{tabular} 
		}
		\subfloat[$N=0$, Eulerian]{
			\begin{tabular}{|l|c|c|c|c|} 
				\hline 
				$ N_{e} $   &$E_{2}^{p}$ 	&$E_{2}^{v}$ 	&$\sigma_{2}^{p}$ 	&$\sigma_{2}^{v}$	 \\ \hline 
				116   & 8.62E-01 & 9.65E-01 &  -  &  -  \\ \hline 
				380   & 4.29E-01 & 5.48E-01 & 1.1 & 0.9 \\ \hline 
				902   & 2.34E-01 & 3.64E-01 & 1.4 & 1.0 \\ \hline 
				1638   & 1.80E-01 & 2.74E-01 & 0.8 & 0.9 \\ \hline 
				%2290   & 3.78E-01 & 2.34E-01 & -5.1 & 1.1 \\ \hline 
			\end{tabular} 	
		}
		\\
		\subfloat[$N=1$, ALE]{
			\begin{tabular}{|l|c|c|c|c|} 
				\hline 
				$ N_{e} $   &$E_{2}^{p}$ 	&$E_{2}^{v}$ 	&$\sigma_{2}^{p}$ 	&$\sigma_{2}^{v}$	 \\ \hline 
				116   & 1.61E-01 & 1.47E-01 &  -  &  -  \\ \hline 
				380   & 5.09E-02 & 4.33E-02 & 1.8 & 2.0 \\ \hline 
				902   & 2.28E-02 & 1.69E-02 & 1.9 & 2.2 \\ \hline 
				1638   & 1.69E-02 & 9.03E-03 & 0.9 & 1.9 \\ \hline 
			\end{tabular} 
		}
		\subfloat[$N=1$, Eulerian]{
			\begin{tabular}{|l|c|c|c|c|} 
				\hline 
				$ N_{e} $   &$E_{2}^{p}$ 	&$E_{2}^{v}$ 	&$\sigma_{2}^{p}$ 	&$\sigma_{2}^{v}$	 \\ \hline 
				116   & 1.35E-01 & 1.04E-01 &  -  &  -  \\ \hline 
				380   & 4.48E-02 & 3.11E-02 & 1.8 & 1.9 \\ \hline 
				902   & 2.09E-02 & 1.23E-02 & 1.8 & 2.2 \\ \hline 
				1638   & 1.18E-02 & 6.65E-03 & 1.8 & 1.9 \\ \hline 
			\end{tabular} 
		}
		\\
		\subfloat[$N=2$, ALE]{
			\begin{tabular}{|l|c|c|c|c|} 
				\hline 
				$ N_{e} $   &$E_{2}^{p}$ 	&$E_{2}^{v}$ 	&$\sigma_{2}^{p}$ 	&$\sigma_{2}^{v}$	 \\ \hline 
				62   & 6.26E-02 & 3.95E-02 &  -  &  -  \\ \hline 
				116   & 2.67E-02 & 1.69E-02 & 2.7 & 2.7 \\ \hline 
				380   & 4.49E-03 & 2.52E-03 & 2.9 & 3.0 \\ \hline 
				902   & 1.36E-03 & 6.50E-04 & 2.8 & 3.2 \\ \hline 
			\end{tabular} 
		}
		\subfloat[$N=2$, Eulerian]{
			\begin{tabular}{|l|c|c|c|c|} 
				\hline 
				$ N_{e} $   &$E_{2}^{p}$ 	&$E_{2}^{v}$ 	&$\sigma_{2}^{p}$ 	&$\sigma_{2}^{v}$	 \\ \hline 
				62   & 4.40E-02 & 2.18E-02 &  -  &  -  \\ \hline 
				116   & 1.87E-02 & 8.43E-03 & 2.8 & 3.1 \\ \hline 
				380   & 3.57E-03 & 1.17E-03 & 2.6 & 3.2 \\ \hline 
				902   & 1.11E-03 & 2.97E-04 & 2.7 & 3.2 \\ \hline 
			\end{tabular} 
		}
		\\
		\subfloat[$N=3$, ALE]{
			\begin{tabular}{|l|c|c|c|c|} 
				\hline 
				$ N_{e} $   &$E_{2}^{p}$ 	&$E_{2}^{v}$ 	&$\sigma_{2}^{p}$ 	&$\sigma_{2}^{v}$	 \\ \hline 
				26   & 5.43E-02 & 3.30E-02 &  -  &  -  \\ \hline 
				62   & 8.05E-03 & 5.48E-03 & 3.9 & 3.7 \\ \hline 
				116   & 4.56E-03 & 1.47E-03 & 1.8 & 4.2 \\ \hline 
				380   & 2.04E-03 & 1.08E-04 & 1.3 & 4.2 \\ \hline 
			\end{tabular} 
		}
		\subfloat[$N=3$, Eulerian]{
			\begin{tabular}{|l|c|c|c|c|} 
				\hline 
				$ N_{e} $   &$E_{2}^{p}$ 	&$E_{2}^{v}$ 	&$\sigma_{2}^{p}$ 	&$\sigma_{2}^{v}$	 \\ \hline 
				26   & 3.63E-02 & 1.54E-02 &  -  &  -  \\ \hline 
				62   & 6.54E-03 & 2.77E-03 & 3.5 & 3.5 \\ \hline 
				116   & 4.64E-03 & 7.00E-04 & 1.1 & 4.4 \\ \hline 
				380   & 2.06E-03 & 4.61E-05 & 1.3 & 4.3 \\ \hline 
			\end{tabular} 
		}
		\label{TEST19_P0P3} 
	\end{center}
\end{table}
For this numerical test, the analysis of the computational times is reported in Figure \ref{TG_Comp_Times}. The new ALE algorithm requires the recomputation of all the matrices listed in Formulae \eqref{eq:MD_3}-\eqref{eq:MD_5_2}, at every time-step of the simulation. For this reason, we expect it to be slower than the Eulerian one. Instead, for the cases $N=0,1$ and $2$, this tendency is inverted, independently from the mesh size. We can state that this is a consequence of two facts. \\
Firstly, in our new geometry-free implementation, the pre-computational time is reduced to almost zero, because the algorithm only needs to read an input binary file for the universal tensors which have already been computed analytically on the reference space--time elements. On the contrary, the pre-computational time of the Eulerian algorithm is dependent on the employed mesh, and it is larger for larger mesh sizes.\\
Secondly, the Lagrangian recalculation of the matrices is written in such a way to be very efficient from the computational point of view: the resulting matrices in \eqref{eq:MD_3}-\eqref{eq:MD_5_2} come from an operation of tensor-reduction between the precomputed universal tensors and the updated coefficients containing the geometric information of the mesh (see Remark \ref{rem_ALE}). From the results of the computational times, it is evident that, up to the case $N=p_{\gamma}=2$, the first data of this multiplication are so small that they can be stored in the first levels of the memory cache, allowing a very fast recalculation of the matrices. But, for the last case, $N=p_{\gamma}=3$, the access to these larger data is not very fast anymore, because their dimension exceeds the dimension of the cache, therefore the simulations turn out to be slower.
\\
Afterwards, an analysis of the required memory can be found in Table \ref{TG_RAM_memory}. The advantage of using the new approach becomes more important, the larger is the mesh size. As it was expected, the required memory increases when increasing the parameters $N$, $p_{\gamma}$ and $N_{e}$, also for the ALE implementation. But an increase of $N_{e}$ only affects the dimension of the variables to be stored (e.g. the array for the nodes of the mesh, the arrays for the space--time degrees of freedom of the solution, \ldots) and not the universal tensors, whose size only depends on $N$ and $p_{\gamma}$. Viceversa, in the Eulerian implementation, we must store all the local matrices \eqref{eq:MD_3}-\eqref{eq:MD_5_2} which occupy a very large part of the memory and which are directly proportional to the mesh size $N_{e}$. For this reason, the Eulerian implementation required a bigger amount of memory, for every numerical test that we have performed. In the tests that we have performed with the coarsest meshes, the ALE implementation required about the half of the memory of the Eulerian one.
\begin{table}[htbp]
	\caption{Unsteady Taylor-Green vortex test: memory usage (in $MB$), required by the ALE and the Eulerian implementations, and percentage of the first with respect to the second, for different polynomial degrees.} 
	\begin{center}
		\scriptsize
		\subfloat[$N=p_{\gamma}=1$]{
			\begin{tabular}{|l|c|c|c|} 
				\hline 
				$ N_{e} $   	& ALE, $V=0$ 	& Eulerian & Percentage	 \\ \hline 
				116   & 25  &  31 & 80.6 \% \\ \hline 
				380   & 39  & 52 & 75.0 \% \\ \hline 
				902   & 65  & 96 & 67.7 \% \\ \hline 
				1638  & 102  & 168 & 60.7 \% \\ \hline
			\end{tabular} 
		}
		\subfloat[$N=p_{\gamma}=2$]{
			\begin{tabular}{|l|c|c|c|} 
				\hline 
				$ N_{e} $   	& ALE, $V=0$ 	& Eulerian & Percentage	 \\ \hline 
				62   & 49  &  67 & 73.1 \% \\ \hline 
				116   & 67  & 106 & 63.2 \% \\ \hline 
				380   & 179  & 311 & 57.6 \% \\ \hline 
				902  & 392  & 709 & 55.3 \% \\ \hline
			\end{tabular} 
		}
		\\
		\subfloat[$N=p_{\gamma}=3$]{
			\begin{tabular}{|l|c|c|c|} 
				\hline 
				$ N_{e} $   	& ALE, $V=0$ 	& Eulerian & Percentage	 \\ \hline 
				26   & 110  &  120 & 91.7 \% \\ \hline 
				62   & 181  & 262 & 69.1 \% \\ \hline 
				116   & 279  & 472 & 59.1 \% \\ \hline 
				380  & 800  & 1497 & 53.4 \% \\ \hline
			\end{tabular} 
		}
		\label{TG_RAM_memory} 
	\end{center}
\end{table} 
\subsection{The lid-driven cavity}
\label{Cavity}
\begin{figure}[htbp]
	\begin{center}
		\subfloat[$Re=100$]{
			\includegraphics[width=0.5\textwidth]{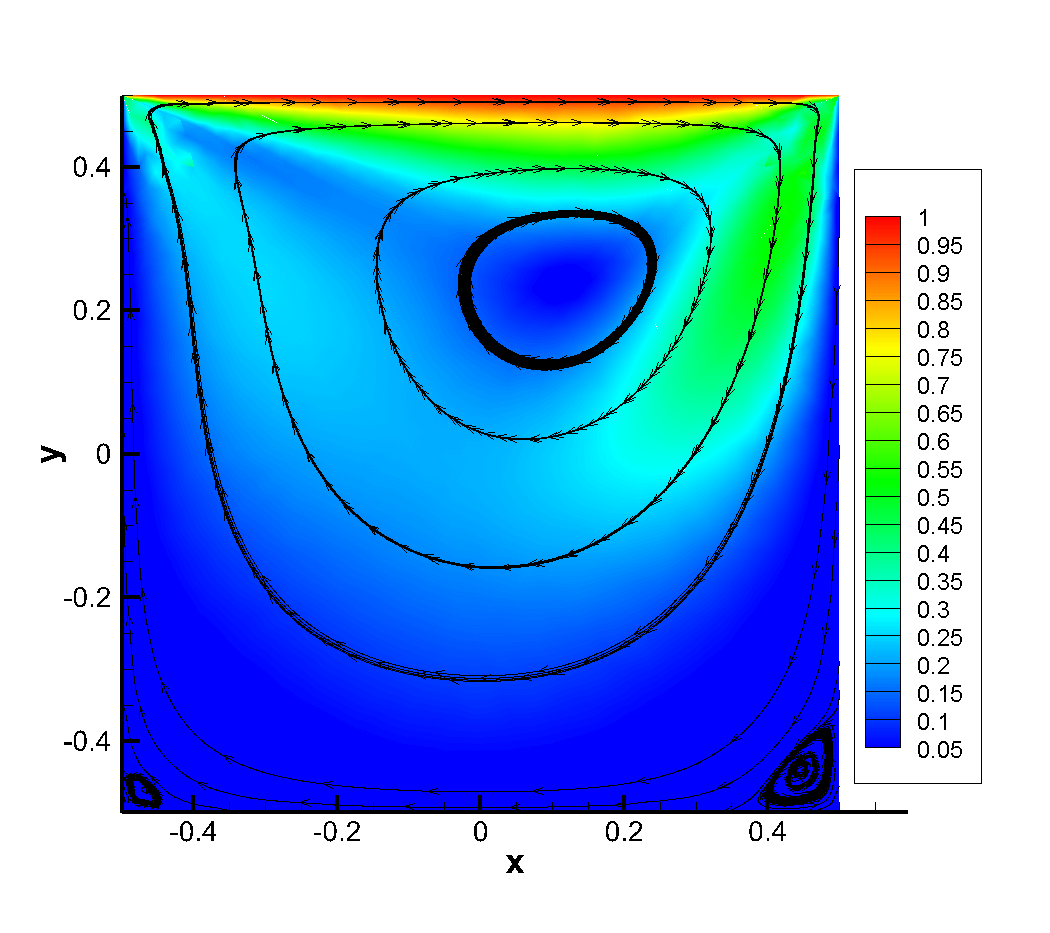}
			\label{LDC_100}
		}
		\subfloat[$Re=400$]{
			\includegraphics[width=0.5\textwidth]{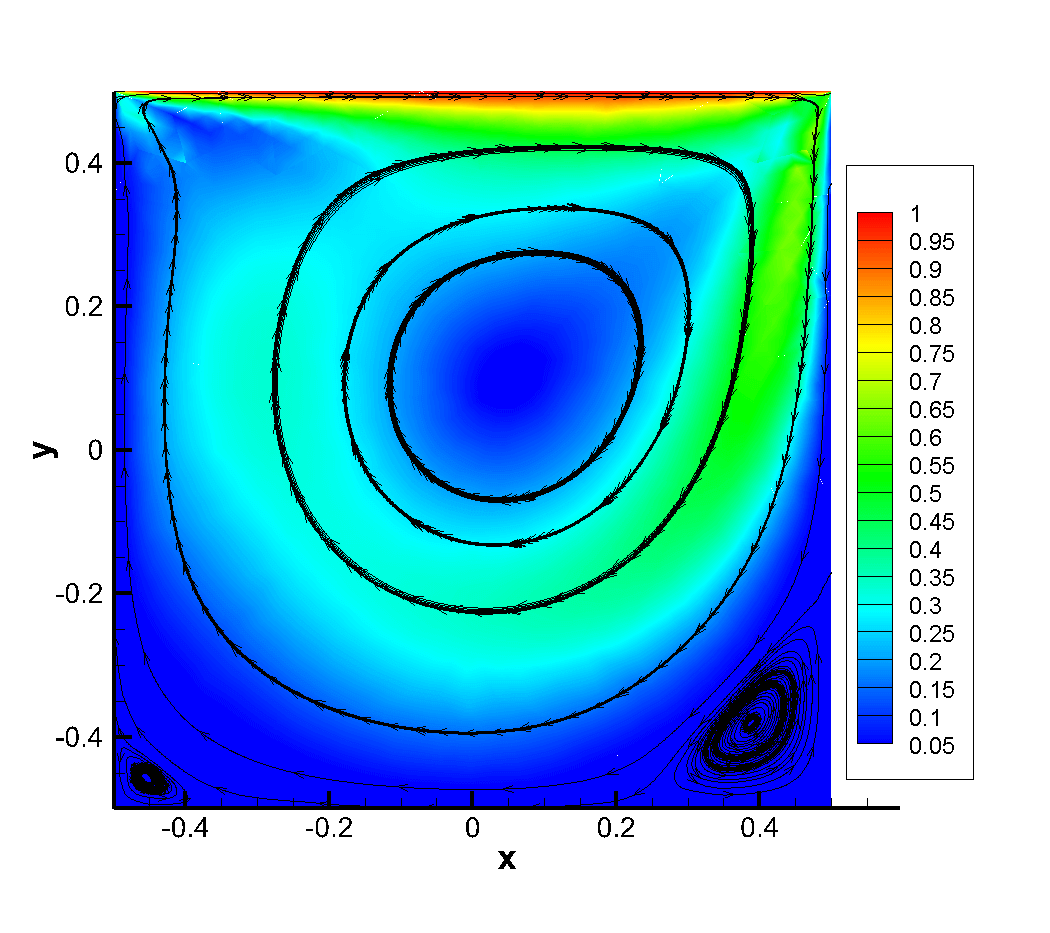}
			\label{LDC_400}
		}
		\\
		\subfloat[$Re=1000$]{
			\includegraphics[width=0.5\textwidth]{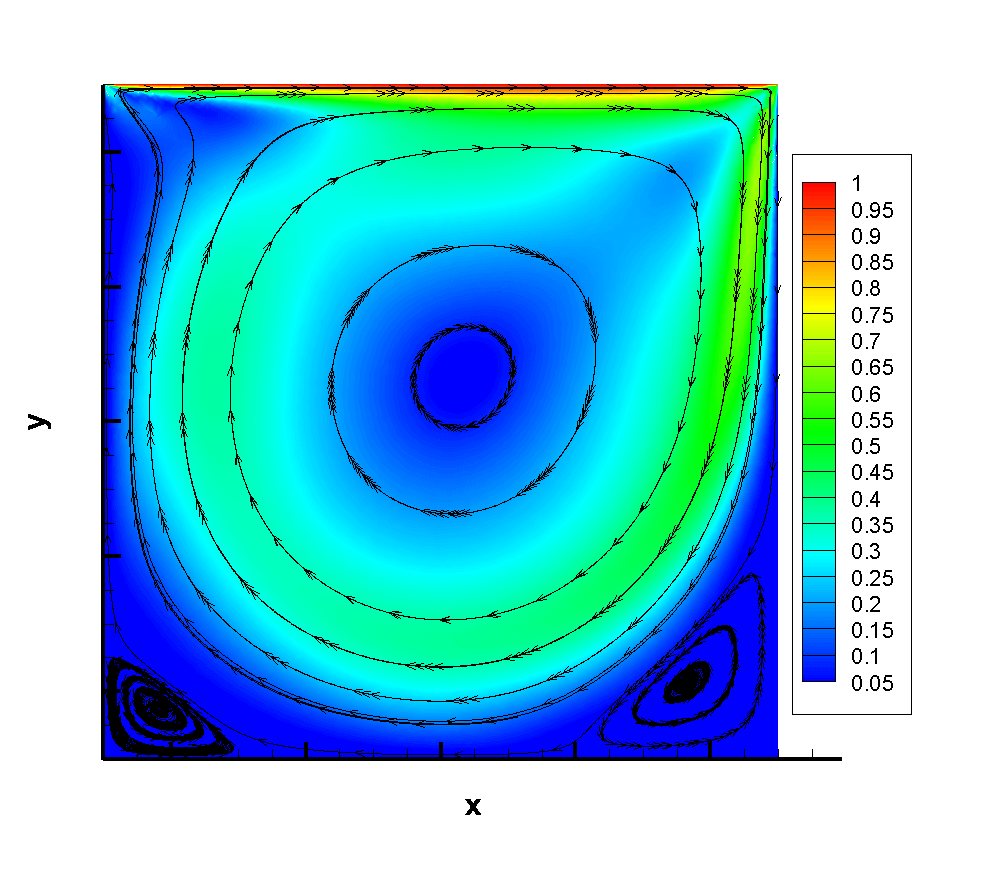}
			\label{LDC_1000}
		}
		\caption{Lid-driven cavity test: plots of the velocity's magnitude and streamlines, at the final simulation time $T$, for different values of the Reynolds number.}
		\label{LDC_u_v_T}
	\end{center}
\end{figure}

\begin{figure}[htbp]
	\begin{center}
		\subfloat[$u(0.5,y,T)$, $Re=100$]{
			\includegraphics[width=0.5\textwidth]{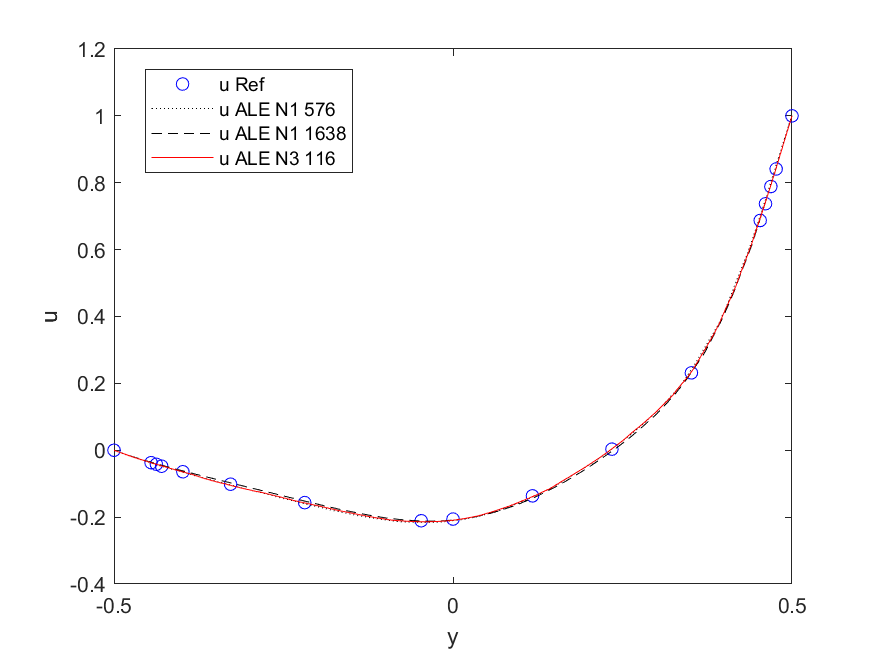}
		}
		\subfloat[$v(x,0.5,T)$, $Re=100$]{
			\includegraphics[width=0.5\textwidth]{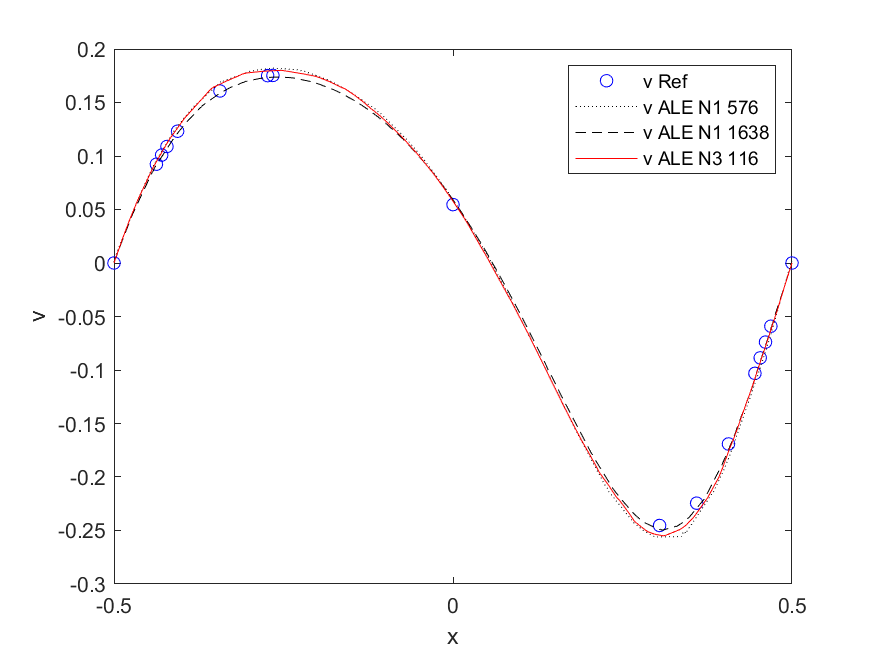}
		}
		\\
		\subfloat[$u(0.5,y,T)$, $Re=400$]{
			\includegraphics[width=0.5\textwidth]{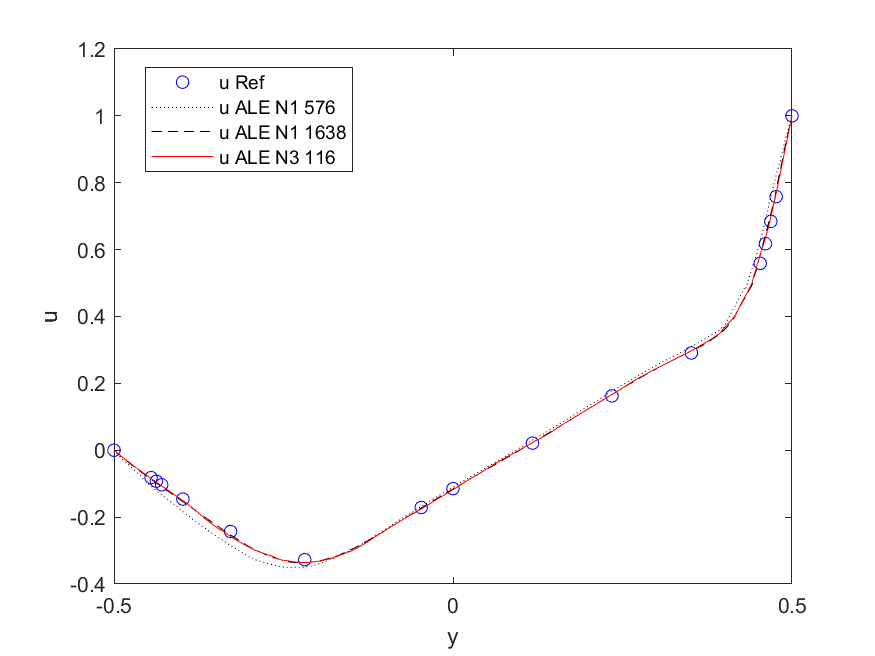}
		}
		\subfloat[$v(x,0.5,T)$, $Re=400$]{
			\includegraphics[width=0.5\textwidth]{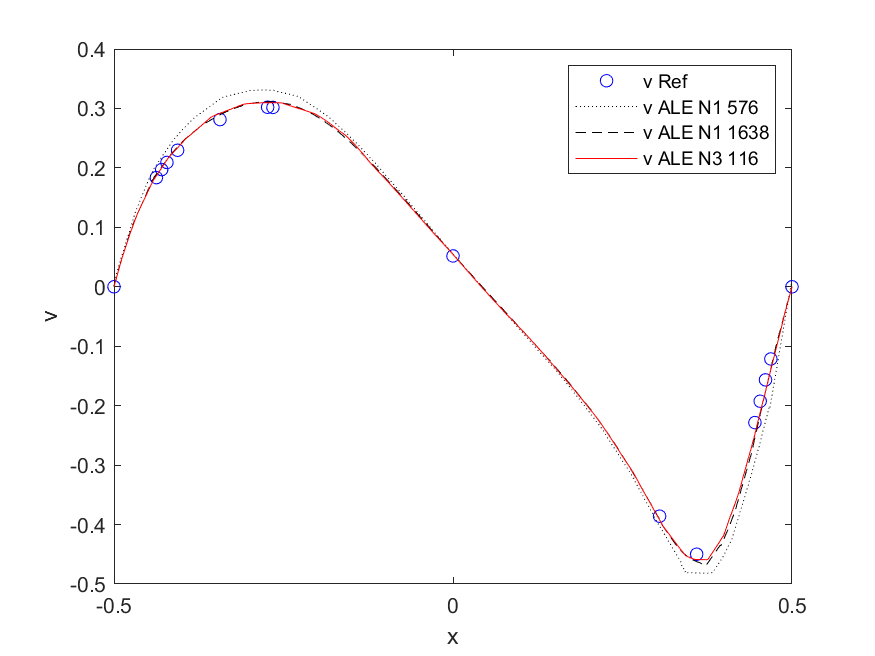}
		}
		\\
		\subfloat[$u(0.5,y,T)$, $Re=1000$]{
			\includegraphics[width=0.5\textwidth]{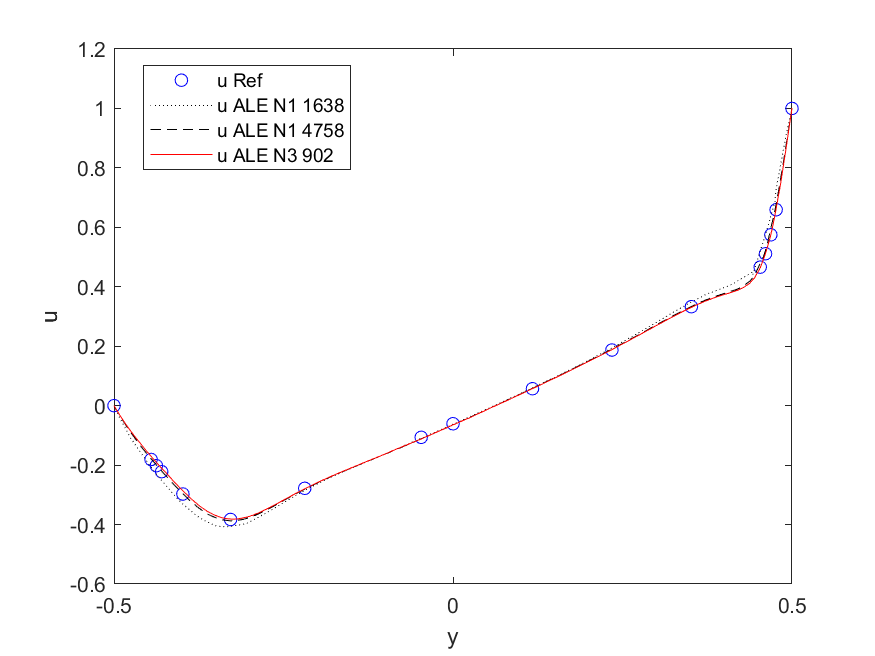}
		}
		\subfloat[$v(x,0.5,T)$, $Re=1000$]{
			\includegraphics[width=0.5\textwidth]{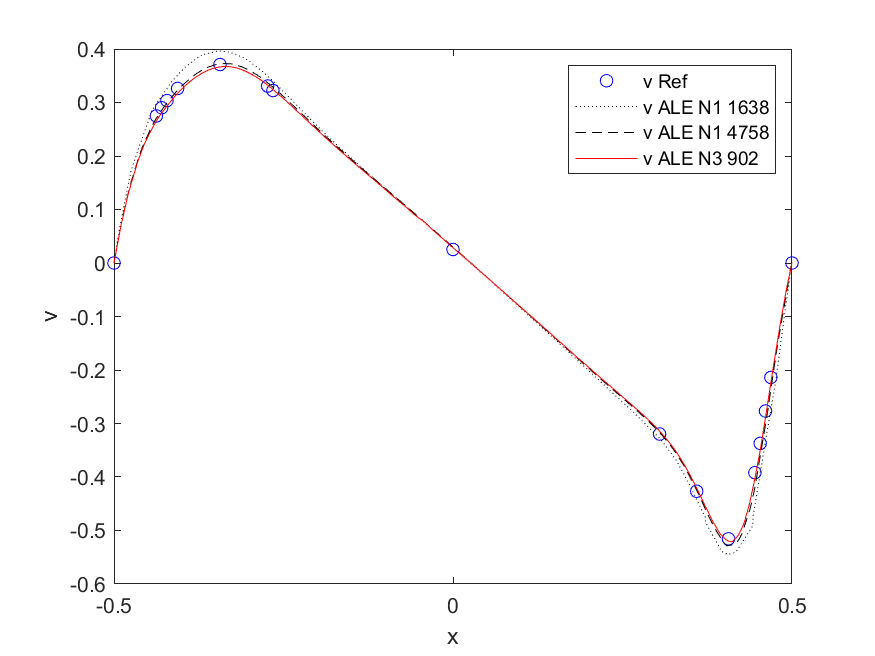}
		}
		\caption{Lid-driven cavity test: plots of the velocity components in the middle lines of the domain, for different values of the Reynolds number.}
		\label{LDC_u_v_middle_lines}
	\end{center}
\end{figure}

The lid-driven cavity test, which was proposed for the first time in \cite{Ghia1982}, has become a standard benchmark for the verification of the numerical solvers of the incompressible Navier-Stokes equations. Indeed, despite its geometrical simplicity, this test has got a complex flow regime which is characterized by the formation of some interesting vortices at the corners of the cavity.\\
Let us consider the square domain $\Omega=[-0.5,0.5] ^2$ and let us define the following boundary conditions (see Figure \ref{Cavity_description}): 
\begin{figure}[htbp]
	\begin{center}
		\includegraphics[width=0.45\textwidth]{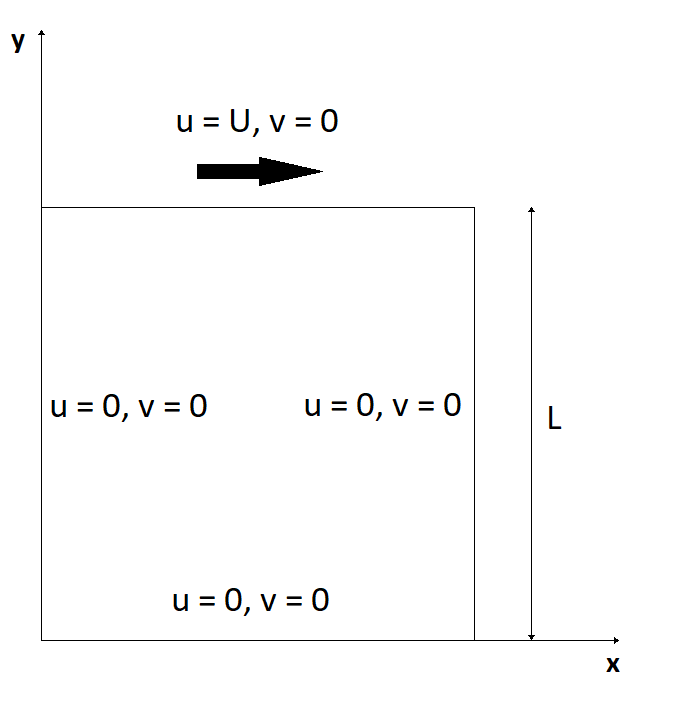}
		\caption{The lid-driven cavity flow test: a schematic representation.}
		\label{Cavity_description}
	\end{center}
\end{figure}
\begin{itemize}
	\item no-slip ("wall"-type) conditions at the bottom, at the left and at the right of the domain: $u=0$, $v=0$ and $\frac{\partial p}{\partial \vec{n}} = 0$ (zero normal derivative of the pressure);
	\item Dirichlet's boundary condition for the velocity at the top of the domain: $u=1$, $v=0$ and  $\frac{\partial p}{\partial \vec{n}} = 0$ (zero normal derivative of the pressure).
\end{itemize}
Let us consider an initial condition of still fluid inside the domain ($u=0$, $v=0$). The final simulation time $T$ is determined by the following tolerance criterion: 
\begin{equation}
\max_{j=1,\ldots,N_d}{|\hat{\mathbf{v}}_j^{n+1} - \hat{\mathbf{v}}_j^{n} |} < \epsilon .
\label{stop_crit}
\end{equation}
The simulation stops when the $L^{\infty}$ norm of the difference of the velocity fields from two consecutive solutions is below a tolerance value, $\epsilon=10^{-6}$. In this way, we can state that the steady solution of the problem has been reached.
\\
For the case $Re=100$, we used a coarse unstructured mesh with $N_{e}=116$ and we got the results depicted in Figure \ref{LDC_100} at around $T=8$, by using $N=3, p_{\gamma}=0$. Then, for the case $Re=400$, we used exactly the same discretization parameters and we obtained the results in Figure \ref{LDC_400} at around $T=20$. Finally, we used $N=3,p_{\gamma}=0$ and an unstructured mesh with $N_{e}=902$ for the case $Re=1000$, and we obtained the results shown in Figure \ref{LDC_1000}, at around $T=25$. From these plots, it is evident that not only the primary, central vortex, but also the two smaller, counter-rotating vortices in the two lower corners of the cavity (BR$_{1}$ and BL$_{1}$, using the notation in \cite{Ghia1982}) are well reproduced by our numerical method. We observe that, thanks to the use of the high-order in space $N=3$, the numerical solver is able to handle the singularities of the pressure and of the velocity field at the two upper corners, even if a quasi-uniform mesh with a modest number of elements is employed.\\
In order to have a quantitative assessment of the solution, we evaluate the velocity components at the two middle lines of the cavity $x=0$ and $y=0$, and we compare them with the reference results provided in \cite{Ghia1982}. The red lines in all the graphs of Figure \ref{LDC_u_v_middle_lines} refer to the high-order $N=3$ simulations that have been described above, showing a very good agreement with the reference data, both for the $u$ and the $v$ functions, for every case of the employed Reynolds number. In the same graphs, also the results given by the low-order method $N=1$ obtained for different meshes (the mesh size $N_{e}$ is reported in the legend of the graphs) are shown. We point out that the high-order $N=3$ scheme is very accurate because the differences with the low-order scheme applied to very refined meshes are almost indistinguishable. However, the benefits of using the high-order scheme can be clearly appreciated from the comparison of the performances with the low-order scheme (see Table \ref{LDC_Performances}).
\begin{table}[htbp] 
	\caption{Lid-driven cavity test: comparison of the performances between the high-order ($N=3$) and low-order ($N=1$) schemes, in the cases $Re=400$ (a) and $Re=1000$ (b).} 
	\begin{center}
		\scriptsize 
		\subfloat[$Re=400$]{
			\begin{tabular}{|l|c|c|c|} 
				\hline 
				& ALE N$3$ $116$ 	& ALE N$1$ $1638$  & Percentage \\ \hline 
				$N_{e}$			& 116	& 1638	& 7.1\% \\ \hline 
				CPU time [$min$]  & 21 & 111 & 18.9\%	\\ \hline 
				RAM [$MB$]   	& 33 & 44 & 75.0\%	\\ \hline 
			\end{tabular} 
		}
		\\
		\subfloat[$Re=1000$]{
			\begin{tabular}{|l|c|c|c|} 
				\hline 
				& ALE N$3$ $902$ 	& ALE N$1$ $4758$  & Percentage \\ \hline 
				$N_{e}$			& 902	& 4758	& 19.0\% \\ \hline 
				CPU time [$min$]  & 492 & 706 & 69.7\%	\\ \hline 
				RAM [$MB$]   	& 86 & 112 & 76.8\%	\\ \hline 
			\end{tabular} 
		}
		\label{LDC_Performances} 
	\end{center}
\end{table} 
\subsection{The Blasius boundary layer}
\label{Blasius}
\begin{figure}[htbp]
	\begin{center}
		\includegraphics[width=0.75\textwidth]{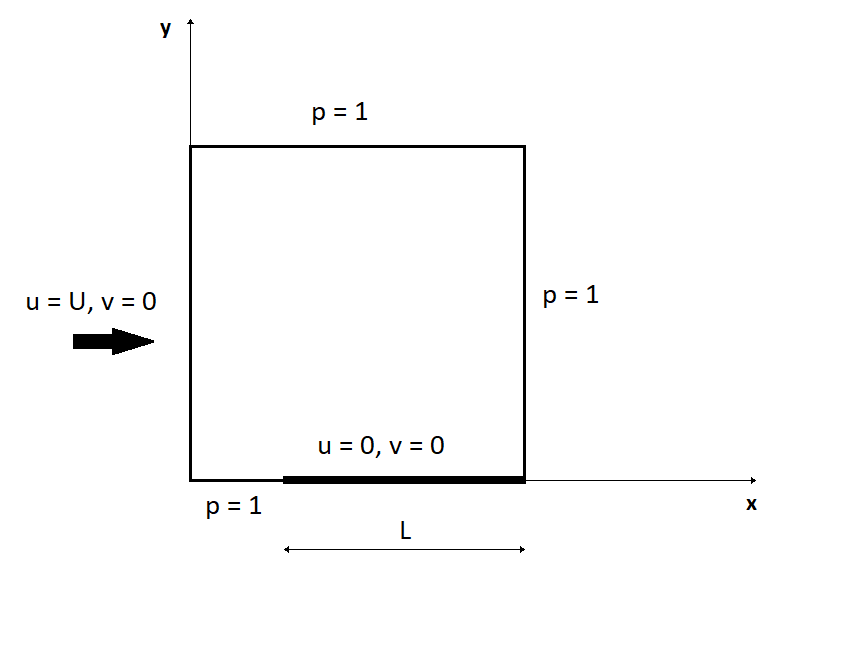}
		\caption{Blasius boundary layer: a scheme of the boundary conditions.}
		\label{Blasius_description}
	\end{center}
\end{figure}

\begin{figure}[htbp]
	\begin{center}
		\subfloat[$x=0.4$]{
			\includegraphics[width=0.75\textwidth]{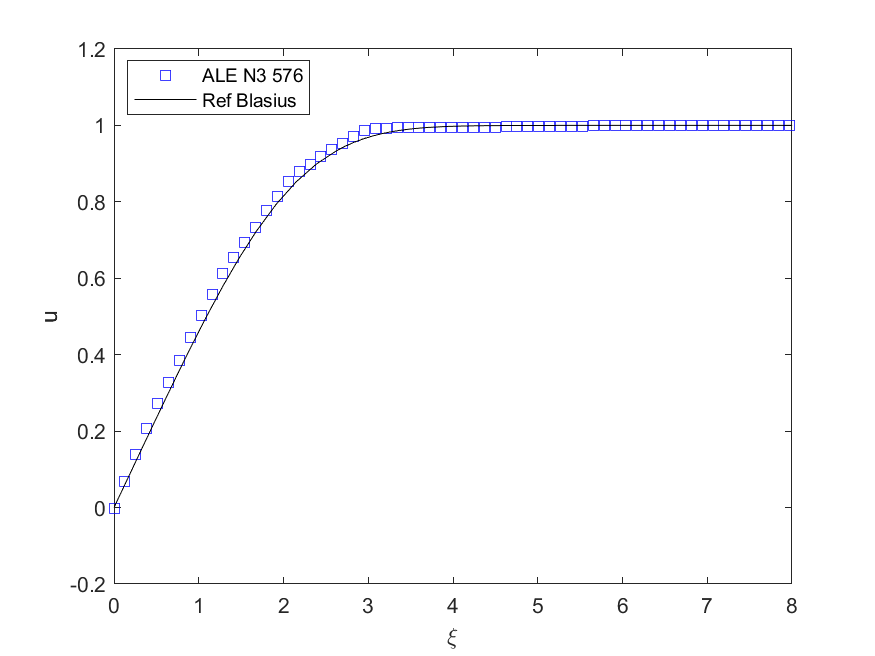}
		}
		\\
		\subfloat[$x=0.7$]{
			\includegraphics[width=0.75\textwidth]{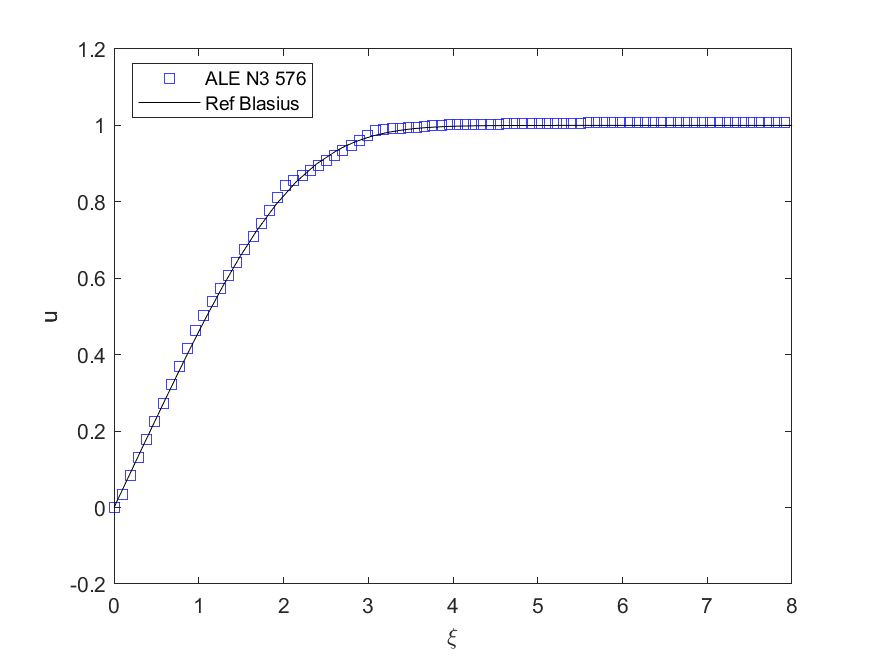}	
		}
		\caption{Blasius boundary layer test: comparison between the reference solution (continuous line) and the numerical solution (squares), at distances (a) $x=0.4$ and (b) $x=0.7$ from the attachment point of the boundary layer.}
		\label{Blasius_corner_3200_P3T0_36}
	\end{center}
\end{figure}

Another classical benchmark for the incompressible Navier-Stokes equations is the formation of the laminar boundary layer over a flat plate. Let us consider a steady laminar flow with velocity field $\vec{v}=(U,0)$ which passes over a flat plate with length $L$. Under the assumptions that there are no gravity forces and the pressure is constant, a self-similar solution of the Navier-Stokes equations was found by Blasius in \cite{Blasius1908}. The Blasius equations (third-order non-linear differential equations) are given by:
\begin{eqnarray}
\left\{
\begin{array}{l}
f'''+ff''=0, \\
f(0)=0, \quad f'(0)=0, \quad \lim_{\xi \rightarrow \infty} f'(\xi)=1,
\end{array}
\right.
\end{eqnarray}
where $\xi=y \sqrt{\frac{U}{2\nu x}}$ is the Blasius coordinate and $f'=\frac{u}{U}$.\\
Let us consider a square domain and let us use the following boundary conditions (see Figure \ref{Blasius_description}): 
\begin{itemize}
	\item ("inflow"-type) Dirichlet's boundary condition for the velocity on the left of the domain: $u=U$, $v=0$ and  $\frac{\partial p}{\partial \vec{n}} = 0$ (zero normal derivative of the pressure);
	\item no-slip ("wall"-type) conditions on the flat plate: $u=0$, $v=0$ and $\frac{\partial p}{\partial \vec{n}} = 0$ (zero normal derivative of the pressure);
	\item ("outflow"-type) Dirichlet's boundary condition for the pressure on the top and on the right of the domain and before the flat plate: $p=1$ and $\frac{\partial \vec{v}}{\partial \vec{n}} = 0$ (zero normal derivative of both the velocity components).
\end{itemize}
As initial conditions, we impose $p=1$ and the farfield velocity field $\vec{v}=(U,0)$ inside the domain, with $U=1$. We use an unstructured mesh with $N_{e}=576$ elements and we employ a polynomial degree $N=3$ in space and $p_{\gamma}=0$ in time. We set the kinematic viscosity coefficient to $\nu=3.125\cdot 10^{-4}$ and we study the resulting boundary layer at the distances $L=0.4$ and $L=0.7$ of the flat plate, so that the local Reynolds number $Re_{L} = \frac{U L}{\nu}$ is equal to $Re_{L}=1280$ and $Re_{L}=2240$, respectively. The simulation stopped at time $T=2.7$, when a steady state was reached, according to the criterion in Formula \eqref{stop_crit}, with tolerance $\epsilon=10^{-6}$.
In Figure \ref{Blasius_corner_3200_P3T0_36}, the results of the numerical solution, in the $\xi - u$ plot, are compared with the reference Blasius solution obtained in \cite{DumbserNSE} by a tenth-order DG Ordinary Differential Equations (ODE) solver with a classical shooting method for the boundary conditions. A very good agreement can be observed, despite the use of a quite coarse mesh. Moreover, the self-similarity property of the solution is verified. %: the numerical solutions at different distances from the attachment point of the plate have got the same shape.

\section{Conclusions}
\label{concl}
In this article, we have presented a new high order accurate staggered semi-implicit space--time discontinuous Galerkin method for the two dimensional incompressible Navier-Stokes equations. The numerical scheme is based on the work proposed in \cite{STINS2D,STINS2DTri,STCNS}, but is oriented towards the Arbitrary Lagragian Eulerian framework thanks to the use of a suitable pair of velocity-pressure finite elements. The same polynomial degree $N$ is employed for representing both the pressure and the velocity, therefore it is named an  equal-order-interpolation finite element method. \\
The new scheme has been verified through several benchmarks in fluid dynamics, from the Couette flow to the cavity flow and the Blasius boundary layer. The numerical results agree very well to the reference solutions in all the numerical tests that we have performed. The high order of accuracy has been verified in the Taylor-Green vortex test, up to $N=3$ and $p_{\gamma}=3$.\\
The consistency in the Eulerian limit has been assessed by comparing the results obtained with the staggered space--time DG algorithm proposed in \cite{STINS2DTri}.
Compared to the Eulerian implementation, the new ALE implementation proposed here is not only computationally more efficient, but also less memory consuming, since all universal matrices are precomputed before runtime, once and for all, on a reference element. 
Indeed, all local matrices in the Galerkin formulation can be updated, with their geometric information, at every time-step $n$, by a cheap tensor-reduction operation.

\begin{acknowledgements}
This work was partially funded by the research project \textsl{STiMulUs}, ERC Grant agreement no. 278267. Financial support has also been provided by the Italian Ministry of Education, University and Research (MIUR) in the frame of the \textit{Departments of Excellence Initiative 2018-2022} attributed to DICAM of the University of Trento (grant L. 232/2016) and via the PRIN2017 project. The authors have also received funding from the University of Trento via the Strategic Initiative \textit{Modeling and Simulation}. 
\end{acknowledgements}

\section*{Conflict of interest}
On behalf of all authors, the corresponding author states that there is no conflict of interest. 

\bibliographystyle{plain}
\bibliography{short,romeo_new}   % name your BibTeX data base

\end{document}